\documentclass[11pt]{article}
\usepackage[english]{babel}
\usepackage[left=1in,right=1in,top=1in,bottom=1.5in]{geometry} 
\usepackage[skip=5pt, indent=10pt]{parskip}
\usepackage{amssymb,amsmath, amsthm}
\usepackage{hyperref}
\usepackage{xcolor}
\usepackage{nicefrac}
\usepackage[margin=1cm]{caption}
\usepackage{natbib}
\usepackage[makeroom]{cancel}
\usepackage{bbm}
\usepackage{mathrsfs}
\usepackage{algorithm}
\usepackage{algpseudocode}
\usepackage{times}
\usepackage{graphicx}

\definecolor{darkred}{HTML}{880000}
\definecolor{darkblue}{HTML}{000088}

\hypersetup{
  colorlinks,
  linkcolor={darkred},
  citecolor={darkblue}
}

\renewcommand{\le}{\leqslant}
\renewcommand{\leq}{\leqslant}

\renewcommand{\geq}{\geqslant}

\newcommand{\eps}{\varepsilon}

\newcommand{\rmd}{{\rm d}}
\newcommand{\Fr}{{\rm F}}
\newcommand{\diag}{{\rm diag}}
\newcommand{\Tr}{{\rm Tr}}
\newcommand{\Var}{{\rm Var}}

\newcommand{\rmvec}{{\rm\bf vec}}
\renewcommand{\Re}{{\rm Re}}

\newcommand{\Be}{{\rm Be}}

\newcommand{\sfp}{\mathsf{p}}
\newcommand{\sfP}{\mathsf{P}}

\newcommand{\E}{\mathbb E}

\newcommand{\p}{\mathbb P}
\newcommand{\R}{\mathbb R}

\newcommand{\1}{\mathbbm 1}

\newcommand{\bu}{{\bf u}}

\newcommand{\bx}{{\bf x}}

\newcommand{\bX}{{\bf X}}

\newcommand{\beps}{\boldsymbol \eps}

\newcommand{\bfeta}{\boldsymbol \eta}
\newcommand{\bzeta}{\boldsymbol \zeta}

\newcommand{\bxi}{\boldsymbol \xi}

\newcommand{\bzero}{{\bf 0}}

\newcommand{\cB}{{\mathcal{B}}}
\newcommand{\cE}{{\mathcal{E}}}

\newcommand{\cK}{{\mathcal{K}}}

\newcommand{\cN}{{\mathcal{N}}}
\newcommand{\cO}{{\mathcal{O}}}
\newcommand{\cT}{{\mathcal{T}}}

\newcommand{\tta}{\mathtt a}
\renewcommand{\i}{\mathtt i}
\newcommand{\ttm}{\mathtt m}
\newcommand{\ttr}{\mathtt r}
\newcommand{\ttR}{\mathtt R}
\newcommand{\tts}{\mathtt s}
\newcommand{\ttv}{\mathtt v}
\newcommand{\ttx}{\mathtt x}

\newcommand{\z}{\mathfrak z}

\newcommand{\myendproof}{\hfill$\square$}

\newtheorem{Th}{Theorem}[section]
\newtheorem{Lem}[Th]{Lemma}

\newtheorem{Prop}[Th]{Proposition}
\newtheorem{Co}[Th]{Corollary}
\newtheorem{Rem}[Th]{Remark}

\newtheorem{As}[Th]{Assumption}

\title{Sharper dimension-free bounds on the Frobenius distance\\between sample covariance and its expectation}

\author{
Nikita Puchkin\thanks{HSE University and IITP RAS, Russian Federation, npuchkin@hse.ru}
\and
Fedor Noskov\thanks{HSE University, fnoskov@hse.ru}
\and
Vladimir Spokoiny
\thanks{WIAS Berlin and Humboldt University, Germany, spokoiny@wias-berlin.de}
}

\date{}

\begin{document}

\maketitle

\begin{abstract}
    We study properties of a sample covariance estimate $\widehat \Sigma$ given a finite sample of $n$ i.i.d. centered random elements in $\R^d$ with the covariance matrix $\Sigma$. We derive dimension-free bounds on the squared Frobenius norm of $(\widehat\Sigma - \Sigma)$ under reasonable assumptions. For instance, we show that $\smash{\|\widehat\Sigma - \Sigma\|_{\rm F}^2}$ differs from its expectation by at most $\smash{\mathcal O({\rm{Tr}}(\Sigma^2) / n)}$ with overwhelming probability, which is a significant improvement over the existing results. This allows us to establish the concentration phenomenon for the squared Frobenius distance between the covariance and its empirical counterpart in the case of moderately large effective rank of $\Sigma$.
\end{abstract}

\section{Introduction}
\label{sec:introduction}

Covariance estimation is one of classical topics in multivariate statistics with applications in plenty of areas, including signal processing \citep{krim96, haghighatshoar18}, bioinformatics \citep{xie03, schafer05, hero12}, image analysis \citep{jorg00, zhang10}, and finance \citep{ledoit03, holtz10, bai11, fan15}. Let $\bX, \bX_1, \dots, \bX_n$ be i.i.d. centered random vectors in $\R^d$ with a covariance matrix $\E \bX \bX^\top = \Sigma$. In the present paper, we study properties of the most natural estimator of $\Sigma$, namely, the sample covariance, which is defined by the formula
\[
    \widehat \Sigma = \frac1n \sum\limits_{i=1}^n \bX_i \bX_i^\top.
\]
The question of how well $\widehat \Sigma$ approximates $\Sigma$ was extensively studied in the last decades. For instance, it arose in the paper of \citet*{kannan97} in the context of volume estimation. The authors used the sample covariance to bring a convex body $\cK \subset \R^d$ into a nearly isotropic position in an intermediate step of their algorithm. Assuming that $\bX_1, \dots, \bX_n$ are uniformly distributed on $\cK$, \citet*{kannan97} showed that, for any $\delta \in (0, 1)$, the operator norm of $(\widehat \Sigma - \Sigma)$ is of order $\cO(d / (\delta \sqrt{n}))$ with probability at least $(1 - \delta)$ (see their proof of Theorem 5.11). This result was improved in a series of works \citep{bourgain99, rudelson99, giannopoulos05, paouris06, adamczak10} until \citet*{adamczak11} obtained that, if the sample size $n$ is large enough, then
\[
    \left\| \widehat \Sigma - \Sigma \right\|
    \lesssim \sqrt{\frac{d}{n}}
\]
with overwhelming probability $(1 - \exp\{-\cO(\sqrt{d})\})$ for a large class of distributions, including the log-concave ones and the measures satisfying the Poincar\'{e} inequality. Here and further in this paper, the relation $\lesssim$ stands for inequality up to a positive multiplicative constant (see the notation section below). A bit sharper bound can be derived for random vectors with lighter tails. For example, in \cite[Corollary 5.50]{vershynin12b}, the author showed that, if $\bX_1, \dots, \bX_n$ are sub-Gaussian, then 
\begin{equation}
    \label{eq:sub-gaussian_worst_case}
    \left\| \widehat\Sigma - \Sigma \right\|
    \lesssim \sqrt{\frac{d + \log(1/\delta)}n} \vee \frac{d + \log(1/\delta)}n
\end{equation}
with probability at least $(1 - \delta)$. Here and further in this paper, $(a \vee b)$ stands for $\max\{a, b\}$. The inequality \eqref{eq:sub-gaussian_worst_case}, in particular, yields that the Frobenius norm of $(\widehat \Sigma - \Sigma)$ is of order $\cO(d / \sqrt{n})$ with high probability.

The bound \eqref{eq:sub-gaussian_worst_case} cannot be improved in the worst-case scenario. If the population covariance $\Sigma$ has $d$ large eigenvalues, then the number of parameters to be estimated is $\Omega(d^2)$. This means that the performance of $\widehat\Sigma$ suffers from the curse of dimensionality. In \citet{tropp12} and \citet{vershynin12b}, the authors noticed that if the distribution of $\bX_1, \dots, \bX_n$ is supported on a Euclidean ball of radius $R$ centered at the origin, then the following inequality holds with high probability:
\begin{equation}
    \label{eq:bounded_case}
    \left\| \widehat\Sigma - \Sigma \right\|
    \lesssim \sqrt{\frac{R \|\Sigma\| \log(d/\delta)}n} \vee \frac{R \log(d/\delta)}n.
\end{equation}
We see that the finite support assumption almost eliminates the dependence on $d$. However, we still face problems when the dimension becomes exponentially large compared to the sample size $n$. Fortunately, \citet[Remark 5.53]{vershynin12b} noticed that if the data lies near a low-dimensional subspace (which is often the case in high-dimensional tasks), then the rate of convergence of the sample covariance matrix depends on the intrinsic dimension, characterized by the effective rank
\[
    \ttr(\Sigma) = \frac{\Tr(\Sigma)}{\|\Sigma\|}.
\]
In \citet{adamczak15} and \citet{koltchinskii17}, the authors proved a dimension-free version of \eqref{eq:sub-gaussian_worst_case} in the Gaussian setup:
\begin{equation}
    \label{eq:gaussian_dimension-free}
    \left\| \widehat\Sigma - \Sigma \right\|
    \lesssim \|\Sigma\| \left( \sqrt{\frac{\ttr(\Sigma) + \log(1/\delta)}n} \vee \frac{\ttr(\Sigma) + \log(1/\delta)}n \right).
\end{equation}
As mentioned in \citep[Theorem 9]{koltchinskii17}, a similar inequality holds for a broad class of sub-Gaussian distributions 
satisfying the $\psi_2$-$L_2$-equivalence condition.
The word ``dimension-free'' means that the right-hand side of \eqref{eq:gaussian_dimension-free} depends only on the operator norm and effective rank of $\Sigma$. Hence, the bound \eqref{eq:gaussian_dimension-free} still makes sense even if the ambient dimension $d$ is huge. Later, this result was extended in the papers of \citet[Theorem 9.2.4 and Exercise 9.2.5]{vershynin18}, \citet[Theorem 1]{zhivotovskiy21}, and \citet{han22}. In \citep[Theorem 3]{zhivotovskiy21}, the author also obtained a dimension-free version of \citep[Theorem 4.1]{adamczak10} for the log-concave case. A bit earlier, \cite{bunea15} studied the behaviour of the Frobenius norm of $(\widehat \Sigma - \Sigma)$ under the same $\psi_2$-$L_2$-equivalence assumption. They proved that (see \citep[Proposition A.3]{bunea15})
\begin{equation}
    \label{eq:bunea_fr_norm_bound}
    \left\| \widehat\Sigma - \Sigma \right\|_{\Fr}
    \lesssim \|\Sigma\| \ttr(\Sigma) \left( \sqrt{\frac{\log(2/\delta)}n} \vee \frac{\log(2/\delta)}n \right)
\end{equation}
with probability at least $(1 - \delta)$. We would like to note that the low effective rank assumption is not the only way to escape the curse of dimensionality. There are a lot of papers devoted to the problem of covariance estimation under  
structural assumptions (for example, sparsity \citep{bickel08, cai12, fan15}, Kronecker product structure \citep{tsiligkaridis13, leng18}, bandable \citep{bickel08b, cai10} or Toeplitz \citep{xiao12, cai13} matrix $\Sigma$ to mention a few). This goes beyond the scope of the present paper. A reader is referred to the comprehensive survey \citep{cai16} on this subject. We consider the situation when $\Sigma$ has no other properties but a small effective rank.

The aforementioned works illustrate how rapidly the topic of covariance estimation evolved in the last decade. Recently developed advanced techniques allowed statisticians to examine more challenging setups, like those with missing observations \citep{lounici14, lounici23, abdalla23}, heavy tails \citep{vershynin12, srivastava13, youssef13, tikhomirov18, mendelson20, abdalla22} and adversarial contamination \citep{abdalla22, minasyan23}. We do not pursue the goal of pushing these results even further. In contrast, we use our technical findings to discover subtle effects that have not been noticed in the classical sub-Gaussian setting so far.

\bigskip

\noindent\textbf{Contribution.}\quad In the present paper, we study properties of the squared Frobenuis distance between the sample covariance and its expectation, which received less attention compared to the operator norm of $(\widehat \Sigma - \Sigma)$. The state-of-the-art result \eqref{eq:bunea_fr_norm_bound} of \citet{bunea15} implies that
\[
    \left\|\widehat \Sigma - \Sigma \right\|_{\Fr}^2 \lesssim \frac{\|\Sigma\|^2 \ttr(\Sigma)^2 \log(2/\delta)}n
\]
with high probability in the case of moderate confidence level $\delta$ (that is, $\log(1/\delta) \lesssim n$). This bound exhibits the correct dependence on the operator norm and on the effective rank of $\Sigma$. For instance, if $\bX_1, \dots, \bX_n$ have the Gaussian distribution $\cN(\bzero, \Sigma)$, then
\[
    \E \left\|\widehat \Sigma - \Sigma \right\|_{\Fr}^2 = \frac{(\Tr(\Sigma))^2 + \Tr(\Sigma^2)}n,
\]
and, according to the Paley-Zygmund inequality, $\|\widehat \Sigma - \Sigma \|_{\Fr}^2$ is not smaller than $\big[(\Tr(\Sigma))^2 + \Tr(\Sigma^2) \big] / (2n)$ on an event of positive probability.
However, the upper bound $\cO((\Tr(\Sigma))^2 / n)$ becomes suboptimal for some distributions, if we are speaking of the difference
$\|\widehat \Sigma - \Sigma\|_{\Fr}^2 - \E \|\widehat \Sigma - \Sigma\|_{\Fr}^2$. A thorough analysis of higher-order derivatives of the cumulant generating function
\begin{equation}
    \label{eq:phi}
    \varphi(U) = \log \E e^{\bxi^\top U \bxi},
    \quad \text{where} \quad
    \bxi = \Sigma^{-1/2} \bX,
\end{equation}
allows us to establish a dimension-free high-probability upper bound
\[
    \left| \left\|\widehat \Sigma - \Sigma \right\|_{\Fr}^2 - \E \left\|\widehat \Sigma - \Sigma \right\|_{\Fr}^2 \right|
    \lesssim \frac{\|\Sigma\|^2}n \max\left\{ \ttr(\Sigma^2) \sqrt{\log(8/\delta)}, \log(8/\delta) \right\},
    \quad \text{where} \quad
    \ttr(\Sigma^2) = \frac{\Tr(\Sigma^2)}{\|\Sigma\|^2},
\]
which holds for a large subclass of sub-Gaussian distributions (see Assumption \ref{as:fourth_moment} and Proposition \ref{prop:orlicz} below), arbitrary $\delta \in (0, 1)$, such that $\log(1/\delta) \lesssim n / \ttr(\Sigma)^2$, and $n \gtrsim \ttr(\Sigma)^6$ (see Theorem \ref{th:fr_norm_upper_bound} and Theorem \ref{th:fr_norm_lower_bound}). We would like to note that $\ttr(\Sigma^2)$, also referred to as stable rank of $\Sigma$, is always not greater than $\ttr(\Sigma)^2$ and sometimes it can be significantly smaller than the squared effective rank. As a consequence, we obtain the following bound on $\|\widehat \Sigma - \Sigma \|_{\Fr}^2$ with a sharp constant in a high-dimensional setting (Corollary \ref{co:concentration}):
\[
    \left| \frac{\|\widehat \Sigma - \Sigma \|_{\Fr}^2}{\E \|\widehat \Sigma - \Sigma \|_{\Fr}^2} - 1 \right|
    \lesssim \frac{ \sqrt{\log(8/\delta)}}{\ttr(\Sigma) - 1}
    = \cO\left( \frac1{\ttr(\Sigma)} \right),
    \quad \ttr(\Sigma), \, n / \ttr(\Sigma)^6 \rightarrow \infty.
\]

\bigskip

\noindent\textbf{Motivation.}\quad
Our results on concentration of $\|\widehat \Sigma - \Sigma\|_{\Fr}$ around its expectation have direct implications in hypothesis testing and construction of confidence bands. The squared Frobenius norm is a popular tool for construction of statistical criterions for testing a simple hypothesis or two-sample tests (see, for instance, \citep{ledoit02, chen10, li12, han20}). The reason for that can be easily explained by Theorems \ref{th:fr_norm_upper_bound} and \ref{th:fr_norm_lower_bound}. 
Assume that our goal is to check whether $\Sigma$ is equal to a given matrix $\Sigma_0$. Suppose that there is $\omega > 0$, such that the inequality
\begin{equation}
    \label{eq:orlicz_norm_condition}
    \left\|\bX^\top V \bX - \Tr(V\Sigma) \right\|_{\psi_1} \leq \omega \left\|\Sigma^{1/2} V \Sigma^{1/2} \right\|_{\Fr}
\end{equation}
holds for all matrices $V \in \R^{d \times d}$. Here $\|\cdot\|_{\psi_1}$ stands for the Orlicz norm (see the definition \eqref{eq:orlicz_norm_def} below). In Section \ref{sec:main_results}, we discuss examples of probability distributions satisfying \eqref{eq:orlicz_norm_condition}. Under this condition,
\[
    \widehat T = \frac1n \sum\limits_{i = 1}^n \|\bX\|^4 + \frac1{\lfloor n/2 \rfloor} \sum\limits_{i = 1}^{\lfloor n/2 \rfloor} \big( \bX_{2i - 1}^\top \bX_{2i} \big)^2
\]
is a consistent estimate $\widehat T$ of $n \, \E \|\widehat\Sigma - \Sigma\|_{\Fr}^2 = \E \|\bX\|^4 - \Tr(\Sigma^2)$. Then, under the null hypothesis $\mathcal H_0 : \Sigma = \Sigma_0$, the statistic $\|\widehat\Sigma - \Sigma_0\|_{\Fr}^2 - \widehat T / n$ will be as large as $\cO(\Tr(\Sigma^2) / n)$. On the other hand, if $\Sigma \neq \Sigma_0$, then we have
\[
    \|\widehat\Sigma - \Sigma_0\|_{\Fr}^2 - \widehat T / n
    = \|\widehat\Sigma - \Sigma\|_{\Fr}^2 - \widehat T / n + 2 \Tr\left[ (\Sigma - \Sigma_0)(\widehat\Sigma - \Sigma) \right] + \|\Sigma - \Sigma_0\|_{\Fr}^2.
\]
Due to \eqref{eq:orlicz_norm_condition}, it holds that $\|\widehat\Sigma - \Sigma\|_{\Fr}^2 - \widehat T / n = \cO(\Tr(\Sigma^2) / n)$ and
\begin{align*}
    \Tr\left[ (\Sigma - \Sigma_0)(\widehat\Sigma - \Sigma) \right]
    &
    = \frac1n \sum\limits_{i = 1}^n \left( \bX_i^\top (\Sigma - \Sigma_0) \bX_i - \Tr\left[(\Sigma - \Sigma_0) \Sigma\right] \right)
    \\&
    = \cO\left( \frac1{\sqrt n} \left\| \bX^\top (\Sigma - \Sigma_0) \bX - \Tr\left[(\Sigma - \Sigma_0) \Sigma\right] \right\|_{\psi_1} \right)
    \\&
    = \cO\left( \frac1{\sqrt n} \left\|\Sigma^{1/2} (\Sigma - \Sigma_0) \Sigma^{1/2} \right\|_{\Fr}\right)
\end{align*}
with high probability. Hence, the test statistic
$\|\widehat\Sigma - \Sigma_0\|_{\Fr}^2 - \widehat T / n$ will be of order
\[
    \|\Sigma - \Sigma_0\|_{\Fr}^2 - \cO\left( \frac1{\sqrt n} \left\|\Sigma^{1/2} (\Sigma - \Sigma_0) \Sigma^{1/2} \right\|_{\Fr}\right) - \cO\left(\frac{\Tr(\Sigma^2)}n \right).
\]
This leads to a detection boundary $\|\Sigma - \Sigma_0\|_{\Fr} \gtrsim \|\Sigma\|_{\Fr} / \sqrt{n}$. This does not follow from the result of \citet{bunea15}.

\bigskip

\noindent\textbf{Notation.}\quad
Throughout the paper, $\|A\|_{\Fr}$ and $\|A\|$ stand for the Frobenius and the operator norm of $A$, respectively. If the matrix $A$ is symmetric positive semidefinite and $A \neq O$, where $O$ is the matrix with zero entries, we denote its effective rank by
\[
    \ttr(A) = \Tr(A) / \|A\|.
\]
We use the standard notation $(A \otimes B)$ for the Kronecker product of matrices $A$ and $B$. Its definition and some useful properties are listed in Appendix \ref{sec:kronecker}. For any matrix $U \in \R^{p \times q}$ with columns $\bu_1, \dots, \bu_q$, the vectorization operator is given by $\rmvec(U) = (\bu_1^\top, \dots, \bu_q^\top)^\top \in \R^{pq}$. Here and further in the paper, the bold font is reserved for vectors, while matrices and scalars are displayed in regular font. For a random variable $\eta$, its Orlicz $\psi_s$-norm, $s \geq 1$, is defined as
\begin{equation}
    \label{eq:orlicz_norm_def}
    \|\eta\|_{\psi_s} = \inf\left\{ t > 0 : \E e^{|\eta|^s / t^s} \leq 2 \right\}.
\end{equation}
The expressions $(a \vee b)$ and $(a \wedge b)$ denote $\max\{a, b\}$ and $\min\{a, b\}$, respectively.
Sometimes, instead of the standard $\cO$ notation, we use $f \lesssim g$ or $g \gtrsim f$, which means that there is a universal constant $c > 0$, such that $f \leq c g$.
Finally, $\cB(\bx, R)$ is the Euclidean ball of radius $R$ centered at $\bx$.

\bigskip

\noindent\textbf{Paper structure.}\quad The rest of the paper is organized as follows. In Section \ref{sec:main_results}, we present our main results and discuss their implications. After that, we illustrate them with some numerical simulations in Section \ref{sec:numerical}. Section \ref{sec:proofs} is devoted to the proofs of the statements from Section \ref{sec:main_results}, some technical details are moved to Appendix.

\section{Main results}
\label{sec:main_results}

In this section, we present rigorous statements of the main results of the paper. Let us start with a couple of auxiliary definitions. Throughout the paper, $\bxi = \Sigma^{-1/2} \bX$ stands for the whitened vector. We assume that the covariance matrix $\Sigma$ is invertible for simplicity. In a general situation, one can always reduce the problem of covariance estimation to the non-degenerate case, adding tiny isotropic Gaussian noise to the observations. For a matrix $U \in \R^{d \times d}$, such that $\E e^{\bxi^\top U \bxi} < \infty$, let $\sfP_{U}$ denote a probability measure, defined as
\begin{equation}
    \label{eq:pu}
	\rmd \sfP_{U}(\bx) = \frac{e^{\bx^\top U \bx}}{\E e^{\bxi^\top U \bxi}} \rmd \p_{\bxi}(\bx),
    \quad \text{where} \quad
    \bxi = \Sigma^{-1/2} \bX.
\end{equation}
For any Borel function $f : \R^d \times \mathbb \R^{d \times d} \rightarrow \mathbb \R$, the expression $\sfP_{U} f(\bxi, U)$ stands for the expectation of $f(\bxi, U)$ with respect to the measure $\sfP_{U}$:
\begin{equation}
    \label{eq:pu_expectation}
	\sfP_{U} f(\bxi, U) = \frac{\E f(\bxi, U) e^{\bxi^\top U \bxi}}{\E e^{\bxi^\top U \bxi}}.
\end{equation}
As we mentioned in the introduction, our proof is based on the study of derivatives of the cumulant generating function $\varphi(U)$, defined in \eqref{eq:phi}. They admit a nice representation in terms of the introduced functional $\sfP_U$. The only property we require is the regularity of the third and the fourth derivatives of $\varphi(U)$, which is guaranteed by the following assumption.
\begin{As}
    \label{as:derivatives}
    There exist positive numbers $\tau$ and $\rho_{\max}$ such that the vector $\bxi = \Sigma^{-1/2} \bX$ satisfies the inequality
    \[
        \sfP_U \left(\bxi^\top V \bxi - \sfP_U \bxi^\top V \bxi \right)^4 \leq \tau^2 \; \|V\|_{\Fr}^4
    \]
    for all $U \in \R^{d \times d}$ such that $\|U\|_{\Fr} \leq \rho_{\max}$ and all $V \in \R^{d \times d}$.
\end{As}
To our knowledge, Assumption \ref{as:derivatives} has not been used in the literature on covariance estimation, so the natural question is what are the examples of random vectors $\bxi$ satisfying this condition. We claim that Assumption \ref{as:derivatives} holds for a broad class of distributions and support our assertion with the next proposition.
\begin{Prop}
    \label{prop:orlicz}
    Assume that there exists $\omega > 0$ such that
    \begin{equation}
        \label{eq:orlicz_norm_inequality}
        \left\|\bxi^\top V \bxi - \Tr(V) \right\|_{\psi_1} \leq \omega \|V\|_{\Fr}
        \quad \text{for all $V \in \R^{d \times d}$}.
    \end{equation}
    Then the random vector $\bxi$ satisfies Assumption \ref{as:derivatives} with $\tau = 64 \omega^2$ and $\rho_{\max} = (6 \omega)^{-1}$.
\end{Prop}
Proposition \ref{prop:orlicz} immediately yields that Assumption \ref{as:derivatives} is fulfilled for all random vectors $\bxi$ satisfying the Hanson-Wright inequality. For instance, this is the case for sub-Gaussian random vectors with independent entries (see, e.g., \citep{hanson71, rudelson13}). Less trivial examples are the distributions with the convex concentration property \citep{adamczak15}. The random vector $\bxi$ is said to have the convex concentration property if there exists a positive constant $K$ such that
\[
    \p\left( \left|g(\bxi) - \E g(\bxi)\right| \geq t \right) \leq 2 e^{-t^2 / K^2}
\]
for any convex $1$-Lipschitz function $g$ and any $t > 0$. We would like to note that random vectors with strongly log-concave density or satisfying the log-Sobolev inequality have such property. It is also easy to verify that the distributions satisfying the Hanson-Wright inequality have the $\psi_2$-$L_2$-equivalence property. However, we would like to note that Assumption \ref{as:derivatives} (and, consequently, \eqref{eq:orlicz_norm_inequality}) 
is a stronger condition than the one considered in \citep{bunea15}. We support our claim with the following proposition.

\begin{Prop}
    \label{prop:counterexample}
    There exists a random vector $\bxi \in \R^d$ satisfying the following properties:
    \begin{itemize}
        \item $\E \bxi = \bzero$ and $\E \bxi \bxi^\top = I_d$;
        \item for any $\bu \in \R^d$ it holds that $\left\| \bxi^\top \bu \right\|_{\psi_2}^2 \leq 16 \, \E \| \bxi^\top \bu \|^2 / 3$;
        \item there exists a matrix $V \in \R^{d \times d}$ such that
        \[
            \E \left(\bxi^\top V \bxi - \E \bxi^\top V \bxi \right)^4 = d^2 \; \left\|V \right\|_{\Fr}^4.
        \]
    \end{itemize}
\end{Prop}
The counterexample in Proposition \ref{prop:counterexample} is quite straightforward. We take $\bxi$ equal to a random vector with i.i.d. Rademacher entries multiplied by a bounded random variable. We postpone the proof of Proposition \ref{prop:counterexample} to Section \ref{sec:prop_counterexample_proof} and proceed with a high-probability upper bound on the Frobenius norm of $(\widehat\Sigma - \Sigma)$.
\begin{Th}
    \label{th:fr_norm_upper_bound}
    Suppose that Assumption \ref{as:derivatives} holds with $\tau \geq 2$. Fix any $\delta \in (0, 1)$ and let
    \[
        \ttR(\Sigma, \delta) = \frac{3\ttr(\Sigma)}{2\sqrt{2}} + \frac12 \sqrt{\ttr(\Sigma^2)} + 2 \sqrt{e \log(3 / \delta)}.
    \]
    Assume that the sample size $n$ satisfies the inequalities
    \begin{equation}
        \label{eq:n_admissible}
        \|\Sigma\| \, \ttR(\Sigma, \delta) \leq \rho_{\max} \sqrt{n}
        \quad \text{and} \quad
        \ttr(\Sigma)^2 \, \ttR(\Sigma, \delta)^2 \left( \frac{\ttR(\Sigma, \delta)^2}{4} + 3 \right) \leq 36 n.
    \end{equation}
    Then, with probability at least $(1 - \delta)$, it holds that
    \begin{equation}
        \label{eq:fr_norm_upper_bound}
        \left\|\widehat \Sigma - \Sigma \right\|_{\Fr}^2 - \E \left\|\widehat \Sigma - \Sigma \right\|_{\Fr}^2
        < \frac{4 \|\Sigma\|^2}n \max\left\{ \sqrt{2 \big( \tau \ttr(\Sigma^2)^2 + \tau^2 \ttr(\Sigma^4) \big) \log(3 / \delta)}, \; 4e \tau \log(3 / \delta) \right\}.
    \end{equation}
\end{Th}
According to Lemma \ref{lem:mean} below, under Assumption \ref{as:derivatives}, it holds that
\[
    n \; \E \left\|\widehat \Sigma - \Sigma \right\|_{\Fr}^2 \leq \big( \Tr(\Sigma) \big)^2 + (\tau  - 1) \; \Tr(\Sigma^2).
\]
Hence, Theorem \ref{th:fr_norm_upper_bound} yields that
\[
    \left\|\widehat \Sigma - \Sigma \right\|_{\Fr}^2 
    \lesssim \frac{\|\Sigma\|^2 \left(\ttr(\Sigma)^2 + \log(3/\delta) \right)}n
\]
with probability at least $(1 - \delta)$, resembling the result of \cite{bunea15} with slightly better dependence on $\delta$ (we have $\ttr(\Sigma)^2 + \log(2/\delta)$, instead of $\ttr(\Sigma)^2 \log(1/\delta)$, cf. \eqref{eq:bunea_fr_norm_bound}). However, as we show in Corollary \ref{co:concentration} a bit later, the expression in the right-hand side of \eqref{eq:fr_norm_upper_bound} may be much smaller than the expectation of $\|\widehat\Sigma - \Sigma\|_{\Fr}^2$.

The result of Theorem \ref{th:fr_norm_upper_bound} is also related to concentration of quadratic forms, because the squared Frobenius norm of $(\widehat \Sigma - \Sigma)$ can be naturally represented as
\[
    n \left\|\widehat\Sigma - \Sigma\right\|_{\Fr}^2 = \rmvec(H)^\top \left( \Sigma \otimes \Sigma \right) \rmvec(H),
\]
where $\otimes$ stands for the Kronecker product of matrices and $\rmvec(H)$ is a random vector, obtained by reshaping of
\[
    H = \frac1{\sqrt n} \sum\limits_{i = 1}^n \left(\bxi_i \bxi_i^\top - I_d \right).
\]
It is easy to observe that the matrix $H$ may have sub-exponential dependent entries under Assumption \ref{as:derivatives}. At the same time, most of the papers studying the quadratic forms work with sub-Gaussian random vectors (see, for instance, \citep{laurent00, hsu12, klochkov20, spokoiny23} and the previously mentioned \citep{hanson71, rudelson13, adamczak15}).
The only exception we are aware of is the paper of \cite{sambale23}, where the author considers quadratic forms of sub-exponential random vectors with independent components. One may also consider $\|\widehat\Sigma - \Sigma\|_{\Fr}^2$ as a polynomial of a sub-Gaussian random vector of degree $4$. Unfortunately, the latest results on concentration of sub-Gaussian polynomials (see, for instance, \cite{schudy12,adamczak15b,gotze21}) cannot recover the upper bound \eqref{eq:fr_norm_upper_bound}.

The final remark on Theorem \ref{th:fr_norm_upper_bound} we want to make is that $\|\widehat\Sigma - \Sigma\|_{\Fr}^2$ should behave like a square of a sub-exponential random variable. However, in the right-hand side of \eqref{eq:fr_norm_upper_bound}, we have the standard sub-exponential tail $\sqrt{\log(2/\delta)} \vee \log(2/\delta)$, rather than $(\log(2/\delta))^2$. This is because the statement of Theorem \ref{th:fr_norm_upper_bound} holds not for all $\delta$ from the interval $(0, 1)$ but only for $\delta \gtrsim e^{-\cO(n)}$, excluding the case of extremely high confidence. This is a consequence of our approach, relating the probability of large deviation of $\|\widehat \Sigma - \Sigma\|_{\Fr}^2$ with the analysis of the restricted exponential moment
\begin{equation}
    \label{eq:restricted_exp_moment}
    \E \exp\left\{ \frac{\lambda n}2 \|\widehat \Sigma - \Sigma\|_{\Fr}^2 \right\} \1\left( \sqrt{n} \left\|\Sigma^{1/2} (\widehat \Sigma - \Sigma) \Sigma^{1/2} \right\|_{\Fr} \leq \z \right),
\end{equation}
where $\z$ is an properly chosen constant (see Lemma \ref{lem:linearization} and Lemma \ref{lem:restricted_moment} below). Since the expectation \eqref{eq:restricted_exp_moment} may blow up as $\lambda > 0$ grows, we have introduce the technical requirements \eqref{eq:n_admissible} to prevent this undesired event.

We proceed with a complementary lower bound on $\|\widehat\Sigma - \Sigma\|_{\Fr}^2$. In contrast to the upper bound, we do not need the random vector $\bxi$ to be sub-Gaussian anymore. The following will be sufficient for our purposes.

\begin{As}
    \label{as:fourth_moment}
    There exists $\alpha > 0$ such that
    \[
        \E \left( \bxi^\top V \bxi - \E \, \bxi^\top V \bxi \right)^4 
        \leq \alpha^2 \; \|V\|_{\Fr}^4
        \quad \text{for all $V \in \R^{d \times d}$.}
    \]
\end{As}

The reason for the significant relaxation of Assumption \ref{as:derivatives} is that in the proof of lower bounds we have to study the exponential moment $\E \exp\{-\lambda \|\widehat \Sigma - \Sigma\|_{\Fr}^2\}$, which is finite for all $\lambda > 0$ whatever the distribution of $\bxi$ is. As a result, the weaker Assumption \ref{as:fourth_moment} does not lead to worse dependence on the parameter $\delta$. A similar effect was observed in \citep{oliveira16}, where the author studied lower tails of the sample covariance.

\begin{Th}
    \label{th:fr_norm_lower_bound}
    Grant Assumption \ref{as:fourth_moment} and let the sample size $n$ be sufficiently large in a sense that
    \begin{equation}
        \label{eq:n_large}
        5 \ttr(\Sigma^2)\left( 3 \ttr(\Sigma)^2 + \frac{2 \ttr(\Sigma^2) \ttr(\Sigma)^2}{\alpha}\right)
        \leq n
        \quad \text{and} \quad
        96 \left(1 + \sqrt{\alpha} \right)^2 \left( 3 \ttr(\Sigma)^2 + \frac{2 \ttr(\Sigma^2) \ttr(\Sigma)^2}\alpha \right) \leq \alpha n.
    \end{equation}
    Then, for any $\delta \in (0, 1)$, with probability at least $(1 - \delta)$, it holds that
    \begin{equation}
        \label{eq:fr_norm_lower_bound}
        \left\|\widehat \Sigma - \Sigma \right\|_{\Fr}^2 - \E \left\|\widehat \Sigma - \Sigma \right\|_{\Fr}^2
        > -\frac{4 \|\Sigma\|^2}n \left( 2 \alpha \log(5/\delta)
        \vee \sqrt{2 \mathfrak R(\Sigma) \log(5/\delta)} \right),
    \end{equation}
    where
    \[
        \mathfrak R(\Sigma) = \alpha \ttr(\Sigma^2)^2 + \alpha^2 \ttr(\Sigma^4) + \frac{5 \alpha^2 \ttr(\Sigma)^2}{2 n} \left( 3 \ttr(\Sigma)^2 +  \frac{2 \ttr(\Sigma^2) \ttr(\Sigma)^2}{\alpha} \right).
    \]
\end{Th}

Theorems \ref{th:fr_norm_upper_bound} and \ref{th:fr_norm_lower_bound} imply that
\begin{equation}
    \label{eq:fr_norm_two_side_bound}
    \left| \left\|\widehat \Sigma - \Sigma \right\|_{\Fr}^2 - \E \left\|\widehat \Sigma - \Sigma \right\|_{\Fr}^2 \right|
    \lesssim \frac{\|\Sigma\|^2 \ttr(\Sigma^2)}n
\end{equation}
with high probability, provided that Assumption \ref{as:derivatives} is fulfilled. For instance, this is the case when $\bxi$ is a Gaussian random vector. To ensure that the bound \eqref{eq:fr_norm_two_side_bound} is tight, it is enough to recall that the variance of $\|\widehat \Sigma - \Sigma \|_{\Fr}^2$ is of order $\|\Sigma\|^4 \ttr(\Sigma^2)^2 / n^2$ in the Gaussian case. We can go even further examine the relative fluctuation of the squared Frobenius norm of $(\widehat \Sigma - \Sigma)$. Summing up the results of Theorem \ref{th:fr_norm_upper_bound} and Theorem \ref{th:fr_norm_lower_bound} and using the union bound, we obtain the following corollary.

\begin{Co}
\label{co:concentration}
Grant Assumption \ref{as:derivatives} and let $\Sigma$, $\rho_{\max}$, $\tau \geq 2$, and $\delta \in (0, 1)$ (possibly depending on $n$) be such that
\begin{equation}
    \label{eq:asymptotic_conditions}
    \frac{\Tr(\Sigma)^2 \vee \log(8 / \delta)}{n \rho_{\max}^2} \rightarrow 0,
    \quad
    \frac{\ttr(\Sigma)^6}n \rightarrow 0,
    \quad \text{and} \quad
    \frac{\ttr(\Sigma) \log(8 / \delta)}{\sqrt n} \rightarrow 0
    \quad \text{as $n \rightarrow \infty$.}
\end{equation}
Then, with probability at least $(1 - \delta)$, it holds that
\begin{equation}
    \label{eq:statistic_in_experiments}
    \left| \frac{\|\widehat \Sigma - \Sigma \|_{\Fr}^2}{\E \|\widehat \Sigma - \Sigma \|_{\Fr}^2} - 1 \right|
    \lesssim \tau \max\left\{ \frac{\sqrt{\log(8 / \delta)}}{\ttr(\Sigma) - 1}, \frac{\log(8 / \delta)}{\ttr(\Sigma)\big( \ttr(\Sigma) - 1 \big)} \right\},
\end{equation}
where $\lesssim$ stands for inequality up to an absolute constant.
\end{Co}
Corollary \ref{co:concentration} establishes that the ratio $\|\widehat\Sigma - \Sigma\|_{\Fr}^2 / \E \|\widehat\Sigma - \Sigma\|_{\Fr}^2$ is close to $1$, when the effective rank of $\Sigma$ is large and $n / \ttr(\Sigma)^6 \rightarrow \infty$ as $n$ tends to infinity. Though its proof is straightforward, we postpone it to Section \ref{sec:co_concentration_proof} and focus on discussion of the result. First, we would like to emphasize that Assumption \ref{as:derivatives} plays a crucial role in this phenomenon. Our numerical experiments, reported in Section \ref{sec:numerical}, show the absence of concentration when we take the random vector $\bxi$ as described in the proof of Proposition \ref{prop:counterexample}. Second, one can apply Corollary \ref{co:concentration} to bound the fluctuation of
\[
    \left\|\widehat \Sigma - \Sigma \right\|_{\Fr} - \E \left\|\widehat \Sigma - \Sigma \right\|_{\Fr}.
\]
\begin{Th}
    \label{th:fr_norm_fluctuation}
    Assume the conditions of Corollary \ref{co:concentration}. In addition, suppose that $\ttr(\Sigma) / \sqrt{\log(8 / \delta)} \rightarrow \infty$ and $\tau \ttr(\Sigma) / n \rightarrow 0$ as $n$ approaches infinity. Then, with probability at least $(1 - \delta)$, it holds that
    \[ 
        \left| \|\widehat \Sigma - \Sigma \|_{\Fr} - \E \|\widehat \Sigma - \Sigma \|_{\Fr} \right| \lesssim \tau \|\Sigma\| \sqrt{\left( 1 \vee \frac{\tau}{\ttr(\Sigma)} \right) \cdot \frac{\log(8 / \delta)}n},
        \quad n \rightarrow \infty.
    \]
    As before, $\lesssim$ denotes inequality up to an absolute constant.
\end{Th}
The proof of Theorem \ref{th:fr_norm_fluctuation} is moved to Section \ref{sec:th_fr_norm_fluctuation_proof}. It is known that, if $\bxi \sim \cN(\bzero, I_d)$ is a standard Gaussian vector in $\R^d$, then the difference $\|\bxi\|^2 - \E \|\bxi\|^2$ is of order $\sqrt{d}$ while $\big| \|\bxi\| - \E \|\bxi\| \big| \lesssim 1$. We observe a similar phenomenon for the Frobenius norm of $(\widehat \Sigma - \Sigma)$. Indeed, one can take $\tau = \cO(1)$ in the Gaussian case. Then, according to Theorem \ref{th:fr_norm_fluctuation}, it holds that
\[
    \left| \|\widehat \Sigma - \Sigma \|_{\Fr} - \E \|\widehat \Sigma - \Sigma \|_{\Fr} \right| \lesssim \|\Sigma\| \sqrt{\frac{\log(8 / \delta)}n}
\]
with probability at least $(1 - \delta)$ even if $\ttr(\Sigma)$ is large. On the other hand, Theorems \ref{th:fr_norm_upper_bound} and \ref{th:fr_norm_lower_bound} only guarantee that
\[
    \left| \left\|\widehat \Sigma - \Sigma \right\|_{\Fr}^2 - \E \left\|\widehat \Sigma - \Sigma \right\|_{\Fr}^2 \right|
    \lesssim \frac{\|\Sigma\|^2}n \max\left\{ \ttr(\Sigma^2) \sqrt{\log(8/\delta)}, \log(8/\delta) \right\},
\]
where the right-hand side is of order $d / n$ in the worst-case scenario $\Sigma = I_d$.

\section{Experiments}
\label{sec:numerical}

In this section, we illustrate the concentration phenomenon discussed in Section~\ref{sec:main_results} with numerical simulations. Our goal is to observe that the ratio $\|\widehat \Sigma - \Sigma \|_{\Fr}^2 / \E \|\widehat \Sigma - \Sigma \|_{\Fr}^2$ concentrates around $1$ when $\ttr(\Sigma)$ tends to infinity, as guaranteed by \eqref{eq:statistic_in_experiments}.

\begin{figure}[ht]
    \centering
    \includegraphics[width=\textwidth]{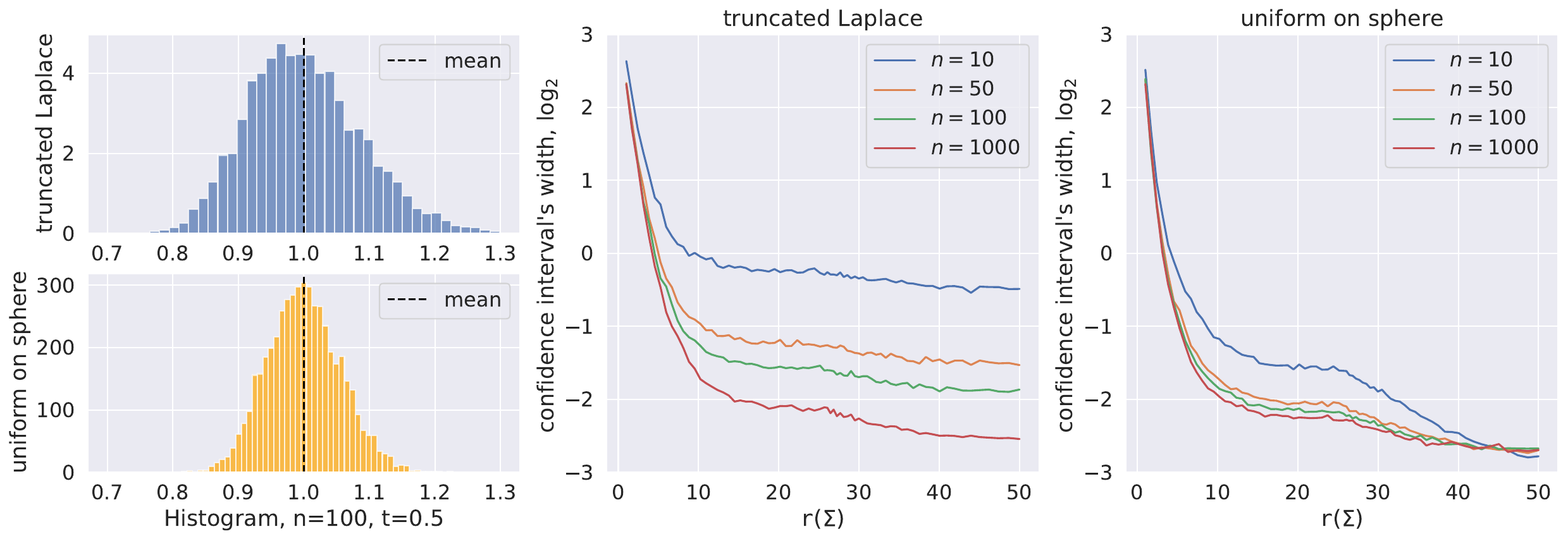}
    \caption{concentration of $\|\widehat{\Sigma} - \Sigma\|_{\Fr}^2 / \E \|\widehat{\Sigma} - \Sigma\|_{\Fr}^2$ around 1 as ${\tt r}(\Sigma)$ tends to infinity.}
    \label{fig:concentration experiment}
\end{figure}

Let us describe the process of population covariance generation.
We set the ambient dimension $d$ equal to $50$. Then, for each $t \in [0, 1]$, we define the covariance matrix $\Sigma_t = U \Lambda_t U^\top, t \in [0,1]$, where $U$ is a drawn from the uniform distribution on the orthogonal group $\mathbb{O}(d)$ and $\Lambda_t$ is defined as follows:
\begin{align*}
    \Lambda_t = 
    \begin{cases}
        \diag \left ( 1, 2 t \left (1 - \frac 1 d \right ), 2t \left ( 1 - \frac 2 d \right ), \ldots, \frac{2 t}{d}  \right ) & \text{ if } t \in [0, 0.5]; \\
        \diag \left ( 1, (1 - 1/d)^{2 (1 - t)}, (1 - 2/d)^{2 (1 - t)}, \dots, (1/d)^{2 (1 - t)}\right ), & \text{ otherwise.}
    \end{cases}
\end{align*}
In our simulations, we take $t$ from the grid $\cT = \{0, 1/69, 2/69; \ldots, 1\}$, so the effective ranks $\mathtt{r} (\Sigma_t)$, $t \in \cT$, are spread regularly over the segment $[1, d]$.

The simulation process goes as follows.
For each $\Sigma = U \Lambda U^\top \in \{\Sigma_t : t \in \cT\}$,  we consider two different distribution laws of $\bX_1, \ldots, \bX_n$.
\begin{enumerate}
    \item \label{item: trunc laplacian}
    The first variant is $\bX = U \Lambda^{1/2} \bxi / \sigma$, where $\bxi = (\xi_1, \ldots, \xi_d)^\top$ is a random vectors with independent entries drawn from the \emph{truncated Laplace} distribution with the density
    \begin{align*}
        \sfp(x) = \frac{1}{Z} e^{-|x|} \1 (|x| \le 6),
    \end{align*}
     where $Z$ is the normalizing constant. The constant $\sigma > 0$ is chosen is a way to ensure that $\E \bX \bX^\top = \Sigma$:
     \[
        \sigma^2 = \int\limits_{-6}^{6} x^2 \sfp(x) \rmd x.
    \]
    \item \label{item: unifrom on sphere}
    The second case we consider is $\bX = \sqrt{d} U \Lambda^{1/2} \bxi$, where $\bxi$ has the uniform distribution on the unit sphere ${\mathbb S}^{d - 1}$.
\end{enumerate}
Using numerical integration, we calculate $\E \Vert \widehat{\Sigma} - \Sigma_t \Vert_{\Fr}^2$ for the truncated Laplace distribution with a reasonable accuracy:
\begin{align*}
    n \E \left\Vert \widehat \Sigma - \Sigma_t \right\Vert_{\Fr}^2 = \E \left\Vert \bX \bX^\top - \Sigma_t \right\Vert_{\Fr}^2 = (\Tr(\Lambda_t))^2 + \left(  K - 2 \right) \Tr(\Lambda_t^2),
\end{align*}
where $K$ is the kurtosis of the truncated Laplace distribution, defined as
\[
    K = \frac{1}{\sigma^4} \int\limits_{-6}^6 x^4 \sfp(x) \rmd x.
\]
Similarly, we compute $\E \Vert \widehat{\Sigma} - \Sigma_t \Vert_{\Fr}^2$ for the uniform distribution on the unit sphere:
\begin{align*}
    n \E \left\Vert \widehat{\Sigma} - \Sigma_t \right\Vert_{\Fr}^2
    = \E \left \Vert \bX \bX^\top - \Sigma_t \right\Vert_{\Fr}^2 = \frac{d}{d + 2} \left (\Tr(\Lambda_t) \right )^2 + \frac{d - 2}{d + 2} \, \Tr (\Lambda^2_t).
\end{align*}
Note that both distributions fulfill the condition of Proposition~\ref{prop:orlicz} since in the former case $\bxi$ is a vector of independent bounded random variables, and in the second case $\bxi$ satisfies the logarithmic Sobolev inequality (see, for instance, \citep[Theorem 5.7.4]{bakry14}). Consequently, they admit concentration bounds provided in Theorem~\ref{th:fr_norm_upper_bound} and Theorem~\ref{th:fr_norm_lower_bound}.

For each $n \in \{10, 50, 100, 1000\}$, $t \in \cT = \{0, 1/69, 2/69; \ldots, 1\}$ and both distributions, we generate $5000$ samples $(\bX_1^{j, t}, \ldots, \bX_n^{j, t})$, $j = 1, \ldots, 5000$, of size $n$. Next, we compute
\begin{align*}
    a_{j, t} := \frac{\Vert \widehat{\Sigma}_{j, t} - \Sigma_t \Vert_{\Fr}^2}{\E \Vert \widehat \Sigma_{j, t} - \Sigma_t \Vert_{\Fr}^2},
\end{align*}
where $\widehat{\Sigma}_{j, t}$ is the empirical covariance matrix based on $\bX_1^{j,t}, \ldots, \bX^{j, t}_n$. For a fixed $t \in \cT$, we calculate the width $w_t$ of the empirical $0.95$-confidence interval for $\{a_{j,t} : 1 \leq j \leq 5000\}$. 

Finally, we plot the dependence of $\log_2(w_t)$ on $\mathtt{r}(\Sigma_t)$. The results are displayed in Figure~\ref{fig:concentration experiment}. We observe that the width of the interval goes to zero as $\mathtt{r}(\Sigma_t)$ increases.

\begin{figure}
    \centering
    \includegraphics[width=\textwidth]{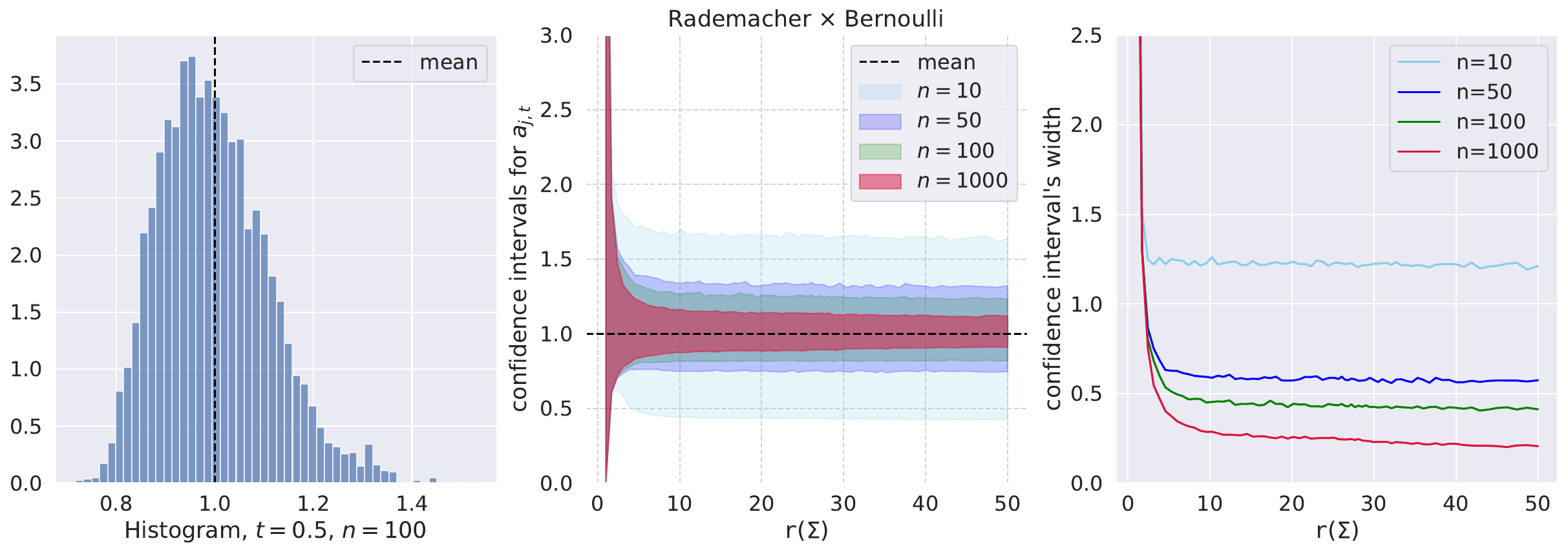}
    \caption{Absence of concetration of $\|\widehat{\Sigma} - \Sigma\|_{\Fr}^2 / \E \|\widehat{\Sigma} - \Sigma\|_{\Fr}^2$ for the counterexample suggested in Proposition~\ref{prop:counterexample}.}
    \label{fig:counterexample}
\end{figure}

Next, we support the necessity of Assumption~\ref{as:derivatives} with numerical experiments. Proposition~\ref{prop:counterexample} states that there exists a random vector $\bxi$ with unit covariance that satisfies $\psi_2$-$L_2$-equivalence but violates Assumption~\ref{as:derivatives}. The desired random vector $\bxi$ is explicitly constructed in the proof of Proposition~\ref{prop:counterexample}. To show that it does not admit the studied concentration phenomenon, we apply the same framework as before in this section. For each $t \in \cT = \{0, 1/69, \ldots, 68/69, 1\}$ and the corresponding covariance matrix $\Sigma_t = U \Lambda_t U^\top$, we consider the random vector $\bX = U \Lambda_t^{1/2} \bxi$. Then, the sample covariance matrix satisfies
\begin{align*}
    n \E \Vert \widehat \Sigma - \Sigma_t \Vert_{\Fr}^2 = \E \Vert \bX \bX^\top - \Sigma_t \Vert_{\Fr}^2 = 2 \Tr^2(\Lambda_t) -  \Tr(\Lambda_t^2).
\end{align*}
As before, for each $n \in \{10, 50, 100, 1000\}$, $t \in \cT = \{0, 1/69, \ldots, 68/69, 1\}$, we generate 5000 samples
\begin{align*}
    a_{j,t} = \frac{
        \Vert \widehat \Sigma_{j,t} - \Sigma_t \Vert_{\Fr}^2
    }{
        \E \Vert \widehat \Sigma_{j,t} - \Sigma_t \Vert_{\Fr}^2
    }, \quad j = 1, \ldots, 5000, 
\end{align*}
where $\widehat \Sigma_{j,t}$ is the empirtical covariance matrix based on $\bX_1^{j,t}, \ldots, \bX_n^{j,t}$. Finally, we plot the dependence of $0.95$-confidence intervals of $a_{j,t}$ on $\ttr(\Sigma_t)$. The results are displayed in Figure~\ref{fig:counterexample}. Clearly, the width of confidence intervals does not tend to zero as $\ttr(\Sigma)$ increases.

\section{Proofs}
\label{sec:proofs}

This section contains proofs of the results presented in Section \ref{sec:main_results}.

\subsection{Proof of Proposition \ref{prop:orlicz}}
\label{sec:prop_orlicz_proof}

The key step of our proof is an upper bound the Orlicz norm of $(\bxi^\top V \bxi - \Tr(V))$ with respect to the measure $\sfP_U$, defined in \eqref{eq:pu} and \eqref{eq:pu_expectation}. Similarly to the standard Orlicz $\psi_1$-norm, we define
\[
    \|f(\bxi)\|_{\psi_1(\sfP_U)} = \inf\left\{t > 0 : \sfP_U e^{|f(\bxi)| / t} \leq 2 \right\}
\]
for any function $f : \R^d \rightarrow \R$ and any $U \in \R^{d \times d}$ such that $\E \exp\{\bxi^\top U \bxi\} < \infty$. The $\psi_1(\sfP_U)$-norm is just the Orlicz norm with respect to a different probability measure, and hence, it inherits all the properties of the usual $\psi_1$-norm. For instance, analogously to the bound on moments of sub-exponential random variables (see, e.g., \citep[proof of Lemma 2, eq. (10)]{zhivotovskiy21}), we have
\begin{equation}
    \label{eq:sub-exp_moments}
    \sfP_U |f(\bxi)|^k \leq 2 \, \Gamma(k + 1) \, \|f(\bxi)\|_{\psi_1(\sfP_U)}^k
    \quad \text{for all $k \in \mathbb N$.}
\end{equation}
If we manage to obtain an upper bound on
\[
    \left\| \bxi^\top V \bxi - \Tr(V) \right\|_{\psi_1(\sfP_U)},
\]
then it is straightforward to confine
\[
    \sfP_U \bxi^\top V \bxi - \Tr(V)
    \quad \text{and} \quad
    \sfP_U \left(\bxi^\top V \bxi - \Tr(V) \right)^4.
\]
After that, the claim of the proposition follows immediately from the triangle inequality:
\[
    \left[ \sfP_U \left(\bxi^\top V \bxi - \sfP_U \bxi^\top V \bxi \right)^4 \right]^{1/4}
    \leq \left[ \sfP_U \left(\bxi^\top V \bxi - \Tr(V) \right)^4 \right]^{1/4} + \left| \sfP_U \bxi^\top V \bxi - \Tr(V) \right|.
\]
For convenience, we split the rest of the proof into several steps.

\bigskip

\noindent\textbf{Step 1: a bound on the Orlicz norm.}\quad
We are going to show that
\begin{equation}
    \label{eq:orlicz_norm_bound}
    \left\| \bxi^\top V \bxi - \Tr(V) \right\|_{\psi_1(\sfP_U)}
    \leq \frac{2\omega \|V\|_{\Fr}}{1 - \omega \|U\|_{\Fr}}.
\end{equation}
For this purpose, it is enough to prove that
\begin{equation}
    \label{eq:t}
    \sfP_U \exp\left\{ \frac{|\bxi^\top V \bxi - \Tr(V)|}t \right\} \leq 2,
    \quad \text{where} \quad
    t = \frac{2\omega \|V\|_{\Fr}}{1 - \omega \|U\|_{\Fr}}.
\end{equation}
According to Jensen's inequality, we have
\[
    \E e^{\bxi^\top U \bxi} \geq e^{\E \bxi^\top U \bxi} = e^{\Tr(U)}.
\]
Then it holds that
\begin{align*}
    \sfP_U \exp\left\{ \frac{|\bxi^\top V \bxi - \Tr(V)|}t \right\}
    &
    \leq \left( \sfP_U \exp\left\{ \frac{2 |\bxi^\top V \bxi - \Tr(V)|}t \right\} \right)^{1/2}
    \\&
    \leq \left( \E \exp\left\{ \frac{2 |\bxi^\top V \bxi - \Tr(V)|}t + \bxi^\top U \bxi - \Tr(U) \right\}  \right)^{1/2}.
\end{align*}
Since $|\bxi^\top V \bxi - \Tr(V)| = \max\{\bxi^\top V \bxi - \Tr(V), \Tr(V) - \bxi^\top V \bxi\}$, the expression in the right-hand side is not greater than
\[
    \left( \sum\limits_{\eps \in \{-1, 1\}} \exp\left\{ \bxi^\top \left( U + 2 \eps V / t \right) \bxi - \Tr(U + 2 \eps V / t) \right\} \right)^{1/2}.
\]
Note that, due to the triangle inequality and \eqref{eq:orlicz_norm_inequality}, for any $\eps \in \{-1, 1\}$, we have
\begin{align*}
    \left\| \bxi^\top \left( U + 2 \eps V / t \right) \bxi - \Tr(U + 2 \eps V / t) \right\|_{\psi_1}
    &
    \leq \left\| \bxi^\top U \bxi - \Tr(U) \right\|_{\psi_1}
    + \frac{2}t \left\| \bxi^\top V \bxi - \Tr(V) \right\|_{\psi_1}
    \\&
    \leq \omega \|U\|_{\Fr} + \frac{2 \omega \|V\|_{\Fr}}t
    = 1,
\end{align*}
where $\|\cdot\|_{\psi_1}$ stands for the standard $\psi_1$-norm (see the notations in Section \ref{sec:introduction}).
Then the H\"{o}lder inequality implies that
\[
    \E \exp\left\{ \bxi^\top \left( U + 2 \eps V / t \right) \bxi - \Tr(U + 2 \eps V / t) \right\}
    \leq 2^{\left\| \bxi^\top \left( U + 2 \eps V / t \right) \bxi - \Tr(U + 2 \eps V / t) \right\|_{\psi_1}}
    \leq 2.
\]
Hence, we obtain that
\begin{align*}
    \sfP_U \exp\left\{ \frac{|\bxi^\top V \bxi - \Tr(V)|}t \right\}
    &
    \leq \left( \sum\limits_{\eps \in \{-1, 1\}} \exp\left\{ \bxi^\top \left( U + 2 \eps V / t \right) \bxi - \Tr(U + 2 \eps V / t) \right\} \right)^{1/2}
    \\&
    \leq \left( \sum\limits_{\eps \in \{-1, 1\}} 2 \right)^{1/2} = 2,
\end{align*}
which yields \eqref{eq:orlicz_norm_bound}.

\bigskip

\noindent\textbf{Step 2: bounds on $\sfP_U \left(\bxi^\top V \bxi - \Tr(V) \right)^4$ and $\left|\sfP_U \bxi^\top V \bxi - \Tr(V)\right|$.}\quad
Applying the inequality \eqref{eq:sub-exp_moments}, we obtain that
\begin{equation}
    \label{eq:pu_fourth_moment_bound}
    \sfP_U \left(\bxi^\top V \bxi - \Tr(V) \right)^4
    \leq 48 \left\|\bxi^\top V \bxi - \Tr(V) \right\|_{\psi_1(\sfP_U)}^4
    \leq 48 \left( \frac{2\omega \|V\|_{\Fr}}{1 - \omega \|U\|_{\Fr}} \right)^4.
\end{equation}
Concerning $\left|\sfP_U \bxi^\top V \bxi - \Tr(V)\right|$, we will prove a tighter bound than in \eqref{eq:sub-exp_moments}. Let $t > 0$ be as defined in \eqref{eq:t}. Then it holds that
\begin{equation}
    \label{eq:pu_expectation_bound}
    \left|\sfP_U \bxi^\top V \bxi - \Tr(V)\right|
    \leq t \log \E \exp\left\{ \frac{\left|\sfP_U \bxi^\top V \bxi - \Tr(V)\right|}t \right\}
    \leq t \log 2
    = \frac{2 \omega \|V\|_{\Fr}}{1 - \omega \|U\|_{\Fr}} \cdot \log 2.
\end{equation}

\bigskip

\noindent\textbf{Step 3: final bound.}\quad
Summing up the inequalities \eqref{eq:pu_fourth_moment_bound} and \eqref{eq:pu_expectation_bound} and using the triangle inequality, we obtain that
\[
    \left[ \sfP_U \left(\bxi^\top V \bxi - \sfP_U \bxi^\top V \bxi \right)^4 \right]^{1/4}
    \leq \frac{2 \omega \|V\|_{\Fr}}{1 - \omega \|U\|_{\Fr}} \left( \log 2 + 48^{1/4} \right).
\]
If $\omega \|U\|_{\Fr} \leq 1/6$, then
\[
    \left[ \sfP_U \left(\bxi^\top V \bxi - \sfP_U \bxi^\top V \bxi \right)^4 \right]^{1/4}
    \leq \frac{12}5 \left( \log 2 + 48^{1/4} \right) \omega \|V\|_{\Fr} 
    \leq 8 \omega \|V\|_{\Fr}.
\]
This implies that, for any $U \in \R^{d \times d}$ satisfying the inequality $6 \omega \|U\|_{\Fr} \leq 1$, it holds that
\[
    \sfP_U \left(\bxi^\top V \bxi - \sfP_U \bxi^\top V \bxi \right)^4 \leq \left(64 \omega^2 \right)^2 \|V\|_{\Fr}^4.
\]
\myendproof

\subsection{Proof of Proposition \ref{prop:counterexample}}
\label{sec:prop_counterexample_proof}
Let $\beps = (\eps_1, \dots, \eps_d)^\top$ be a random vector with i.i.d. Rademacher entries and let $\eta \sim \Be(1/2)$ be independent of $\beps$. We are going to prove that a random vector $\bxi = \sqrt{2} \eta \, \beps$ satisfies Proposition \ref{prop:counterexample}. First, it is easy to observe that
\[
    \E \bxi = \sqrt{2} \E \eta \, \E \beps = \bzero
    \quad \text{and} \quad
    \E \bxi \bxi^\top = 2 \E \eta^2 \, \E \beps \beps^\top = I_d.
\]
Second, let us fix any $\bu \in \R^d$ and show that
\[
    \E \exp\left\{ \frac{3(\bxi^\top \bu)^2}{16 \|\bu\|^2} \right\} \leq 2,
\]
or equivalently, that the $\psi_2$-norm of $\bxi^\top \bu$ does not exceed $16 \|\bu\|^2 / 3 = 16 \E (\bxi^\top \bu)^2 / 3$. Indeed, introducing a random variable $\gamma \sim \cN(0, 1)$, which is independent of $\eta$ and $\beps$, we obtain that
\[
    \E \exp\left\{ \frac{3(\bxi^\top \bu)^2}{16 \|\bu\|^2} \right\}
    = \E \exp\left\{ \frac{3 \eta^2 (\beps^\top \bu)^2}{8 \|\bu\|^2} \right\}
    \leq \E \exp\left\{ \frac{3 (\beps^\top \bu)^2}{8 \|\bu\|^2} \right\}
    = \E \exp\left\{ \frac{\sqrt{3} \gamma \beps^\top \bu}{2 \|\bu\|} \right\}.
\]
Hoeffding's lemma, applied to independent Rademacher random variables $\eps_1, \dots, \eps_d$, implies that
\[
    \E \exp\left\{ \frac{\sqrt{3} \gamma \beps^\top \bu}{2 \|\bu\|} \right\}
    \leq \E \prod\limits_{i = 1}^d \exp\left\{ \frac{3 \gamma^2 u_i^2}{8 \|\bu\|^2} \right\}
    = \E \exp\left\{ \frac{3 \gamma^2}{8} \right\}
    = \left(1 - \frac34 \right)^{-1/2}
    = 2.
\]
Thus, we proved that $\left\| \bxi^\top \bu \right\|_{\psi_2}^2 \leq 16 \, \E \| \bxi^\top \bu \|^2 / 3$. On the other hand, for $V = I_d$, we have
\begin{align*}
    \E \left(\bxi^\top V \bxi - \E \bxi^\top V \bxi \right)^4
    &
    = \E \left(\|\bxi\|^2 - d \right)^4
    = \E \left(2 \eta^2 \|\beps\|^2 - d \right)^4
    = \E \left(2 \eta^2 d - d \right)^4
    \\&
    = \E \left(2 \eta^2 - 1 \right)^4 d^4
    = d^4
    = d^2 \left\|V \right\|_{\Fr}^2.
\end{align*}
The proof is finished.

\subsection{Proof of Theorem \ref{th:fr_norm_upper_bound}}
\label{sec:th_fr_norm_upper_bound_proof}

This section is devoted to the proof of Theorem \ref{th:fr_norm_upper_bound}. Let us denote
\[
    H = \sqrt n \; \Sigma^{-1/2} \left( \widehat \Sigma - \Sigma \right) \Sigma^{-1/2}
    = \frac1{\sqrt n} \sum\limits_{i = 1}^n \left( (\Sigma^{-1/2} \bX_i) (\Sigma^{-1/2} \bX_i)^\top - I \right)
    = \frac1{\sqrt n} \sum\limits_{i = 1}^n \left( \bxi_i \bxi_i^\top - I \right).
\]
Then
\[
    \left\| \widehat \Sigma - \Sigma \right\|_{\Fr}^2 = \frac1n \left\|\Sigma^{1/2} H \Sigma^{1/2} \right\|_{\Fr}^2,
\]
and it is enough to show that
\begin{align*}
    &
    \left\|\Sigma^{1/2} H \Sigma^{1/2} \right\|_{\Fr}^2 - \E \left\|\Sigma^{1/2} H \Sigma^{1/2} \right\|_{\Fr}^2
    \\&
    \leq 4 \|\Sigma\|^2 \max\left\{ \sqrt{2 \big( \tau \ttr(\Sigma^2) + \tau^2 \ttr(\Sigma^4) \big) \log(2 / \delta)}, \; 4e \tau \log(2 / \delta) \right\}
\end{align*}
with probability at least $(1 - \delta)$. In order to do this, we first prove that
\begin{equation}
    \label{eq:H_large_deviation}
    \left\|\Sigma^{1/2} H \Sigma^{1/2} \right\|_{\Fr}^2 - \E \left\|\Sigma^{1/2} H \Sigma^{1/2} \right\|_{\Fr}^2 \leq \max\left\{ 4 \ttv \sqrt{2 \log(2 / \delta)}, \; 16e (\kappa \vee \|\Sigma\|^2) \log(2 / \delta) \right\}
\end{equation}
with high probability, where
\begin{align}
    \label{eq:v_kappa}
    \ttv^2
    &\notag
    = \left\| \E \rmvec(\bX \bX^\top - \Sigma) \rmvec(\bX \bX^\top - \Sigma)^\top \right\|_{\Fr}^2,
    \\
    \kappa
    &
    = \sup\limits_{\|U\|_{\Fr} = 1} \E \left[ \bX^\top U \bX - \Tr(U^\top \Sigma) \right]^2,
\end{align}
and then use Lemma \ref{lem:variance} to derive the final bound.
In the rest of the proof we study exponential moments of $\|\Sigma^{1/2} H \Sigma^{1/2}\|_{\Fr}^2$ and derive the large deviation bound \eqref{eq:H_large_deviation}.
Since the proof is quite technical, we split it into several steps to improve readability of the paper. The proofs of some auxiliary results are moved to Appendix.

\medskip
\noindent{\bf Step 1: linearization.}\quad
We start with the linearization trick. Let $\Gamma \in \R^{d \times d}$ be a random matrix with i.i.d. standard Gaussian entries. Then it is straightforward to check that
\begin{equation}
    \label{eq:standard_linearization}
    \E \exp\left\{\lambda \|\Sigma^{1/2} H \Sigma^{1/2}\|_{\Fr}^2 / 2 \right\}
    = \E \exp\left\{\sqrt{\lambda} \Tr(\Gamma^\top \Sigma^{1/2} H \Sigma^{1/2}) \right\}.
\end{equation}
The idea is that the power of the exponent in the right-hand side is linear in $H$, and the analysis of linear statistics is more convenient.
Unfortunately, the linearization trick in the form \eqref{eq:standard_linearization} will not bring us to the desired result, because the elements of $H$ are sub-exponential, and hence, the conditional moment
\[
    \E_H \exp\left\{\sqrt{\lambda} \Tr(\Gamma^\top \Sigma^{1/2} H \Sigma^{1/2}) \right\}
    \equiv \E \left[ \exp\left\{\sqrt{\lambda} \Tr(\Gamma^\top \Sigma^{1/2} H \Sigma^{1/2}) \right\} \,\Big\vert\, \Gamma \right]
\]
exists not for all $\Gamma$. Nevertheless, the idea of linearization is still useful, and we just have to tailor the right-hand side of \eqref{eq:standard_linearization} for our setup. This brings us to the following lemma.

\begin{Lem}
    \label{lem:linearization}
    Let $\z$ and $\lambda$ be any positive numbers, and let $\Gamma$ be a random matrix with i.i.d. standard Gaussian entries. Then, for any matrix $A \in \R^{d \times d}$ satisfying the inequality
    \[
        \|\Sigma A \Sigma\|_{\Fr} \leq \z,
    \]
    and for any $\beta > 1$, it holds that
    \begin{align*}
        &
        \exp\left\{\lambda \|\Sigma^{1/2} A \Sigma^{1/2}\|_{\Fr}^2 / 2 \right\}
        \\&
        \leq \frac{\beta}{\beta - 1} \E_\Gamma \, \exp\left\{\sqrt{\lambda} \Tr(\Gamma^\top \Sigma^{1/2} A \Sigma^{1/2}) \right\} \1\left( \|\Sigma^{1/2} \Gamma \Sigma^{1/2}\|_{\Fr} \leq \z \sqrt{\lambda} +  \sqrt{\beta} \Tr(\Sigma) \right).
    \end{align*}
\end{Lem}

The proof of Lemma \ref{lem:linearization} is moved to Appendix \ref{sec:lem_linearization_proof}.
For any deterministic matrix $U \in \R^{d \times d}$, let us denote 
\begin{equation}
    \label{eq:Phi}
    \Phi(U) = \log \E e^{\Tr(H^\top U)} = \sum\limits_{i = 1}^n \log \E e^{\bxi_i^\top U \bxi_i / \sqrt{n}}.
\end{equation}
Applying Lemma \ref{lem:linearization} with $\beta = 2$ and using the Fubini theorem, we obtain that
\begin{align*}
    &
    \E e^{\lambda \|\Sigma^{1/2} H \Sigma^{1/2}\|_{\Fr}^2 / 2} \1\left( \|\Sigma H \Sigma\|_{\Fr} \leq \z \right)
    \\&
    \leq 2 \; \E e^{\sqrt{\lambda} \Tr(\Gamma^\top \Sigma^{1/2} H \Sigma^{1/2})} \1\left( \|\Sigma^{1/2} \Gamma \Sigma^{1/2}\|_{\Fr} \leq \z \sqrt{\lambda} + \sqrt{2} \, \Tr(\Sigma) \right)
    \\&
    = 2 \; \E e^{\Phi(\sqrt{\lambda} \Sigma^{1/2} \Gamma \Sigma^{1/2})} \1\left( \|\Sigma^{1/2} \Gamma \Sigma^{1/2}\|_{\Fr} \leq \z \sqrt{\lambda} + \sqrt{2} \, \Tr(\Sigma) \right).
\end{align*}

\medskip
\noindent{\bf Step 2: Taylor's expansion.}\quad
It is easy to observe that the function $\Phi$ from \eqref{eq:Phi} can be expressed through the function $\varphi$, defined in \eqref{eq:phi}:
\begin{equation}
    \label{eq:Phi_phi}
    \Phi(U) = \sum\limits_{i = 1}^n \log \E e^{\bxi_i^\top U \bxi_i / \sqrt{n}} = n \varphi(U / \sqrt n).
\end{equation}
On this step, we use the smoothness of $\varphi(U)$, guaranteed by Assumption \ref{as:derivatives}, to derive an upper bound on
\[
    \E e^{\Phi(\sqrt{\lambda} \Sigma^{1/2} \Gamma \Sigma^{1/2})} \1\left( \|\Sigma^{1/2} \Gamma \Sigma^{1/2}\|_{\Fr} \leq \z \sqrt{\lambda} +  \sqrt{2} \, \Tr(\Sigma) \right)
\]
for a fixed $\z > 0$. The next lemma helps to bound the derivatives of $\varphi$.
\begin{Lem}
    \label{lem:derivatives}
    Suppose that the random vector $\bxi = \Sigma^{-1/2} \bX$ satisfies Assumption \ref{as:derivatives}. Then, for any $U \in \R^{d \times d}$, such that $\|U\|_{\Fr} \leq \rho_{\max}$, the derivatives of the cumulant generating function
    \[
        \varphi(U) = \log \E \exp\left\{\bxi^\top U \bxi \right\}
    \]
    satisfy the inequalities
    \[
        0 \leq
        \left\langle \nabla^2 \varphi(U), V^{\otimes 2} \right\rangle
        = \sfP_U \left(\bxi^\top V \bxi - \sfP_U \bxi^\top V \bxi \right)^2
        \leq \tau \|V\|_{\Fr}^2,
    \]
    \[
        - \tau^{3/2} \|V\|_{\Fr}^3
        \leq \left\langle \nabla^3 \varphi(U), V^{\otimes 3} \right\rangle
        = \sfP_U \left(\bxi^\top V \bxi - \sfP_U \bxi^\top V \bxi \right)^3
        \leq \tau^{3/2} \|V\|_{\Fr}^3,
    \]
    and
    \[
        -2 \tau^2 \|V\|_{\Fr}^4
        \leq \left\langle \nabla^4 \varphi(U), V^{\otimes 4} \right\rangle
        \leq \tau^2 \|V\|_{\Fr}^4.
    \]
\end{Lem}

\begin{Rem}
    Throughout the paper, we use the standard notation
    \[
        \left\langle \nabla \varphi(U), V \right\rangle \equiv \Tr \left(\nabla \varphi(U)^\top V \right),
    \]
    where $\nabla \varphi(U)$ is the gradient of the scalar function $\varphi : \R^{d \times d} \rightarrow \R$ with respect to the matrix $U \in \R^{d \times d}$.
    For $k \geq 2$, the notation $\left\langle \nabla^k \varphi(U), V^{\otimes k} \right\rangle$ is defined recursively:
    \[
        \left\langle \nabla^k \varphi(U), V^{\otimes k} \right\rangle
        \equiv \Big\langle \nabla \left\langle \nabla^{k-1} \varphi(U), V^{\otimes (k - 1)} \right\rangle, V \Big\rangle.
    \]
\end{Rem}
The proof of Lemma \ref{lem:derivatives} is moved to Appendix \ref{sec:derivatives}. Due to \eqref{eq:Phi_phi}, we immediately obtain that
\begin{align}
    \label{eq:Phi_derivatives}
    \left\langle \nabla^2 \Phi(O), U^{\otimes 2} \right\rangle
    &\notag
    = \left\langle \nabla^2 \Phi(O), U^{\otimes 2} \right\rangle
    = \E \left(\bxi^\top U \bxi - \E \bxi^\top U \bxi \right)
    = \E \left[\Tr\left( (\bxi\bxi^\top - I) U \right) \right]^2,
    \\
    \left| \left\langle \nabla^3 \Phi(O), U^{\otimes 3} \right\rangle \right|
    &
    = \frac1{\sqrt{n}} \left| \left\langle \nabla^3 \varphi(O), U^{\otimes 3} \right\rangle \right|
    = \frac1{\sqrt{n}} \E \left[\Tr\left( (\bxi\bxi^\top - I) U \right) \right]^3
    \leq \frac{\tau^{3/2}}{\sqrt n} \|U\|_{\Fr}^3,
    \\
    \left| \left\langle \nabla^4 \Phi(V), U^{\otimes 4} \right\rangle \right|
    &\notag
    = \frac1{n} \left| \left\langle \nabla^4 \varphi(V), U^{\otimes 4} \right\rangle \right|
    \leq \frac{\tau^2}{n} \|U\|_{\Fr}^4
    \quad \text{for all $\|V\|_{\Fr} \leq \rho_{\max}$.}
\end{align}
Here $O \in \R^{d \times d}$ denotes the matrix with zero entries.
Let us introduce
\begin{equation}
    \label{eq:m}
    \ttm = n \; \E \left\|\widehat \Sigma - \Sigma \right\|_{\Fr}^2 = \E \left\|\Sigma^{1/2} H \Sigma^{1/2} \right\|_{\Fr}^2
\end{equation}
for brevity.
Using the inequalities \eqref{eq:Phi_derivatives}, we prove the following result.

\begin{Lem}
    \label{lem:restricted_moment}
    Assume \ref{as:derivatives} and let $\lambda$ and $\z$ be any positive numbers satisfying the inequalities
    \[
        \z \sqrt{\lambda} +  \sqrt{2} \, \Tr(\Sigma) \leq \rho_{\max} \sqrt{n} 
    \]
    and
    \[
        \lambda \kappa + G(\lambda, \z) \|\Sigma\|^2 \leq \frac12,
    \]
    where
    \begin{equation}
        \label{eq:G}
        G(\lambda, \z)
        = \frac{\lambda \big( \z \lambda + \sqrt{2\lambda} \, \Tr(\Sigma) \big)^2}{36 n} \left( \tau^3 \big(\z \lambda +  \sqrt{2 \lambda} \, \Tr(\Sigma) \big)^2 + 3 \tau^2 \right).
    \end{equation}
    Then it holds that
    \begin{align*}
        &
        \E e^{\Phi(\sqrt{\lambda} \Sigma^{1/2} \Gamma \Sigma^{1/2})} \1\left( \|\Sigma^{1/2} \Gamma \Sigma^{1/2}\|_{\Fr} \leq \z \sqrt{\lambda} + \sqrt{2} \, \Tr(\Sigma) \right)
        \\&
        \leq \exp\left\{ \frac{\lambda \ttm}2 + 2 \lambda^2 \ttv^2 + \frac12 G(\lambda, \z) \Tr(\Sigma)^2 + 2 G(\lambda, \z)^2 \left\| \Sigma \right\|_{\Fr}^4 \right\},
    \end{align*}
    where $\ttm$, $\ttv$ and $\kappa$ are defined in \eqref{eq:m} and \eqref{eq:v_kappa}.
\end{Lem}
The proof of Lemma \ref{lem:restricted_moment} is deferred to Appendix \ref{sec:lem_restricted_moment_proof}.

\medskip
\noindent{\bf Step 3: peeling argument.}\quad
On this step, we transform the restricted exponential moment bound from Lemma \ref{lem:restricted_moment} into a large deviation bound on $\|\Sigma^{1/2} H \Sigma^{1/2}\|_{\Fr}$.
Let $\z$ and $\lambda$ be positive numbers to be specified later. Then it holds that
\begin{align*}
    \p\left( \| \Sigma^{1/2} H \Sigma^{1/2} \|_{\Fr} \geq \z \right)
    &
    = \sum\limits_{k = 1}^\infty \p\left( e^{k - 1} \z \leq \| \Sigma^{1/2} H \Sigma^{1/2} \|_{\Fr} < e^k \z \right)
    \\&
    \leq \sum\limits_{k = 1}^\infty \p\left( \| \Sigma^{1/2} H \Sigma^{1/2} \|_{\Fr} \geq e^{k-1} \z \text{ and } 
    \| \Sigma H \Sigma\|_{\Fr}
    \leq e^k \z \|\Sigma\| \right).
\end{align*}
Obviously, for any $k \geq 1$, on the event $\{\| \Sigma^{1/2} H \Sigma^{1/2} \|_{\Fr} \geq e^{k-1} \z \}$, we have
\[
    \1\left( \| \Sigma^{1/2} H \Sigma^{1/2} \|_{\Fr} \geq e^{k-1} \z \right) \leq \exp\left\{ -\frac{e^{2k-2} \lambda \z^2}{2 \cdot e^k} + \frac{\lambda \| \Sigma^{1/2} H \Sigma^{1/2} \|_{\Fr}^2}{2 \cdot e^k}\right\}
\]
This yields the following version of the Markov inequality:
\begin{align*}
    &
    \sum\limits_{k = 1}^\infty \p\left( \| \Sigma^{1/2} H \Sigma^{1/2} \|_{\Fr} \geq e^{k-1} \z \text{ and } 
    \| \Sigma H \Sigma\|_{\Fr}
    \leq e^k \z \|\Sigma\| \right)
    \\&
    \leq \sum\limits_{k = 1}^\infty \E \exp\left\{ -\frac{e^{k-2} \lambda \z^2}2 + \frac{\lambda \| \Sigma^{1/2} H \Sigma^{1/2} \|_{\Fr}^2}{2 \cdot e^k} \right\} \1 \big( \| \Sigma H \Sigma\|_{\Fr}
    \leq e^k \z \|\Sigma\|\big).
\end{align*}
Lemma \ref{lem:linearization} and Lemma \ref{lem:restricted_moment} imply that
\begin{align*}
    &
    \sum\limits_{k = 1}^\infty \E \exp\left\{ -\frac{e^{k-2}  \lambda \z^2}2 + \frac{\lambda \| \Sigma^{1/2} H \Sigma^{1/2} \|_{\Fr}^2}{2 \cdot e^k} \right\} \1 \big( \| \Sigma H \Sigma\|_{\Fr}
    \leq e^k \z \|\Sigma\|\big)
    \\&
    \leq 2 \sum\limits_{k = 1}^\infty \exp\left\{ -\frac{e^{k-2} \lambda \z^2}2 + \frac{\lambda \ttm}{2 \cdot e^k} + \frac{2 \lambda^2 \ttv^2}{e^{2k}} + \frac12 G\left(e^{-k} \lambda, e^k \z \|\Sigma\| \right) \Tr(\Sigma)^2 + 2 G\left(e^{-k} \lambda, e^k \z \|\Sigma\| \right)^2 \left\| \Sigma \right\|_{\Fr}^4 \right\}.
\end{align*}
It is straighforward to check that, by the definition of $G(\lambda, \z)$ (see \eqref{eq:G}), we have
\[
    G\left(e^{-k} \lambda, e^k \z \|\Sigma\| \right)
    \leq G\left(\lambda, \z \|\Sigma\| \right).
\]
Hence,
\begin{align*}
    &
    2 \sum\limits_{k = 1}^\infty \exp\left\{ -\frac{e^{k-2} \lambda \z^2}2 + \frac{\lambda \ttm}{2 \cdot e^k} + \frac{2 \lambda^2 \ttv^2}{e^{2k}} + \frac12 G\left(e^{-k} \lambda, e^k \z \|\Sigma\| \right) \Tr(\Sigma)^2 + 2 G\left(e^{-k} \lambda, e^k \z \|\Sigma\| \right)^2 \left\| \Sigma \right\|_{\Fr}^4 \right\}
    \\&
    \leq 2 \exp\left\{ -\frac{\lambda (\z^2 - \ttm)}{2e} + \frac{2 \lambda^2 \ttv^2}{e^2} + \frac12 G\left(\lambda, \z \|\Sigma\| \right) \Tr(\Sigma)^2 + 2 G\left(\lambda, \z \|\Sigma\| \right)^2 \left\| \Sigma \right\|_{\Fr}^4 \right\}
    \\&\quad
    + 2 \sum\limits_{k = 2}^\infty \exp\left\{ -\frac{(e^{k-1} - 1) \lambda \z^2}{2e} -\frac{\lambda (\z^2 - \ttm)}{2e} + \frac{2 \lambda^2 \ttv^2}{e^2} + \frac12 G\left(\lambda, \z \|\Sigma\| \right) \Tr(\Sigma)^2 + 2 G\left(\lambda, \z \|\Sigma\| \right)^2 \left\| \Sigma \right\|_{\Fr}^4 \right\}.
\end{align*}
Let us introduce a function $g : \R_+ \rightarrow \R_+$, defined by
\begin{equation}
    \label{eq:g}
    g(u)
    = 1 + \sum\limits_{k=1}^\infty \exp\left\{ -\frac{(e^k - 1) u}{2e} \right\}
    = \sum\limits_{k=0}^\infty \exp\left\{ -\frac{(e^k - 1) u}{2e} \right\}.
\end{equation}
Then
\begin{align*}
    &
    2 \sum\limits_{k = 1}^\infty \exp\left\{ -\frac{e^{k-2} \lambda \z^2}2 + \frac{\lambda \ttm}{2 \cdot e^k} + \frac{2 \lambda^2 \ttv^2}{e^{2k}} + \frac12 G\left(e^{-k} \lambda, e^k \z \|\Sigma\| \right) \Tr(\Sigma)^2 + 2 G\left(e^{-k} \lambda, e^k \z \|\Sigma\| \right)^2 \left\| \Sigma \right\|_{\Fr}^4 \right\}
    \\&
    \leq 2 g(\lambda \z^2) \exp\left\{ -\frac{\lambda (\z^2 - \ttm)}{2e} + \frac{2 \lambda^2 \ttv^2}{e^2} + \frac12 G\left(\lambda, \z \|\Sigma\| \right) \Tr(\Sigma)^2 + 2 G\left(\lambda, \z \|\Sigma\| \right)^2 \left\| \Sigma \right\|_{\Fr}^4 \right\},
\end{align*}
and thus,
\begin{align}
    \label{eq:large_deviation_bound}
    &\notag
    \p\left( \| \Sigma^{1/2} H \Sigma^{1/2} \|_{\Fr} \geq \z \right)
    \\&
    \leq 2 g(\lambda \z^2) \exp\left\{ -\frac{\lambda (\z^2 - \ttm)}{2e} + \frac{2 \lambda^2 \ttv^2}{e^2} + \frac12 G\left(\lambda, \z \|\Sigma\| \right) \Tr(\Sigma)^2 + 2 G\left(\lambda, \z \|\Sigma\| \right)^2 \left\| \Sigma \right\|_{\Fr}^4 \right\}.
\end{align}

\medskip

\noindent{\bf Step 4: choosing $\lambda$ and $\z$.}\quad
On this step, we specify $\lambda$ and $\z$ and ensure that they satisfy the conditions of Lemma \ref{lem:restricted_moment}. Let us take
\begin{equation}
    \label{eq:lambda_z}
    \lambda = \frac{e(\z^2 - \ttm)}{8 \ttv^2} \land \frac1{4 (\kappa \vee \tau \|\Sigma\|^2)}
    \quad \text{and} \quad
    \z^2 = \ttm + \max\left\{ 4 \ttv \sqrt{2 \log(2 / \delta)}, 16e (\kappa \vee \tau \|\Sigma\|^2) \log \frac{3}\delta \right\}.
\end{equation}
Such $\lambda$ minimizes the expression
\[
    -\frac{\lambda (\z^2 - \ttm)}{2 e} + \frac{2 \lambda^2 \ttv^2}{e^2}
    \quad \text{over} \quad
    \lambda \in \left[ 0, \frac1{4 (\kappa \vee \tau \|\Sigma\|^2)} \right].
\]
With $\lambda$ and $\z$, defined in \eqref{eq:lambda_z}, we have
\begin{equation}
    \label{eq:lambda_substitution}
    -\frac{\lambda (\z^2 - \ttm)}{2e} + \frac{2 \lambda^2 \ttv^2}{e^2}
    = -\min \left\{ \frac{(\z^2 - \ttm)^2}{32 \ttv^2}, \frac{\z^2 - \ttm}{16 e (\kappa \vee \tau \|\Sigma\|^2)} \right\}
\end{equation}
and
\[
    \z^2 \lambda
    = \ttm \lambda + (\z^2 - \ttm) \lambda
    = \ttm \lambda + 4e \log(3 / \delta) \geq 4e \log 3.
\]
As a consequence, we obtain that
\begin{equation}
    \label{eq:g_value}
    2 g(\z^2 \lambda) \leq 2 g\left( 4e \log 3 \right) < 2.1,
\end{equation}
because, according to the definition of $g$, it is a non-increasing function on $(0, +\infty)$.

Let us show that $\lambda$ and $\z \|\Sigma\|$ satisfy the conditions of Lemma \ref{lem:restricted_moment}, that is,
\begin{equation}
    \label{eq:admissible_z_lambda}
    \z \|\Sigma\| \sqrt{\lambda} +  \sqrt{2} \, \Tr(\Sigma) \leq \rho_{\max} \sqrt{n} 
\end{equation}
and
\begin{equation}
    \label{eq:admissible_lambda}
    \lambda \kappa + G(\lambda, \z \|\Sigma\|) \|\Sigma\|^2 \leq \frac12,
\end{equation}
where the function $G$ is defined in \eqref{eq:G}.
Since $\ttm \lambda \leq 0.25 \, \ttm / (\kappa \vee \tau \|\Sigma\|^2)$, we have
\begin{align*}
    \z \sqrt{\lambda} \|\Sigma\| + \sqrt{2} \, \Tr(\Sigma)
    &
    \leq \|\Sigma\| \left( \sqrt{\frac{\ttm}{4 (\kappa \vee \tau \|\Sigma\|^2)} + 4e \log(3 / \delta)} + \sqrt{2} \, \ttr(\Sigma) \right)
    \\&
    \leq \|\Sigma\| \left( \frac12 \sqrt{\frac{\ttm}{\tau \|\Sigma\|^2}} + 2 \sqrt{e \log(3 / \delta)} + \sqrt{2} \, \ttr(\Sigma) \right).
\end{align*}
The following lemma allows us to simplify the expression in the right-hand side.
\begin{Lem}
    \label{lem:mean}
    Suppose that Assumption \ref{as:fourth_moment} holds with some $\alpha > 0$. Then $\ttm = n \; \E \|\widehat \Sigma - \Sigma \|_{\Fr}^2$ satisfies the inequalities
    \[
        \left( \Tr(\Sigma) \right)^2 - \Tr(\Sigma^2) \leq \ttm \leq \left( \Tr(\Sigma) \right)^2 + (\alpha  - 1) \; \Tr(\Sigma^2).
    \]
\end{Lem}
A reader can find the proof of Lemma \ref{lem:mean} in Appendix \ref{sec:lem_mean_proof}. We would like to note that Assumption \ref{as:fourth_moment} is weaker than Assumption \ref{as:derivatives} and, under the conditions of Theorem \ref{th:fr_norm_upper_bound}, we have that Assumption \ref{as:fourth_moment} is fulfilled with $\alpha = \tau$, which yields
\[
    \ttm \leq \left( \Tr(\Sigma) \right)^2 + (\tau  - 1) \; \Tr(\Sigma^2).
\]
Applying Lemma \ref{lem:mean} and taking into account that $\tau \geq 2$, we obtain that
\begin{align*}
    \z \sqrt{\lambda} \|\Sigma\| + \sqrt{2} \, \Tr(\Sigma)
    &
    \leq \|\Sigma\| \left( \frac12 \sqrt{\frac{\ttm}{\tau \|\Sigma\|^2}} + 2 \sqrt{e \log(3 / \delta)} + \sqrt{2} \, \ttr(\Sigma) \right)
    \\&
    \leq \|\Sigma\| \left( \frac12 \sqrt{\frac{\ttr(\Sigma)^2}\tau + \ttr(\Sigma^2)} + 2 \sqrt{e \log(3 / \delta)} + \sqrt{2} \, \ttr(\Sigma) \right)
    \\&
    \leq \|\Sigma\| \left( \frac{\ttr(\Sigma)}{2 \sqrt{\tau}} + \frac12 \sqrt{\ttr(\Sigma^2)} + 2 \sqrt{e \log(3 / \delta)} + \sqrt{2} \, \ttr(\Sigma) \right)
    \\&
    \leq \|\Sigma\| \left[ \frac{3\ttr(\Sigma)}{2\sqrt{2}} + \frac12 \sqrt{\ttr(\Sigma^2)} + 2 \sqrt{e \log(3 / \delta)} \right]
    = \|\Sigma\| \; \ttR(\Sigma, \delta),
\end{align*}
and then \eqref{eq:admissible_z_lambda} is fulfilled due to \eqref{eq:n_admissible}.
Moreover, the inequality \eqref{eq:n_admissible} implies that
\begin{align}
    \label{eq:phi_trace_squared}
    G(\lambda, \z \|\Sigma\|) \, \Tr(\Sigma)^2
    &\notag
    = \frac{\lambda^2 \|\Sigma\|^4 \, \ttr(\Sigma)^2 \, \ttR(\Sigma, \delta)^2}{36 n} \left( \tau^3 \lambda \|\Sigma\|^2 \, \ttR(\Sigma, \delta)^2 + 3 \tau^2 \right)
    \\&
    \leq \frac{\|\Sigma\|^4 \, \ttr(\Sigma)^2 \, \ttR(\Sigma, \delta)^2}{144 n (\kappa^2 \vee \tau^2 \|\Sigma\|^4)} \left( \frac{\tau^3 \|\Sigma\|^2 \ttR(\Sigma, \delta)^2}{4 (\kappa \vee \tau \|\Sigma\|^2)} + 3 \tau^2 \right)
    \\&\notag
    \leq \frac{\ttr(\Sigma)^2 \, \ttR(\Sigma, \delta)^2}{144 n} \left( \frac{\ttR(\Sigma, \delta)^2}{4} + 3 \right)
    \leq \frac14.
\end{align}
Similarly, it holds that
\begin{equation}
    \label{eq:phi_frobenius_squared}
    G(\lambda, \z \|\Sigma\|) \, \|\Sigma\|_{\Fr}^2
    = G(\lambda, \z \|\Sigma\|) \, \Tr(\Sigma^2)
    \leq G(\lambda, \z \|\Sigma\|) \, \Tr(\Sigma)^2
    \leq \frac14
\end{equation}
and
\[
    G(\lambda, \z \|\Sigma\|) \, \|\Sigma\|^2 \leq \frac14.  
\]
The last inequality means that \eqref{eq:admissible_lambda} is fulfilled, because
\[
    \lambda \kappa + G(\lambda, \z \|\Sigma\|) \|\Sigma\|^2 
    = \frac{\kappa}{4 (\kappa \vee \tau \|\Sigma\|^2)} + G(\lambda, \z \|\Sigma\|) \|\Sigma\|^2
    \leq \frac14 + \frac14 = \frac12.
\]
Finally, taking into account \eqref{eq:large_deviation_bound}, \eqref{eq:g_value}, \eqref{eq:phi_trace_squared}, and \eqref{eq:phi_frobenius_squared}, we obtain that
\begin{align*}
    &
    \p\left( \| \Sigma^{1/2} H \Sigma^{1/2} \|_{\Fr} \geq \z \right)
    \\&
    \leq 2 g(\lambda \z^2) \exp\left\{ -\frac{\lambda (\z^2 - \ttm)}{2e} + \frac{2 \lambda^2 \ttv^2}{e^2} + \frac12 G\left(\lambda, \z \|\Sigma\| \right) \Tr(\Sigma)^2 + 2 G\left(\lambda, \z \|\Sigma\| \right)^2 \left\| \Sigma \right\|_{\Fr}^4 \right\}
    \\&
    \leq 2.1 \exp\left\{-\frac{\lambda (\z^2 - \ttm)}{2e} + \frac{2 \lambda^2 \ttv^2}{e^2} + \frac12 \cdot \frac14 + 2 \cdot \frac1{16} \right\}
    \\&
    < 3 \exp\left\{-\frac{\lambda (\z^2 - \ttm)}{2e} + \frac{2 \lambda^2 \ttv^2}{e^2} \right\},
\end{align*}
and then the choice of $\lambda$ and $\z$ (see \eqref{eq:lambda_z} and \eqref{eq:lambda_substitution}) ensures that
\begin{align}
    \label{eq:large_deviation_bound_2}
    \p\left( \| \Sigma^{1/2} H \Sigma^{1/2} \|_{\Fr} \geq \z \right)
    &\notag
    < 3 \exp\left\{-\frac{\lambda (\z^2 - \ttm)}{2e} + \frac{2 \lambda^2 \ttv^2}{e^2} \right\}
    \\&
    \leq 3 \exp\left\{-\min \left\{ \frac{(\z^2 - \ttm)^2}{32 \ttv^2}, \frac{\z^2 - \ttm}{16 e (\kappa \vee \tau \|\Sigma\|^2)} \right\}\right\} = \delta.
\end{align}

\medskip

\noindent{\bf Step 5: a bound on $\ttv^2$ and $\kappa$.}\quad
The inequality \eqref{eq:large_deviation_bound_2} yields that
\begin{align*}
    n \left\|\widehat \Sigma - \Sigma \right\|_{\Fr}^2 - n \; \E \left\|\widehat \Sigma - \Sigma \right\|_{\Fr}^2
    &
    = \left\|\Sigma^{1/2} H \Sigma^{1/2} \right\|_{\Fr}^2 - \E \left\|\Sigma^{1/2} H \Sigma^{1/2} \right\|_{\Fr}^2
    \\&
    = \left\|\Sigma^{1/2} H \Sigma^{1/2} \right\|_{\Fr}^2 - \ttm
    \\&
    \leq \max\left\{ 4 \ttv \sqrt{2 \log(2 / \delta)}, \; 16e (\kappa \vee \tau \|\Sigma\|^2) \log(2 / \delta) \right\}
\end{align*}
with probability at least $(1 - \delta)$. Here we used the definition of $\ttm$ (see eq. \eqref{eq:m}). The goal of the final step is to obtain upper bounds on $\ttv$ and $\kappa$, using Assumption \ref{as:derivatives}. In Appendix \ref{sec:lem_variance_proof}, we prove the following result.
\begin{Lem}
    \label{lem:variance}
    Let Assumption \ref{as:fourth_moment} be satisfied with some $\alpha > 0$. Then it holds that
    \[
        \ttv^2
        = \left\| \E \rmvec(\bX \bX^\top - \Sigma) \rmvec(\bX \bX^\top - \Sigma)^\top \right\|_{\Fr}^2
        \leq \alpha \left( \Tr(\Sigma^2) \right)^2 + (\alpha^2 - \alpha) \; \Tr(\Sigma^4)
    \]
    and
    \[
        \kappa = \sup\limits_{\|U\|_{\Fr} = 1} \E \left[ \bX^\top U \bX - \Tr(U^\top \Sigma) \right]^2 \leq \alpha \, \|\Sigma\|^2.
    \]
\end{Lem}
Let us recall that, under the conditions of Theorem \ref{th:fr_norm_upper_bound}, Assumption \ref{as:fourth_moment} is fulfilled with $\alpha = \tau$. Hence, we have
\[
    \ttv^2 \leq \tau \left( \Tr(\Sigma^2) \right)^2 + (\tau^2 - \tau) \; \Tr(\Sigma^4)
    \quad \text{and} \quad
    \kappa \leq \tau \, \|\Sigma\|^2.
\]
Applying Lemma \ref{lem:variance}, we get the assertion of the theorem: with probability at least $(1 - \delta)$ it holds that
\[
    n \left\|\widehat \Sigma - \Sigma \right\|_{\Fr}^2 - n \; \E \left\|\widehat \Sigma - \Sigma \right\|_{\Fr}^2 < 4 \|\Sigma\|^2 \max\left\{ \sqrt{2 \big( \tau \ttr(\Sigma^2)^2 + \tau^2 \ttr(\Sigma^4) \big) \log(2 / \delta)}, \; 4e \tau \log(2 / \delta) \right\}.
\]
The proof is finished.

\myendproof

\subsection{Proof of Theorem \ref{th:fr_norm_lower_bound}}
\label{sec:th_fr_norm_lower_bound_proof}

Since, for any $t > 0$ and any $\lambda > 0$, it holds that
\begin{align}
    \label{eq:markov_inequality}
    &\notag
    \p\left( n \|\widehat \Sigma - \Sigma \|_{\Fr}^2 - n \, \E \|\widehat \Sigma - \Sigma \|_{\Fr}^2 \leq -t \right)
    \\&
    \leq \exp\left\{-\frac{\lambda t}2 + \frac{\lambda n \, \E \|\widehat \Sigma - \Sigma \|_{\Fr}^2}2 \right\} \cdot \E \exp\left\{ - \frac{\lambda n \|\widehat \Sigma - \Sigma \|_{\Fr}^2}2 \right\}
\end{align}
due to the Markov inequality, we are interested in upper bounds on the exponential moment
\[
    \E \exp\left\{ - \frac{\lambda n \|\widehat \Sigma - \Sigma \|_{\Fr}^2}2 \right\}.
\]
We would like to note that, in contrast to Theorem \ref{th:fr_norm_upper_bound}, the exponential moment of interest exists for all $\lambda > 0$ in this case. This slightly simplifies the proof. In what follows, we will show that, under the conditions of Theorem \ref{th:fr_norm_lower_bound}, for any $t > 0$, it holds that
\[
    \p\left( n \|\widehat \Sigma - \Sigma \|_{\Fr}^2 - n \, \E \|\widehat \Sigma - \Sigma \|_{\Fr}^2 \leq -t \right)
    \leq 5 \exp\left\{ -\left( \frac{t^2}{8 \tts^2} \wedge \frac{t}{8 \alpha \|\Sigma\|^2} \right) \right\}
\]
with
\[
    \tts^2 =  4 \ttv^2 + \frac{10 \alpha^2 \Tr(\Sigma)^2}n \left( \Tr(\Sigma)^2 + 2 \|\Sigma\|_{\Fr}^2 \left[ \frac{\ttr(\Sigma)^2}{\alpha} + \ttr(\Sigma^2) \right] \right),
\]
where $\ttv$ is defined in \eqref{eq:v_kappa}. This will yield the desired high probability bound \eqref{eq:fr_norm_lower_bound}.
Similarly to the proof of Theorem \ref{th:fr_norm_upper_bound}, we split our derivations in several steps for the ease of presentation.

\medskip

\noindent\textbf{Step 1: linearization.}\quad
As before, let us denote
\[
    H = \sqrt n \; \Sigma^{-1/2} \left( \widehat \Sigma - \Sigma \right) \Sigma^{-1/2}
    = \frac1{\sqrt n} \sum\limits_{i = 1}^n \left( \bxi_i \bxi_i^\top - I \right),
\]
where $\bxi_i = \Sigma^{-1/2} \bX_i$ for all $i \in \{1, \dots, n\}$ are the standardized vectors.
Let $\Gamma \in \R^{d\times d}$ be a random matrix with i.i.d. standard Gaussian entries. Applying the same linearization trick as in the proof of Theorem \ref{th:fr_norm_upper_bound} (Step 1), we obtain that
\[
    \E \exp\left\{ - \frac{\lambda n \|\widehat \Sigma - \Sigma \|_{\Fr}^2}2 \right\}
    = \E \exp\left\{ - \frac{\lambda \| \Sigma^{1/2} H \Sigma^{1/2} \|_{\Fr}^2}2 \right\}
    = \E \exp\left\{\i \sqrt{\lambda} \; \Tr(\Gamma^\top \Sigma^{1/2} H \Sigma^{1/2}) \right\}.
\]
Hence, we have to focus on the properties of the characteristic function of $H$. For any $U \in \R^{d \times d}$, let
\[
    \psi(U) = \log \E e^{\i \bxi^\top U \bxi}
\]
be the logarithm of the Fourier transform of $\bxi^\top U \bxi$. Applying the Fubini theorem, we get that
\begin{align*}
    \E \exp\left\{\i \sqrt{\lambda} \; \Tr(\Gamma^\top \Sigma^{1/2} H \Sigma^{1/2}) \right\}
    &
    = \E \exp\left\{\i \sqrt{\frac{\lambda}n} \; \sum\limits_{j = 1}^n \Tr\left(\Gamma^\top \Sigma^{1/2} \left( \bxi_j \bxi_j^\top - I \right) \Sigma^{1/2} \right) \right\}
    \\&
    = \E_\Gamma \exp\left\{ n \psi\left( \sqrt{\frac{\lambda}n} \; \Sigma^{1/2} \Gamma \Sigma^{1/2} \right) - \i \sqrt{n \lambda} \; \Tr(\Sigma^{1/2} \Gamma \Sigma^{1/2}) \right\}.
\end{align*}
Thus, we proved that
\[
    \E \exp\left\{ - \frac{\lambda n \|\widehat \Sigma - \Sigma \|_{\Fr}^2}2 \right\}
    = \E_\Gamma \exp\left\{ n \psi\left(\sqrt{\frac{\lambda}n} \; \Sigma^{1/2} \Gamma \Sigma^{1/2} \right) - \i \sqrt{n \lambda} \; \Tr(\Sigma^{1/2} \Gamma \Sigma^{1/2}) \right\}.
\]
On the other hand, the absolute value of $\E \exp\left\{ - \lambda n \|\widehat \Sigma - \Sigma \|_{\Fr}^2 / 2 \right\}$ is equal to the exponential moment itself. Then, due to the Jensen inequality, we have
\begin{align*}
    \E \exp\left\{ - \frac{\lambda n \|\widehat \Sigma - \Sigma \|_{\Fr}^2}2 \right\}
    &
    = \left| \E \exp\left\{ - \frac{\lambda n \|\widehat \Sigma - \Sigma \|_{\Fr}^2}2 \right\} \right|
    \\&
    = \left| \E_\Gamma \exp\left\{ n \psi\left( \sqrt{\frac{\lambda}n} \; \Sigma^{1/2} \Gamma \Sigma^{1/2} \right) - \i \sqrt{n \lambda} \; \Tr(\Sigma^{1/2} \Gamma \Sigma^{1/2}) \right\} \right|
    \\&
    \leq \E_\Gamma \left| \exp\left\{ n \psi\left( \sqrt{\frac{\lambda}n} \; \Sigma^{1/2} \Gamma \Sigma^{1/2} \right) - \i \sqrt{n \lambda} \; \Tr(\Sigma^{1/2} \Gamma \Sigma^{1/2}) \right\} \right|.
\end{align*}
Let us introduce a positive number $\rho > 0$ such that
\begin{equation}
    \label{eq:rho}
    \rho^2 = \Tr(\Sigma)^2 + 2 \|\Sigma\|_{\Fr}^2 \left( \frac{\ttr(\Sigma)^2}{\alpha} + \ttr(\Sigma^2) \right)
\end{equation}
and define an event $\cE_\rho$ as follows:
\[
    \cE_{\rho} = \left\{ \left| \Tr(\Gamma \Sigma) \right| \leq \rho
    \text{ and }
    \left\|\Sigma^{1/2} \Gamma \Sigma^{1/2} \right\|_{\Fr} \leq \rho \right\}.
\]
Then it holds that
\begin{align*}
    &
    \E_\Gamma \left| \exp\left\{ n \psi\left( \sqrt{\frac{\lambda}n} \; \Sigma^{1/2} \Gamma \Sigma^{1/2} \right) - \i \sqrt{n \lambda} \; \Tr(\Sigma^{1/2} \Gamma \Sigma^{1/2}) \right\} \right|
    \\&
    = \E_\Gamma \left[ \left| \exp\left\{ n \psi\left( \sqrt{\frac{\lambda}n} \; \Sigma^{1/2} \Gamma \Sigma^{1/2} \right) - \i \sqrt{n \lambda} \; \Tr(\Sigma^{1/2} \Gamma \Sigma^{1/2}) \right\} \right| \1\left( \cE_\rho \right) \right]
    \\&\quad
    + \E_\Gamma \left[ \left| \exp\left\{ n \psi\left( \sqrt{\frac{\lambda}n} \; \Sigma^{1/2} \Gamma \Sigma^{1/2} \right) - \i \sqrt{n \lambda} \; \Tr(\Sigma^{1/2} \Gamma \Sigma^{1/2}) \right\} \right| \1\left( \cE_\rho^c \right) \right],
\end{align*}
where $\cE_\rho^c$ stands for the complement of $\cE_\rho$. On $\cE_\rho^c$, we apply the bound
\begin{align*}
    &
    \E_\Gamma \left[ \left| \exp\left\{ n \psi\left( \sqrt{\frac{\lambda}n} \; \Sigma^{1/2} \Gamma \Sigma^{1/2} \right) - \i \sqrt{n \lambda} \; \Tr(\Sigma^{1/2} \Gamma \Sigma^{1/2}) \right\} \right| \1\left( \cE_\rho^c \right) \right]
    \\&
    = \E_\Gamma \left[ \left| \E_{\bxi} \exp\left\{\i \sqrt{\lambda} \; \Tr(\Gamma^\top \Sigma^{1/2} H \Sigma^{1/2}) \right\} \right| \1\left( \cE_\rho^c \right) \right]
    \\&
    \leq \E_\Gamma \left[ 1 \cdot \1\left( \cE_\rho^c \right) \right]
    = \p\left( \cE_\rho^c \right).
\end{align*}
On the other hand, on $\cE_\rho$ we use the inequality
\begin{align*}
    &
    \E_\Gamma \left[ \left| \exp\left\{ n \psi\left( \sqrt{\frac{\lambda}n} \; \Sigma^{1/2} \Gamma \Sigma^{1/2} \right) - \i \sqrt{n \lambda} \; \Tr(\Sigma^{1/2} \Gamma \Sigma^{1/2}) \right\} \right| \1\left( \cE_\rho \right) \right]
    \\&
    \leq \E_\Gamma \left[ \exp\left\{ n \; \Re\left[ \psi\left( \sqrt{\frac{\lambda}n} \; \Sigma^{1/2} \Gamma \Sigma^{1/2} \right) \right] \right\} \1\left( \cE_\rho \right) \right].
\end{align*}
Hence, we showed that
\begin{equation}
    \label{eq:exp_moment_split}
    \E \exp\left\{ - \frac{\lambda n \|\widehat \Sigma - \Sigma \|_{\Fr}^2}2 \right\}
    \leq \E_\Gamma \left[ \exp\left\{ n \; \Re\left[ \psi\left( \sqrt{\frac{\lambda}n} \; \Sigma^{1/2} \Gamma \Sigma^{1/2} \right) \right] \right\} \1\left( \cE_\rho \right) \right]
    + \p\left( \cE_\rho^c \right).
\end{equation}
In what follows, we bound the terms in the right-hand side of \eqref{eq:exp_moment_split} one by one, starting with the second summand.

\medskip

\noindent\textbf{Step 2: a bound on the remainder term.}\quad
According to the definition of $\cE_\rho$, it is enough to bound the probabilities $\p \left(|\Tr(\Gamma \Sigma)| > \rho \right)$ and $\p\left( \|\Sigma^{1/2} \Gamma \Sigma^{1/2} \|_{\Fr} > \rho \right)$. Note that $\Tr(\Gamma \Sigma)$ is a centered Gaussian random variable with variance $\|\Sigma\|_{\Fr}^2$. Then, due to the Hoeffding inequality, it holds that
\begin{equation}
    \label{eq:gamma_trace_bound}
    \p \left( |\Tr(\Gamma \Sigma)| > \rho \right) \leq 2 \exp \left\{-\frac{\rho^2}{2 \|\Sigma\|_{\Fr}^2} \right\}.
\end{equation}
The bound on $\p\left( \|\Sigma^{1/2} \Gamma \Sigma^{1/2}\|_{\Fr} > \rho \right)$ follows from the standard results on large deviations of Gaussian quadratic forms. Let us represent the squared Frobenius norm of $\Sigma^{1/2} \Gamma \Sigma^{1/2}$ in the following form, using \eqref{eq:kronecker_trace}:
\[
    \|\Sigma^{1/2} \Gamma \Sigma^{1/2}\|_{\Fr}^2
    = \Tr\left( \Sigma^{1/2} \Gamma^\top \Sigma \Gamma \Sigma^{1/2} \right)
    = \Tr\left( \Gamma^\top \Sigma \Gamma \Sigma \right)
    = \rmvec(\Gamma)^\top \left( \Sigma \otimes \Sigma \right) \rmvec(\Gamma).
\]
Then, according to \cite{laurent00} (see the proof of Lemma 1), for any $\mu \in (0, 1/2)$, it holds that
\begin{align*}
    \p\left( \|\Sigma^{1/2} \Gamma \Sigma^{1/2}\|_{\Fr} > \rho \right)
    &
    = \p\left( \rmvec(\Gamma)^\top \left( \Sigma \otimes \Sigma \right) \rmvec(\Gamma) > \rho^2 \right)
    \\&
    \leq \exp\left\{-\mu \rho^2 + \mu \Tr(\Sigma \otimes \Sigma) + \frac{\|\Sigma \otimes \Sigma\|_{\Fr}^2 \; \mu^2}{1 - 2 \|\Sigma \otimes \Sigma\| \, \mu} \right\}.
\end{align*}
Due to the properties of the Kronecker product (see \eqref{eq:kronecker_eigenvalues}), we have
\[
    \Tr(\Sigma \otimes \Sigma) = \left( \Tr(\Sigma) \right)^2,
    \quad
    \|\Sigma \otimes \Sigma\|_{\Fr}
    = \|\Sigma\|_{\Fr}^2
    \quad \text{and} \quad
    \|\Sigma \otimes \Sigma\|
    = \|\Sigma\|^2.
\]
Hence,
\begin{equation}
    \label{eq:gamma_fr_norm_bound}
    \p\left( \|\Sigma^{1/2} \Gamma \Sigma^{1/2}\|_{\Fr} > \rho \right)
    \leq \exp\left\{-\mu \left(\rho^2 - \left[ \Tr(\Sigma) \right]^2 \right) + \frac{\|\Sigma\|_{\Fr}^4 \; \mu^2}{1 - 2 \|\Sigma\|^2 \mu} \right\}
    \quad \text{for all $\mu \in (0, 1/2)$.}
\end{equation}
Let us take $\mu \in (0, 1/2)$ satisfying the condition
\[
    1 - 2 \|\Sigma\|^2 \mu = \frac{\|\Sigma\|_{\Fr}^4}{\|\Sigma\|_{\Fr}^4 + \|\Sigma\|^2 (\rho^2 - \left[ \Tr(\Sigma) \right]^2 )},
    \quad \text{that is,} \quad
    \mu = \frac{(\rho^2 - \left[ \Tr(\Sigma) \right]^2 ) / 2}{\|\Sigma\|_{\Fr}^4 + \|\Sigma\|^2 (\rho^2 - \left[ \Tr(\Sigma) \right]^2 )}.
\]
Then it is straightforward to check that
\begin{align*}
    -\mu \left(\rho^2 - \left[ \Tr(\Sigma) \right]^2 \right) + \frac{\|\Sigma\|_{\Fr}^4 \; \mu^2}{1 - 2 \|\Sigma\|^2 \mu}
    &
    = - \mu \left(\rho^2 - \left[ \Tr(\Sigma) \right]^2 \right) + \mu^2 \left( \|\Sigma\|_{\Fr}^4 + \|\Sigma\|^2 \rho^2 \right)
    \\&
    = - \frac{(\rho^2 - \left[ \Tr(\Sigma) \right]^2 )^2 / 4}{\|\Sigma\|_{\Fr}^4 + \|\Sigma\|^2 (\rho^2 - \left[ \Tr(\Sigma) \right]^2 )}.
\end{align*}
Substituting this equality into \eqref{eq:gamma_fr_norm_bound}, we finally obtain that
\[
    \p\left( \|\Sigma^{1/2} \Gamma \Sigma^{1/2}\|_{\Fr} > \rho \right)
    \leq \exp\left\{ - \frac{(\rho^2 - \left[ \Tr(\Sigma) \right]^2 )^2 / 4}{\|\Sigma\|_{\Fr}^4 + \|\Sigma\|^2 (\rho^2 - \left[ \Tr(\Sigma) \right]^2 )} \right\}.
\]
This inequality and \eqref{eq:gamma_trace_bound} immediately yield that 
\begin{align}
    \label{eq:e_complement}
    \p\left( \cE_\rho^c \right)
    &\notag
    \leq \p \left( |\Tr(\Gamma \Sigma)| > \rho \right) + \p\left( \|\Sigma^{1/2} \Gamma \Sigma^{1/2}\|_{\Fr} > \rho \right)
    \\&
    \leq 2 \exp \left\{-\frac{\rho^2}{2 \|\Sigma\|_{\Fr}^2} \right\} + \exp\left\{ - \frac{(\rho^2 - \left[ \Tr(\Sigma) \right]^2 )^2 / 4}{\|\Sigma\|_{\Fr}^4 + \|\Sigma\|^2 (\rho^2 - \left[ \Tr(\Sigma) \right]^2 )} \right\}.
\end{align}

\medskip

\noindent\textbf{Step 3: Taylor's expansion.}\quad
Next, we elaborate on the first term in the right-hand side of \eqref{eq:exp_moment_split}. For any $U \in \R^{d\times d}$, whenever $\psi(U)$ is well-defined, the gradient $\nabla \psi(U)$ is given by
\[
    \left\langle \nabla \psi(U), V \right\rangle = \frac{\i \E \left[ \bxi^\top V \bxi \; e^{\i \bxi^\top U \bxi} \right]}{\E e^{\i \bxi^\top U \bxi}}
    = \i \sfP_{\i U} \, \bxi^\top V \bxi.
\]
Here we extended the notation $\sfP_U f(\bxi, U)$ (see Section \ref{sec:main_results}) to the complex-valued matrices. Note that the extension is formal and $\sfP_{\i U}$ is not a probability measure anymore. Nevertheless, we still can compute the higher-order derivatives of $\psi(U)$ in similar way as in the proof Lemma \ref{lem:derivatives}:
\begin{align*}
    &
    \left\langle \nabla^2 \psi(O), V^{\otimes 2} \right\rangle
    = - \E \left(\bxi^\top V \bxi - \E \bxi^\top V \bxi \right)^2,
    \\&
    \left\langle \nabla^3 \psi(O), V^{\otimes 3} \right\rangle
    = -\i \; \E \left(\bxi^\top V \bxi - \E \bxi^\top V \bxi \right)^3,
    \\&
    \left\langle \nabla^4 \psi(U), V^{\otimes 4} \right\rangle
    = \sfP_{\i U} \left(\bxi^\top V \bxi - \sfP_{\i U} \bxi^\top V \bxi \right)^4 - 3 \left( \sfP_{\i U} \left(\bxi^\top V \bxi - \sfP_{\i U} \bxi^\top V \bxi \right)^2 \right)^2.
\end{align*}
As before, $O \in \R^{d \times d}$ stands the matrix with zero entries.
Applying Taylor's formula with the Lagrange remainder term, we obtain that
\begin{align*}
    \psi\left( \sqrt{\frac{\lambda}n} \; \Sigma^{1/2} \Gamma \Sigma^{1/2} \right)
    &
    = \sqrt{\frac{\lambda}n} \left\langle \nabla \psi(O), \Sigma^{1/2} \Gamma \Sigma^{1/2} \right\rangle
    + \frac{\lambda}{2n} \left\langle \nabla^2 \psi(O), \left(\Sigma^{1/2} \Gamma \Sigma^{1/2} \right)^{\otimes 2} \right\rangle
    \\&
    + \frac16 \left( \frac{\lambda}{n} \right)^{3/2} \left\langle \nabla^3 \psi(O), \left(\Sigma^{1/2} \Gamma \Sigma^{1/2} \right)^{\otimes 3} \right\rangle
    + \frac{\lambda^2}{24 n^2} \left\langle \nabla^4 \psi(\Theta), \left(\Sigma^{1/2} \Gamma \Sigma^{1/2} \right)^{\otimes 4} \right\rangle
\end{align*}
for some $\Theta$ such that
\[
    \|\Theta\|_{\Fr}
    \leq \sqrt{\frac\lambda n} \; \left\|\Sigma^{1/2} \Gamma \Sigma^{1/2} \right\|_{\Fr}
    \leq \rho \sqrt{\frac\lambda n}
    \quad \text{on $\cE_\rho$.}
\]
Taking into account that $\langle \nabla\psi(O), \Sigma^{1/2} \Gamma \Sigma^{1/2} \rangle$ and $\langle \nabla^3 \psi(O), (\Sigma^{1/2} \Gamma \Sigma^{1/2})^{\otimes 3} \rangle$ are imaginary, we immediately obtain
\[
    \Re\left[ \psi\left( \sqrt{\frac{\lambda}n} \; \Sigma^{1/2} \Gamma \Sigma^{1/2} \right) \right]
    = \frac{\lambda}{2n} \left\langle \nabla^2 \psi(O), \left(\Sigma^{1/2} \Gamma \Sigma^{1/2} \right)^{\otimes 2} \right\rangle + \frac{\lambda^2}{24 n^2} \left\langle \nabla^4 \psi(\Theta), \left(\Sigma^{1/2} \Gamma \Sigma^{1/2} \right)^{\otimes 4} \right\rangle.
\]
Thus, our next goal is to bound the fourth derivative of $\psi$. To do so, we use the following lemma.
\begin{Lem}
    \label{lem:kth_moment}
    Grant Assumption \ref{as:fourth_moment} and suppose that $|\E e^{\i \bxi^\top U \bxi}| \geq 1 / \beta$ for some $\beta > 1$. Then, for any $k \in [1, 4]$, it holds that
    \[
        \left| \sfP_{\i U} \left( \bxi^\top V \bxi - \sfP_{\i U} \bxi^\top V \bxi \right)^k \right| \leq 2^{k - 1} \beta \left( 1 + \beta^k \right) \alpha^{k/2} \|V\|_{\Fr}^k
        \quad \text{for all $V \in \R^{d \times d}$.}
    \]
\end{Lem}

The proof of Lemma \ref{lem:kth_moment} is deferred to Appendix \ref{sec:lem_kth_moment_proof}. Note that on the event $\cE_\rho$, we have
\begin{align*}
    \left| \E e^{\i \bxi^\top \Theta \bxi} - 1 \right|
    &
    = 2 \left| \E \sin\left(\bxi^\top \Theta \bxi / 2 \right) e^{\i \bxi^\top \Theta \bxi / 2} \right|
    \leq \E \left| \bxi^\top \Theta \bxi \right|
    \\&
    \leq |\Tr(\Theta)| + \sqrt{\alpha} \|\Theta\|_{\Fr}
    \leq \rho \left(1 + \sqrt{\alpha} \right) \sqrt{\frac\lambda n}.
\end{align*}
In what follows, we will choose $\lambda \leq 1 / (2 \alpha \|\Sigma\|^2)$. Invoking the definition of $\rho$ (see eq. \eqref{eq:rho}) and using the condition \eqref{eq:n_large}, we observe that
\begin{align*}
    \rho^2 \left(1 + \sqrt{\alpha} \right)^2 \frac\lambda n
    &
    \leq \frac{3 \left(1 + \sqrt{\alpha} \right)^2}{2 \alpha \|\Sigma\|^2 n} \left( \Tr(\Sigma)^2 + 2 \|\Sigma\|_{\Fr}^2 \left[ \frac{\ttr(\Sigma)^2}\alpha + \ttr(\Sigma^2) \right] \right)
    \\&
    = \frac{3 \left(1 + \sqrt{\alpha} \right)^2}{2 \alpha n} \left( \ttr(\Sigma)^2 + 2 \ttr(\Sigma^2) \left[ \frac{\ttr(\Sigma)^2}\alpha + \ttr(\Sigma^2) \right] \right)
    \\&
    \leq \frac{3 \left(1 + \sqrt{\alpha} \right)^2}{2 \alpha n} \left( 3 \ttr(\Sigma)^2 + \frac{2 \ttr(\Sigma^2) \ttr(\Sigma)^2}\alpha \right)
    \leq \frac1{64}.
\end{align*}
In the last line, we used the fact that $\ttr(\Sigma^2) \leq \ttr(\Sigma)$.
Hence, the conditions of Lemma \ref{lem:kth_moment} are fulfilled with
\[
    \beta
    = \left( 1 - \rho \left(1 + \sqrt{\alpha} \right) \sqrt{\frac\lambda n} \right)^{-1}
    \leq \left( 1 - 1/8 \right)^{-1}
    = \frac87
    < 2^{1/5},
\]
and consequently, on $\cE_{\rho}$, it holds that
\begin{align*}
    \left\langle \nabla^4 \psi(\Theta), \left(\Sigma^{1/2} \Gamma \Sigma^{1/2} \right)^{\otimes 4} \right\rangle
    &
    \leq \left( 8 \beta (1 + \beta^4) + 12 \beta^2 \left(1 + \beta^2 \right)^2 \right) \alpha^2 \|\Sigma^{1/2} \Gamma \Sigma^{1/2}\|_{\Fr}^4
    \\&
    < 111 \alpha^2 \rho^2 \|\Sigma^{1/2} \Gamma \Sigma^{1/2}\|_{\Fr}^2.
\end{align*}
This implies that
\begin{align*}
    \Re\left[ \psi\left( \sqrt{\frac{\lambda}n} \; \Sigma^{1/2} \Gamma \Sigma^{1/2} \right) \right]
    &
    = \frac{\lambda}{2n} \left\langle \nabla^2 \psi(O), \left(\Sigma^{1/2} \Gamma \Sigma^{1/2} \right)^{\otimes 2} \right\rangle + \frac{\lambda^2}{24 n^2} \left\langle \nabla^4 \psi(\Theta), \left(\Sigma^{1/2} \Gamma \Sigma^{1/2} \right)^{\otimes 4} \right\rangle
    \\&
    \leq -\frac{\lambda}{2n} \E_{\bxi} \left( \bxi^\top \Sigma^{1/2} \Gamma \Sigma^{1/2} \bxi - \Tr(\Gamma \Sigma) \right)^2 + \frac{111 \alpha^2 \lambda^2}{24 n^2} \left\|\Sigma^{1/2} \Gamma \Sigma^{1/2} \right\|_{\Fr}^4
    \\&
    \leq -\frac{\lambda}{2n} \E_{\bX} \left( \bX^\top \Gamma \bX - \Tr(\Gamma \Sigma) \right)^2 + \frac{5 \alpha^2 \lambda^2 \rho^2}{n^2} \left\|\Sigma^{1/2} \Gamma \Sigma^{1/2} \right\|_{\Fr}^2
\end{align*}
on the same event. Hence, it holds that
\begin{align*}
    &
    \exp\left\{ n \; \Re\left[ \psi\left( \sqrt{\frac{\lambda}n} \; \Sigma^{1/2} \Gamma \Sigma^{1/2} \right) \right] \right\} \1\left( \|\Sigma^{1/2} \Gamma \Sigma^{1/2}\|_{\Fr} \leq \rho \right)
    \\&
    \leq \exp\left\{ -\frac{\lambda}2 \; \E_{\bX} \left( \bX^\top \Gamma \bX - \Tr(\Gamma \Sigma) \right)^2 + \frac{5 \alpha^2 \lambda^2 \rho^2}{n} \left\|\Sigma^{1/2} \Gamma \Sigma^{1/2} \right\|_{\Fr}^2 \right\}.
\end{align*}
The logarithm of the right-hand side is quadratic in $\Gamma$, and thus, we can easily bound its exponential moment. The precise statement is given in the following lemma.

\begin{Lem}
    \label{lem:exp_moment_bound}
    Assume that
    \[
        \frac{10 \alpha^2 \lambda^2 \rho^2 \|\Sigma\|^2}{n} \leq \frac12.
    \]
    Then it holds that
    \begin{align*}
        &
        \E_\Gamma \exp\left\{ -\frac{\lambda}2 \; \E_{\bX} \left( \bX^\top \Gamma \bX - \Tr(\Gamma \Sigma) \right)^2 + \frac{5 \alpha^2 \lambda^2 \rho^2}{n} \left\|\Sigma^{1/2} \Gamma \Sigma^{1/2} \right\|_{\Fr}^2 \right\}
        \\&
        \leq \exp\left\{ -\frac{\lambda \ttm}2 + 2 \lambda^2 \ttv^2 + \frac{5 \alpha^2 \lambda^2 \rho^2}{n} \left( \Tr(\Sigma) \right)^2 + \frac{200 \alpha^4 \lambda^4 \rho^4 \|\Sigma\|_{\Fr}^4}{n^2} \right\}.
    \end{align*}
\end{Lem}
We present the proof of Lemma \ref{lem:exp_moment_bound} in Appendix \ref{sec:lem_exp_moment_bound_proof}. Hence, on the third step, we proved that
\begin{equation}
    \label{eq:real_part_bound}
    \Re\left[ \psi\left( \sqrt{\frac{\lambda}n} \; \Sigma^{1/2} \Gamma \Sigma^{1/2} \right) \right]
    \leq \exp\left\{ -\frac{\lambda \ttm}2 + 2 \lambda^2 \ttv^2 + \frac{5 \alpha^2 \lambda^2 \rho^2}{n} \left( \Tr(\Sigma) \right)^2 + \frac{200 \alpha^4 \lambda^4 \rho^4 \|\Sigma\|_{\Fr}^4}{n^2} \right\}.
\end{equation}

\medskip

\noindent\textbf{Step 4: final bound.}
Summing up the inequalities \eqref{eq:markov_inequality}, \eqref{eq:exp_moment_split}, \eqref{eq:e_complement} and \eqref{eq:real_part_bound}, we obtain that the following bound holds for any positive $t$ and $\lambda$:
\begin{align}
    \label{eq:left_deviation_bound}
    &\notag
    \p\left( n \|\widehat \Sigma - \Sigma \|_{\Fr}^2 - n \, \E \|\widehat \Sigma - \Sigma \|_{\Fr}^2 \leq -t \right)
    \\&
    \leq \exp\left\{ -\frac{\lambda t}2 + 2 \lambda^2 \ttv^2 + \frac{5 \alpha^2 \lambda^2 \rho^2}{n} \left( \Tr(\Sigma) \right)^2 + \frac{200 \alpha^4 \lambda^4 \rho^4 \|\Sigma\|_{\Fr}^4}{n^2} \right\}
    \\&\notag\quad
    + e^{-\lambda t / 2 + \lambda \ttm / 2} \left( 2 \exp \left\{-\frac{\rho^2}{2 \|\Sigma\|_{\Fr}^2} \right\} + \exp\left\{ - \frac{(\rho^2 - \left[ \Tr(\Sigma) \right]^2 )^2 / 4}{\|\Sigma\|_{\Fr}^4 + \|\Sigma\|^2 (\rho^2 - \left[ \Tr(\Sigma) \right]^2 )} \right\} \right).    
\end{align}
In what follows, we will take $\lambda$ from $\big[0, 1 / (2 \alpha \|\Sigma\|^2) \big]$. In view of Lemma \ref{lem:mean}, this implies that
\[
    \frac{\lambda \ttm}2
    \leq \frac14 \left( \frac{\ttr(\Sigma)^2}{\alpha} + \ttr(\Sigma^2) \right).
\]
The next lemma helps us to simplify the inequality \eqref{eq:left_deviation_bound} substantially.

\begin{Lem}
    \label{lem:rho_inequalities}
    Let $\rho > 0$ be as defined in \eqref{eq:rho}. Then it holds that
    \[
        \frac{2 \rho^2}{\|\Sigma\|_{\Fr}^2} \wedge \frac{(\rho^2 - \left[ \Tr(\Sigma) \right]^2 )^2}{\|\Sigma\|_{\Fr}^4 + \|\Sigma\|^2 (\rho^2 - \left[ \Tr(\Sigma) \right]^2 )} \geq \frac{\ttr(\Sigma)^2}{\alpha} + \ttr(\Sigma^2).
    \]
\end{Lem}
We defer the proof of Lemma \ref{lem:rho_inequalities} to Appendix \ref{sec:lem_rho_inequalities_proof}. With this lemma at hand, we can conclude that
\begin{align*}
    &
    \p\left( n \|\widehat \Sigma - \Sigma \|_{\Fr}^2 - n \, \E \|\widehat \Sigma - \Sigma \|_{\Fr}^2 \leq -t \right)
    \\&
    \leq e^{-\lambda t / 2} \left( 3 + \exp\left\{2 \lambda^2 \ttv^2 + \frac{5 \alpha^2 \lambda^2 \rho^2}{n} \left( \Tr(\Sigma) \right)^2 + \frac{200 \alpha^4 \lambda^4 \rho^4 \|\Sigma\|_{\Fr}^4}{n^2} \right\} \right).
\end{align*}
This expression can be simplified even further, because
\begin{align*}
    \frac{10 \alpha^2 \lambda^2 \rho^2 \|\Sigma\|_{\Fr}^2}n
    &
    \leq \frac{10 \alpha^2 \|\Sigma\|_{\Fr}^2}{4 \alpha^2 \|\Sigma\|^4 n} \left[ \Tr(\Sigma)^2 + 2 \|\Sigma\|_{\Fr}^2 \left( \frac{\ttr(\Sigma)^2}{\alpha} + \ttr(\Sigma^2) \right) \right]
    \\&
    = \frac{5 \ttr(\Sigma^2)}{2n} \left( 3 \ttr(\Sigma)^2 + \frac{2 \ttr(\Sigma^2) \ttr(\Sigma)^2}{\alpha}\right)
    \leq \frac12 < \sqrt{\frac{\log 2}2}
\end{align*}
due to the condition of the theorem. Hence,
\begin{align*}
    \p\left( n \|\widehat \Sigma - \Sigma \|_{\Fr}^2 - n \, \E \|\widehat \Sigma - \Sigma \|_{\Fr}^2 \leq -t \right)
    &
    \leq e^{-\lambda t / 2} \left( 3 + 2 \exp\left\{2 \lambda^2 \ttv^2 + \frac{5 \alpha^2 \lambda^2 \rho^2}{n} \left( \Tr(\Sigma) \right)^2 \right\} \right)
    \\&
    \leq 5 \exp\left\{-\frac{\lambda t}2 + 2 \lambda^2 \ttv^2 + \frac{5 \alpha^2 \lambda^2 \rho^2}{n} \left( \Tr(\Sigma) \right)^2 \right\}.
\end{align*}
It remains to minimize the expression
\[
    -\frac{\lambda t}2 + 2 \lambda^2 \ttv^2 + \frac{5 \alpha^2 \lambda^2 \rho^2}{n} \left( \Tr(\Sigma) \right)^2
    \quad \text{over $\lambda \in \left[0, \frac1{2 \alpha \|\Sigma\|^2} \right]$.}
\]
For this purpose, we take
\[
    \lambda = \frac{t}{8 \ttv^2 + 20 \alpha^2 \rho^2 \Tr(\Sigma)^2 / n} \wedge \frac1{2 \alpha \|\Sigma\|^2}.
\]
It is straightforward to check that, with such value of $\lambda$, we have
\[
    -\frac{\lambda t}2 + 2 \lambda^2 \ttv^2 + \frac{5 \alpha^2 \lambda^2 \rho^2}{n} \left( \Tr(\Sigma) \right)^2
    = -\min\left\{ \frac{t^2}{32 \ttv^2 + 80 \alpha^2 \rho^2 \Tr(\Sigma)^2 / n}, \frac{t}{8 \alpha \|\Sigma\|^2}\right\}.
\]
Thus, we proved that
\[
    \p\left( n \|\widehat \Sigma - \Sigma \|_{\Fr}^2 - n \, \E \|\widehat \Sigma - \Sigma \|_{\Fr}^2 \leq -t \right)
    \leq 5\exp\left\{ -\left( \frac{t^2}{32 \ttv^2 + 80 \alpha^2 \rho^2 \Tr(\Sigma)^2 / n} \wedge \frac{t}{8 \alpha \|\Sigma\|^2}\right) \right\},
\]
as we announced in the beginning of the proof.
In other words, for any $\delta > 0$, with probability at least $(1 - \delta)$, it holds that
\[
    n \, \E \|\widehat \Sigma - \Sigma \|_{\Fr}^2 - n \|\widehat \Sigma - \Sigma \|_{\Fr}^2
    < \left( 8 \alpha \|\Sigma\|^2 \log(5/\delta) \right) \vee \sqrt{\left(32 \ttv^2 + \frac{80 \alpha^2 \rho^2 \Tr(\Sigma)^2}n \right) \log(5/\delta)}.
\]
Taking into account that
\[
    \ttv^2 \leq \alpha \left( \Tr(\Sigma^2) \right)^2 + (\alpha^2 - \alpha) \; \Tr(\Sigma^4)
\]
due to Lemma \ref{lem:variance} and substituting $\rho$ with
\begin{align*}
    \rho^2
    &
    = \Tr(\Sigma)^2 + 2 \|\Sigma\|_{\Fr}^2 \left( \frac{\ttr(\Sigma)^2}{\alpha} + \ttr(\Sigma^2) \right)
    \\&
    = \|\Sigma\|^2 \left( \ttr(\Sigma)^2 + 2 \ttr(\Sigma^2) \left[ \frac{\ttr(\Sigma)^2}{\alpha} + \ttr(\Sigma^2) \right] \right)
    \\&
    \leq \|\Sigma\|^2 \left( 3 \ttr(\Sigma)^2 +  \frac{2 \ttr(\Sigma^2) \ttr(\Sigma)^2}{\alpha} \right),
\end{align*}
we get the desired bound:
\[
    \E \|\widehat \Sigma - \Sigma \|_{\Fr}^2 - \|\widehat \Sigma - \Sigma \|_{\Fr}^2
    < \frac{4 \|\Sigma\|^2}n \left( 2 \log(5/\delta)
    \vee \sqrt{2 \mathfrak R(\Sigma) \log(5/\delta)} \right),
\]
where
\[
    \mathfrak R(\Sigma) = \alpha \ttr(\Sigma^2)^2 + \alpha^2 \ttr(\Sigma^4) + \frac{5 \alpha^2 \ttr(\Sigma)^2}{2 n} \left( 3 \ttr(\Sigma)^2 +  \frac{2 \ttr(\Sigma^2) \ttr(\Sigma)^2}{\alpha} \right).
\]
\myendproof

\subsection{Proof of Corollary \ref{co:concentration}}
\label{sec:co_concentration_proof}

Let us remind a reader that, if Assumption \ref{as:derivatives} is satisfied, then Assumption \ref{as:fourth_moment} holds as well. In addition, the relations \eqref{eq:asymptotic_conditions}, together with the fact that $\tau \geq 2$, yield that the conditions of Theorems \ref{th:fr_norm_upper_bound} and \ref{th:fr_norm_lower_bound} (including the inequalities \eqref{eq:n_admissible} and \eqref{eq:n_large}) are fulfilled. Applying Theorem \ref{th:fr_norm_upper_bound} with the confidence level $\delta_1 = 3 \delta / 8$ and Theorem \ref{th:fr_norm_lower_bound} with $\delta_2 = 5 \delta / 8$ and using the union bound, we obtain that
\begin{equation}
    \label{eq:concentration}
    \left| \left\|\widehat \Sigma - \Sigma \right\|_{\Fr}^2 - \E \left\|\widehat \Sigma - \Sigma \right\|_{\Fr}^2 \right|
    \lesssim \frac{\tau \|\Sigma\|^2}n \max\left\{ \ttr(\Sigma^2) \sqrt{\log(8/\delta)}, \log(8/\delta) \right\} 
\end{equation}
with probability at least $(1 - \delta)$. Here we used the fact that $\ttr(\Sigma^4) \leq \ttr(\Sigma^2)$. From now on, we restrict our attention on the event where \eqref{eq:concentration} takes place. According to Lemma \ref{lem:mean}, it holds that
\[
    n \, \E \|\widehat\Sigma - \Sigma\|_{\Fr}^2
    \geq \left( \Tr(\Sigma) \right)^2 - \Tr(\Sigma^2)
    = \|\Sigma\|^2 \left(\ttr(\Sigma)^2 - \ttr(\Sigma^2) \right).
\]
Hence, we have
\begin{align*}
    \left| \frac{\|\widehat \Sigma - \Sigma \|_{\Fr}^2}{\E \|\widehat \Sigma - \Sigma \|_{\Fr}^2} - 1 \right|
    &
    \lesssim \tau \max\left\{ \frac{\ttr(\Sigma^2) \sqrt{\log(8/\delta)}}{\ttr(\Sigma)^2 - \ttr(\Sigma^2)}, \frac{\log(8/\delta)}{\ttr(\Sigma)^2 - \ttr(\Sigma^2)} \right\}.
\end{align*}
Finally, the inequalities $1 \leq \ttr(\Sigma^2) \leq \ttr(\Sigma)$ yield that
\[
    \left| \frac{\|\widehat \Sigma - \Sigma \|_{\Fr}^2}{\E \|\widehat \Sigma - \Sigma \|_{\Fr}^2} - 1 \right|
    \lesssim
    \tau \max\left\{ \frac{ \sqrt{\log(8/\delta)}}{\ttr(\Sigma) - 1}, \frac{\log(8/\delta)}{\ttr(\Sigma) \left( \ttr(\Sigma) - 1 \right)} \right\}.
\]
The proof is finished.

\myendproof

\subsection{Proof of Theorem \ref{th:fr_norm_fluctuation}}
\label{sec:th_fr_norm_fluctuation_proof}
Let us introduce
\begin{equation}
    \label{eq:upsilon}
    \upsilon = \frac{\|\widehat \Sigma - \Sigma\|_{\Fr}^2 - \E \|\widehat \Sigma - \Sigma\|_{\Fr}^2}{\E \|\widehat \Sigma - \Sigma\|_{\Fr}^2}
\end{equation}
and rewrite the difference $\|\widehat \Sigma - \Sigma\|_{\Fr} - \E \|\widehat \Sigma - \Sigma\|_{\Fr}$ in the following form:
\[
    \left| \|\widehat \Sigma - \Sigma\|_{\Fr} - \E \|\widehat \Sigma - \Sigma\|_{\Fr} \right|
    = \left| \frac{\|\widehat \Sigma - \Sigma\|_{\Fr}^2 - \left(\E \|\widehat \Sigma - \Sigma\|_{\Fr}\right)^2}{\|\widehat \Sigma - \Sigma\|_{\Fr} + \E \|\widehat \Sigma - \Sigma\|_{\Fr}} \right|.
\]
Applying the triangle inequality, we obtain that
\begin{equation}
    \label{eq:fr_norm_difference_abs_value_bound}
    \left| \|\widehat \Sigma - \Sigma\|_{\Fr} - \E \|\widehat \Sigma - \Sigma\|_{\Fr} \right|
    \leq \left| \frac{\|\widehat \Sigma - \Sigma\|_{\Fr}^2 - \E \|\widehat \Sigma - \Sigma\|_{\Fr}^2}{\|\widehat \Sigma - \Sigma\|_{\Fr} + \E \|\widehat \Sigma - \Sigma\|_{\Fr}} \right|
    + \frac{\Var\left(\|\widehat \Sigma - \Sigma\|_{\Fr}\right)}{\|\widehat \Sigma - \Sigma\|_{\Fr} + \E \|\widehat \Sigma - \Sigma\|_{\Fr}}.
\end{equation}
We are going to show that the expression in the right-hand side does not exceed
\[
    \sqrt{\E \left\|\widehat \Sigma - \Sigma \right\|_{\Fr}^2} \cdot \frac{|\upsilon| + \E \upsilon^2}{1 - \E \upsilon^2/2 + \sqrt{1 + \upsilon}}.
\]
For this purpose, we elaborate on $\|\widehat \Sigma - \Sigma\|_{\Fr}$, $\E \|\widehat \Sigma - \Sigma\|_{\Fr}$, and $\Var\big(\|\widehat \Sigma - \Sigma\|_{\Fr} \big)$.
For convenience, we split the rest of the proof into several steps.

\medskip

\noindent\textbf{Step 1: $\|\widehat \Sigma - \Sigma\|_{\Fr}^2 - \E \|\widehat \Sigma - \Sigma\|_{\Fr}^2$ and $\|\widehat \Sigma - \Sigma\|_{\Fr}^2$.}\quad
The first step is obvious. According to the definition of $\upsilon$ (see \eqref{eq:upsilon}), it holds that
\begin{equation}
    \label{eq:squared_fr_norm_upsilon}
    \left\|\widehat \Sigma - \Sigma\right\|_{\Fr}^2 - \E \left\|\widehat \Sigma - \Sigma\right\|_{\Fr}^2 = \upsilon \; \E \left\|\widehat \Sigma - \Sigma\right\|_{\Fr}^2
\end{equation}
and
\begin{equation}
    \label{eq:fr_norm_upsilon}
    \left\|\widehat \Sigma - \Sigma\right\|_{\Fr} = 
    \sqrt{\left\|\widehat \Sigma - \Sigma\right\|_{\Fr}^2}
    = \sqrt{(1 + \upsilon) \, \E \left\|\widehat \Sigma - \Sigma\right\|_{\Fr}^2}.
\end{equation}

\medskip

\noindent\textbf{Step 2: $\E \|\widehat \Sigma - \Sigma\|_{\Fr}$.}\quad Our next goal is to provide a lower bound on $\E \|\widehat \Sigma - \Sigma\|_{\Fr}$. Since $\sqrt{1 + x} \geq 1 + x/2 - x^2/2$ for any $x \geq -1$, the equality \eqref{eq:fr_norm_upsilon} yields that
\begin{align}
    \label{eq:fr_norm_expectation_upsilon}
    \E \left\|\widehat \Sigma - \Sigma\right\|_{\Fr}
    &\notag
    = \sqrt{\E \left\|\widehat \Sigma - \Sigma\right\|_{\Fr}^2} \cdot \E \sqrt{1 + \upsilon}
    \\&
    \geq \sqrt{\E \left\|\widehat \Sigma - \Sigma\right\|_{\Fr}^2} \left(1 + \frac{\E \upsilon}2 - \frac{\E \upsilon^2}2 \right)
    \\&\notag
    = \sqrt{\E \left\|\widehat \Sigma - \Sigma\right\|_{\Fr}^2} \left(1 - \frac{\E \upsilon^2}2 \right).
\end{align}
Here we used the fact that $\E \upsilon = 0$, which follows from the definition of $\upsilon$ (see \eqref{eq:upsilon}).

\medskip

\noindent\textbf{Step 3: $\Var\big(\|\widehat \Sigma - \Sigma\|_{\Fr} \big)$.}\quad Finally, let us elaborate on the variance of $\|\widehat \Sigma - \Sigma\|_{\Fr}$. The inequality \eqref{eq:fr_norm_expectation_upsilon} implies that
\begin{align}
    \label{eq:fr_norm_variance_upsilon}
    \Var\left( \|\widehat \Sigma - \Sigma\|_{\Fr} \right)
    &\notag
    = \E \left\|\widehat \Sigma - \Sigma\right\|_{\Fr}^2 - \left( \E \left\|\widehat \Sigma - \Sigma\right\|_{\Fr} \right)^2
    \\&
    \leq \E \left\|\widehat \Sigma - \Sigma\right\|_{\Fr}^2 - \E \left\|\widehat \Sigma - \Sigma\right\|_{\Fr}^2 \left( 1 - \frac{\E \upsilon^2}2 \right)^2
    \\&\notag
    \leq \E \upsilon^2 \; \E \left\|\widehat \Sigma - \Sigma\right\|_{\Fr}^2. 
\end{align}

\medskip

\noindent\textbf{Step 4: final bound.}
\quad
Summing up \eqref{eq:fr_norm_difference_abs_value_bound}, \eqref{eq:squared_fr_norm_upsilon}, \eqref{eq:fr_norm_upsilon}, \eqref{eq:fr_norm_expectation_upsilon}, and \eqref{eq:fr_norm_variance_upsilon}, we obtain that
\begin{equation}
    \label{eq:fr_norm_difference_abs_value_upsilon}
    \left| \|\widehat \Sigma - \Sigma\|_{\Fr} - \E \|\widehat \Sigma - \Sigma\|_{\Fr} \right|
    \leq \sqrt{\E \left\|\widehat \Sigma - \Sigma \right\|_{\Fr}^2} \cdot \frac{|\upsilon| + \E \upsilon^2}{1 - \E \upsilon^2/2 + \sqrt{1 + \upsilon}}.
\end{equation}
Since Assumption \ref{as:derivatives} implies that Assumption \ref{as:fourth_moment} holds with $\alpha = \tau$, according to Lemma \ref{lem:mean}, we have
\begin{equation}
    \label{eq:squared_fr_norm_expectation_two-side_bound}
    \frac{\|\Sigma\|^2}n \left( \ttr(\Sigma)^2 - \ttr(\Sigma^2) \right)
    \leq \E \left\|\widehat \Sigma - \Sigma \right\|_{\Fr}^2 \leq \frac{\|\Sigma\|^2}n \left( \ttr(\Sigma)^2 + (\tau  - 1) \; \ttr(\Sigma^2) \right).
\end{equation}
Moreover, under the conditions of the theorem, we can apply Corollary \ref{co:concentration}, which yields that
\begin{equation}
    \label{eq:upsilon_bound}
    |\upsilon|
    \lesssim \tau \max\left\{ \frac{\sqrt{\log(8 / \delta)}}{\ttr(\Sigma) - 1}, \frac{\log(8 / \delta)}{\ttr(\Sigma)\big( \ttr(\Sigma) - 1 \big)} \right\}
    \lesssim \frac{\tau \sqrt{\log(8 / \delta)}}{\ttr(\Sigma)}
\end{equation}
with probability at least $1 - \delta$. Finally, the following lemma helps us to bound $\E \upsilon^2$.
\begin{Lem}
    \label{lem:squared_fr_norm_variance_bound}
    Under Assumption \ref{as:fourth_moment}, it holds that
    \[
        \Var\left( \|\widehat\Sigma - \Sigma\|_{\Fr}^2 \right)
        \leq \frac{8 \Tr(\Sigma^2)^2 + 8\alpha \Tr(\Sigma^4)}{n^2}
        + \frac{8 \alpha^2 \Tr(\Sigma^2)^2 + 48 \alpha \Tr(\Sigma)^2 \Tr(\Sigma^2)}{n^3}.
    \]
\end{Lem}
The proof of Lemma \ref{lem:squared_fr_norm_variance_bound} is moved to Appendix \ref{sec:lem_squared_fr_norm_variance_bound_proof}. It immediately implies that
\[
    \Var\left( \|\widehat\Sigma - \Sigma\|_{\Fr}^2 \right)
    \lesssim \frac{\|\Sigma\|^4 \ttr(\Sigma^2) \big( \ttr(\Sigma) \vee \tau \big)}{n^2}
\]
and then
\begin{equation}
    \label{eq:upsilon_variance_bound}
    \E \upsilon^2 = \frac{\Var \big( \|\widehat\Sigma - \Sigma\|_{\Fr}^2 \big)}{\big( \E \|\widehat\Sigma - \Sigma\|_{\Fr}^2 \big)^2}
    \lesssim \frac{\|\Sigma\|^4 \ttr(\Sigma^2) \big( \ttr(\Sigma) \vee \tau \big) / n^2}{\|\Sigma\|^4 \ttr(\Sigma)^4 / n^2}
    \lesssim \frac1{\ttr(\Sigma)^2} \vee \frac{\tau}{\ttr(\Sigma)^3}.
\end{equation}
The inequalities  \eqref{eq:fr_norm_difference_abs_value_upsilon}, \eqref{eq:squared_fr_norm_expectation_two-side_bound}, \eqref{eq:upsilon_bound}, and \eqref{eq:upsilon_variance_bound} yield the assertion of the theorem:
\begin{align*}
    \left| \|\widehat \Sigma - \Sigma\|_{\Fr} - \E \|\widehat \Sigma - \Sigma\|_{\Fr} \right|
    &
    \lesssim \sqrt{\E \left\|\widehat \Sigma - \Sigma \right\|_{\Fr}^2} \cdot \frac{\tau \sqrt{\log(8 / \delta)}}{\ttr(\Sigma)}
    \\&
    \lesssim \sqrt{\|\Sigma\|^2 \ttr(\Sigma) \big( \ttr(\Sigma) \vee \tau \big)} \cdot \frac{\tau \sqrt{\log(8 / \delta)}}{\ttr(\Sigma)}
    \\&
    \lesssim \tau \|\Sigma\| \sqrt{\left( 1 \vee \frac{\tau}{\ttr(\Sigma)} \right) \cdot \frac{\log(8 / \delta)}n}.
\end{align*}
\myendproof

\bibliographystyle{abbrvnat}
\bibliography{references}

\begin{thebibliography}{62}
\providecommand{\natexlab}[1]{#1}
\providecommand{\url}[1]{\texttt{#1}}
\expandafter\ifx\csname urlstyle\endcsname\relax
  \providecommand{\doi}[1]{doi: #1}\else
  \providecommand{\doi}{doi: \begingroup \urlstyle{rm}\Url}\fi

\bibitem[Abdalla(2023)]{abdalla23}
P.~Abdalla.
\newblock Covariance estimation under missing observations and {$L_4$-$L_2$}
  moment equivalence.
\newblock Preprint, arXiv:2305.12981, 2023.

\bibitem[Abdalla and Zhivotovskiy(2022)]{abdalla22}
P.~Abdalla and N.~Zhivotovskiy.
\newblock Covariance estimation: Optimal dimension-free guarantees for
  adversarial corruption and heavy tails.
\newblock Preprint, arXiv:2205.08494, 2022.

\bibitem[Adamczak(2015)]{adamczak15}
R.~Adamczak.
\newblock {A note on the Hanson-Wright inequality for random vectors with
  dependencies}.
\newblock \emph{Electronic Communications in Probability}, 20:\penalty0 1--13,
  2015.

\bibitem[Adamczak and Wolff(2015)]{adamczak15b}
R.~Adamczak and P.~Wolff.
\newblock Concentration inequalities for non-{L}ipschitz functions with bounded
  derivatives of higher order.
\newblock \emph{Probability Theory and Related Fields}, 162\penalty0
  (3):\penalty0 531--586, 2015.

\bibitem[Adamczak et~al.(2010)Adamczak, Litvak, Pajor, and
  Tomczak-Jaegermann]{adamczak10}
R.~Adamczak, A.~E. Litvak, A.~Pajor, and N.~Tomczak-Jaegermann.
\newblock Quantitative estimates of the convergence of the empirical covariance
  matrix in log-concave ensembles.
\newblock \emph{Journal of the American Mathematical Society}, 23\penalty0
  (2):\penalty0 535--561, 2010.

\bibitem[Adamczak et~al.(2011)Adamczak, Litvak, Pajor, and
  Tomczak-Jaegermann]{adamczak11}
R.~Adamczak, A.~E. Litvak, A.~Pajor, and N.~Tomczak-Jaegermann.
\newblock Sharp bounds on the rate of convergence of the empirical covariance
  matrix.
\newblock \emph{Comptes Rendus Math\'{e}matique. Acad\'{e}mie des Sciences.
  Paris}, 349\penalty0 (3-4):\penalty0 195--200, 2011.

\bibitem[Bai and Shi(2011)]{bai11}
J.~Bai and S.~Shi.
\newblock Estimating high dimensional covariance matrices and its applications.
\newblock \emph{Annals of Economics and Finance}, 12:\penalty0 199--215, 2011.

\bibitem[Bakry et~al.(2014)Bakry, Gentil, and Ledoux]{bakry14}
D.~Bakry, I.~Gentil, and M.~Ledoux.
\newblock \emph{Analysis and geometry of {M}arkov diffusion operators}, volume
  348 of \emph{Grundlehren der mathematischen Wissenschaften (Fundamental
  Principles of Mathematical Sciences)}.
\newblock Springer, Cham, 2014.

\bibitem[Bickel and Levina(2008{\natexlab{a}})]{bickel08}
P.~J. Bickel and E.~Levina.
\newblock Covariance regularization by thresholding.
\newblock \emph{The Annals of Statistics}, 36\penalty0 (6):\penalty0
  2577--2604, 2008{\natexlab{a}}.

\bibitem[Bickel and Levina(2008{\natexlab{b}})]{bickel08b}
P.~J. Bickel and E.~Levina.
\newblock Regularized estimation of large covariance matrices.
\newblock \emph{The Annals of Statistics}, 36\penalty0 (1):\penalty0 199--227,
  2008{\natexlab{b}}.

\bibitem[Boucheron et~al.(2013)Boucheron, Lugosi, and Massart]{boucheron13}
S.~Boucheron, G.~Lugosi, and P.~Massart.
\newblock \emph{{Concentration Inequalities: A Nonasymptotic Theory of
  Independence}}.
\newblock Oxford University Press, 2013.

\bibitem[Bourgain(1999)]{bourgain99}
J.~Bourgain.
\newblock Random points in isotropic convex sets.
\newblock In \emph{Convex geometric analysis ({B}erkeley, {CA}, 1996)},
  volume~34 of \emph{Mathematical Sciences Research Institute Publications},
  pages 53--58. Cambridge University Press, Cambridge, 1999.

\bibitem[Bunea and Xiao(2015)]{bunea15}
F.~Bunea and L.~Xiao.
\newblock {On the sample covariance matrix estimator of reduced effective rank
  population matrices, with applications to fPCA}.
\newblock \emph{Bernoulli}, 21\penalty0 (2):\penalty0 1200--1230, 2015.

\bibitem[Cai and Zhou(2012)]{cai12}
T.~T. Cai and H.~H. Zhou.
\newblock {Optimal rates of convergence for sparse covariance matrix
  estimation}.
\newblock \emph{The Annals of Statistics}, 40\penalty0 (5):\penalty0
  2389--2420, 2012.

\bibitem[Cai et~al.(2010)Cai, Zhang, and Zhou]{cai10}
T.~T. Cai, C.-H. Zhang, and H.~H. Zhou.
\newblock Optimal rates of convergence for covariance matrix estimation.
\newblock \emph{The Annals of Statistics}, 38\penalty0 (4):\penalty0
  2118--2144, 2010.

\bibitem[Cai et~al.(2013)Cai, Ren, and Zhou]{cai13}
T.~T. Cai, Z.~Ren, and H.~H. Zhou.
\newblock Optimal rates of convergence for estimating {T}oeplitz covariance
  matrices.
\newblock \emph{Probability Theory and Related Fields}, 156\penalty0
  (1):\penalty0 101--143, 2013.

\bibitem[Cai et~al.(2016)Cai, Ren, and Zhou]{cai16}
T.~T. Cai, Z.~Ren, and H.~H. Zhou.
\newblock Estimating structured high-dimensional covariance and precision
  matrices: Optimal rates and adaptive estimation.
\newblock \emph{Electronic Journal of Statistics}, 10\penalty0 (1):\penalty0
  1--59, 2016.

\bibitem[Dahmen et~al.(2000)Dahmen, Keysers, Pitz, and Ney]{jorg00}
J.~Dahmen, D.~Keysers, M.~Pitz, and H.~Ney.
\newblock Structured covariance matrices for statistical image object
  recognition.
\newblock In \emph{Mustererkennung 2000, 22. DAGM-Symposium}, pages 99--106.
  Springer-Verlag, 2000.

\bibitem[Fan et~al.(2015)Fan, Rigollet, and Wang]{fan15}
J.~Fan, P.~Rigollet, and W.~Wang.
\newblock {Estimation of functionals of sparse covariance matrices}.
\newblock \emph{The Annals of Statistics}, 43\penalty0 (6):\penalty0
  2706--2737, 2015.

\bibitem[Giannopoulos et~al.(2005)Giannopoulos, Hartzoulaki, and
  Tsolomitis]{giannopoulos05}
A.~Giannopoulos, M.~Hartzoulaki, and A.~Tsolomitis.
\newblock Random points in isotropic unconditional convex bodies.
\newblock \emph{Journal of the London Mathematical Society}, 72\penalty0
  (3):\penalty0 779--798, 2005.

\bibitem[G{\"o}tze et~al.(2021)G{\"o}tze, Sambale, and Sinulis]{gotze21}
F.~G{\"o}tze, H.~Sambale, and A.~Sinulis.
\newblock {Concentration inequalities for polynomials in
  $\alpha$-sub-exponential random variables}.
\newblock \emph{Electronic Journal of Probability}, 26:\penalty0 1--22, 2021.

\bibitem[Haghighatshoar and Caire(2018)]{haghighatshoar18}
S.~Haghighatshoar and G.~Caire.
\newblock Low-complexity massive mimo subspace estimation and tracking from
  low-dimensional projections.
\newblock \emph{IEEE Transactions on Signal Processing}, 66\penalty0
  (7):\penalty0 1832--1844, 2018.

\bibitem[Han(2022)]{han22}
Q.~Han.
\newblock Exact spectral norm error of sample covariance.
\newblock Preprint, arXiv:2207.13594, 2022.

\bibitem[Han and Wu(2020)]{han20}
Y.~Han and W.~B. Wu.
\newblock {Test for high dimensional covariance matrices}.
\newblock \emph{The Annals of Statistics}, 48\penalty0 (6):\penalty0
  3565--3588, 2020.

\bibitem[Hanson and Wright(1971)]{hanson71}
D.~L. Hanson and F.~T. Wright.
\newblock {A Bound on Tail Probabilities for Quadratic Forms in Independent
  Random Variables}.
\newblock \emph{The Annals of Mathematical Statistics}, 42\penalty0
  (3):\penalty0 1079--1083, 1971.

\bibitem[Hero and Rajaratnam(2012)]{hero12}
A.~Hero and B.~Rajaratnam.
\newblock Hub discovery in partial correlation graphs.
\newblock \emph{IEEE Transactions on Information Theory}, 58\penalty0
  (9):\penalty0 6064--6078, 2012.

\bibitem[Holtz(2010)]{holtz10}
M.~Holtz.
\newblock Sparse grid quadrature in high dimensions with applications in
  finance and insurance.
\newblock In \emph{Lecture Notes in Computational Science and Engineering},
  2010.

\bibitem[Hsu et~al.(2012)Hsu, Kakade, and Zhang]{hsu12}
D.~Hsu, S.~Kakade, and T.~Zhang.
\newblock {A tail inequality for quadratic forms of subgaussian random
  vectors}.
\newblock \emph{Electronic Communications in Probability}, 17:\penalty0 1--6,
  2012.

\bibitem[Kannan et~al.(1997)Kannan, Lov\'{a}sz, and Simonovits]{kannan97}
R.~Kannan, L.~Lov\'{a}sz, and M.~Simonovits.
\newblock Random walks and an ${O}^*({N}^5)$ volume algorithm for convex
  bodies.
\newblock \emph{Random Structures \& Algorithms}, 11\penalty0 (1):\penalty0
  1–--50, 1997.

\bibitem[Klochkov and Zhivotovskiy(2020)]{klochkov20}
Y.~Klochkov and N.~Zhivotovskiy.
\newblock {Uniform Hanson-Wright type concentration inequalities for unbounded
  entries via the entropy method}.
\newblock \emph{Electronic Journal of Probability}, 25:\penalty0 1--30, 2020.

\bibitem[Koltchinskii and Lounici(2017)]{koltchinskii17}
V.~Koltchinskii and K.~Lounici.
\newblock {Concentration inequalities and moment bounds for sample covariance
  operators}.
\newblock \emph{Bernoulli}, 23\penalty0 (1):\penalty0 110--133, 2017.

\bibitem[Krim and Viberg(1996)]{krim96}
H.~Krim and M.~Viberg.
\newblock Two decades of array signal processing research: the parametric
  approach.
\newblock \emph{IEEE Signal Processing Magazine}, 13\penalty0 (4):\penalty0
  67--94, 1996.

\bibitem[Laurent and Massart(2000)]{laurent00}
B.~Laurent and P.~Massart.
\newblock {Adaptive estimation of a quadratic functional by model selection}.
\newblock \emph{The Annals of Statistics}, 28\penalty0 (5):\penalty0
  1302--1338, 2000.

\bibitem[Ledoit and Wolf(2002)]{ledoit02}
O.~Ledoit and M.~Wolf.
\newblock {Some hypothesis tests for the covariance matrix when the dimension
  is large compared to the sample size}.
\newblock \emph{The Annals of Statistics}, 30\penalty0 (4):\penalty0
  1081--1102, 2002.

\bibitem[Ledoit and Wolf(2003)]{ledoit03}
O.~Ledoit and M.~Wolf.
\newblock Improved estimation of the covariance matrix of stock returns with an
  application to portfolio selection.
\newblock \emph{Journal of Empirical Finance}, 10\penalty0 (5):\penalty0
  603--621, 2003.

\bibitem[Leng and Pan(2018)]{leng18}
C.~Leng and G.~Pan.
\newblock {Covariance estimation via sparse Kronecker structures}.
\newblock \emph{Bernoulli}, 24\penalty0 (4B):\penalty0 3833--3863, 2018.

\bibitem[Li and Chen(2012)]{li12}
J.~Li and S.~X. Chen.
\newblock {Two sample tests for high-dimensional covariance matrices}.
\newblock \emph{The Annals of Statistics}, 40\penalty0 (2):\penalty0 908--940,
  2012.

\bibitem[Lounici(2014)]{lounici14}
K.~Lounici.
\newblock {High-dimensional covariance matrix estimation with missing
  observations}.
\newblock \emph{Bernoulli}, 20\penalty0 (3):\penalty0 1029--1058, 2014.

\bibitem[Lounici and Pacreau(2023)]{lounici23}
K.~Lounici and G.~Pacreau.
\newblock Robust covariance estimation with missing values and cell-wise
  contamination.
\newblock Preprint, arXiv:2306.00752, 2023.

\bibitem[Mendelson and Zhivotovskiy(2020)]{mendelson20}
S.~Mendelson and N.~Zhivotovskiy.
\newblock {Robust covariance estimation under $L_{4}-L_{2}$ norm equivalence}.
\newblock \emph{The Annals of Statistics}, 48\penalty0 (3):\penalty0
  1648--1664, 2020.

\bibitem[Minasyan and Zhivotovskiy(2023)]{minasyan23}
A.~Minasyan and N.~Zhivotovskiy.
\newblock Statistically optimal robust mean and covariance estimation for
  anisotropic gaussians.
\newblock Preprint, arXiv:2301.09024, 2023.

\bibitem[Oliveira(2016)]{oliveira16}
R.~I. Oliveira.
\newblock The lower tail of random quadratic forms with applications to
  ordinary least squares.
\newblock \emph{Probability Theory and Related Fields}, 166\penalty0
  (3):\penalty0 1175--1194, 2016.

\bibitem[Paouris(2006)]{paouris06}
G.~Paouris.
\newblock Concentration of mass on convex bodies.
\newblock \emph{Geometric and Functional Analysis}, 16\penalty0 (5):\penalty0
  1021--1049, 2006.

\bibitem[Rudelson(1999)]{rudelson99}
M.~Rudelson.
\newblock Random vectors in the isotropic position.
\newblock \emph{Journal of Functional Analysis}, 164\penalty0 (1):\penalty0
  60--72, 1999.

\bibitem[Rudelson and Vershynin(2013)]{rudelson13}
M.~Rudelson and R.~Vershynin.
\newblock {Hanson-Wright inequality and sub-gaussian concentration}.
\newblock \emph{Electronic Communications in Probability}, 18:\penalty0 1--9,
  2013.

\bibitem[Sambale(2023)]{sambale23}
H.~Sambale.
\newblock Some notes on concentration for $\alpha$-subexponential random
  variables.
\newblock In \emph{High Dimensional Probability IX}, pages 167--192. Springer
  International Publishing, 2023.

\bibitem[Sch{\"a}fer and Strimmer(2005)]{schafer05}
J.~Sch{\"a}fer and K.~Strimmer.
\newblock A shrinkage approach to large-scale covariance matrix estimation and
  implications for functional genomics.
\newblock \emph{Statistical applications in genetics and molecular biology},
  4:\penalty0 article 32, 2005.

\bibitem[Schudy and Sviridenko(2012)]{schudy12}
W.~Schudy and M.~Sviridenko.
\newblock Concentration and moment inequalities for polynomials of independent
  random variables.
\newblock In \emph{Proceedings of the Twenty-Third Annual ACM-SIAM Symposium on
  Discrete Algorithms}, SODA '12, pages 437--446. Society for Industrial and
  Applied Mathematics, 2012.

\bibitem[Song Xi~Chen and Zhong(2010)]{chen10}
L.-X.~Z. Song Xi~Chen and P.-S. Zhong.
\newblock Tests for high-dimensional covariance matrices.
\newblock \emph{Journal of the American Statistical Association}, 105\penalty0
  (490):\penalty0 810--819, 2010.

\bibitem[Spokoiny(2023)]{spokoiny23}
V.~Spokoiny.
\newblock Sharp deviation bounds and concentration phenomenon for the squared
  norm of a sub-gaussian vector.
\newblock Preprint, arXiv:2305.07885, 2023.

\bibitem[Srivastava and Vershynin(2013)]{srivastava13}
N.~Srivastava and R.~Vershynin.
\newblock Covariance estimation for distributions with $2+\varepsilon$ moments.
\newblock \emph{The Annals of Probability}, 41:\penalty0 3081--3111, 2013.

\bibitem[Tikhomirov(2018)]{tikhomirov18}
K.~Tikhomirov.
\newblock Sample covariance matrices of heavy-tailed distributions.
\newblock \emph{International Mathematics Research Notices}, 2018\penalty0
  (20):\penalty0 6254--6289, 2018.

\bibitem[Tropp(2012)]{tropp12}
J.~A. Tropp.
\newblock User-friendly tail bounds for sums of random matrices.
\newblock \emph{Foundations of Computational Mathematics}, 12\penalty0
  (4):\penalty0 389--434, 2012.

\bibitem[Tsiligkaridis and Hero(2013)]{tsiligkaridis13}
T.~Tsiligkaridis and A.~O. Hero.
\newblock Covariance estimation in high dimensions via {K}ronecker product
  expansions.
\newblock \emph{IEEE Transactions on Signal Processing}, 61\penalty0
  (21):\penalty0 5347--5360, 2013.

\bibitem[Vershynin(2012{\natexlab{a}})]{vershynin12}
R.~Vershynin.
\newblock How close is the sample covariance matrix to the actual covariance
  matrix?
\newblock \emph{Journal of Theoretical Probability}, 25\penalty0 (3):\penalty0
  655--686, 2012{\natexlab{a}}.

\bibitem[Vershynin(2012{\natexlab{b}})]{vershynin12b}
R.~Vershynin.
\newblock Introduction to the non-asymptotic analysis of random matrices.
\newblock In \emph{Compressed Sensing: Theory and Applications}, pages
  210--268, Cambridge, 2012{\natexlab{b}}. Cambridge University Press.

\bibitem[Vershynin(2018)]{vershynin18}
R.~Vershynin.
\newblock \emph{High-Dimensional Probability: An Introduction with Applications
  in Data Science}.
\newblock Cambridge Series in Statistical and Probabilistic Mathematics.
  Cambridge University Press, 2018.

\bibitem[Xiao and Wu(2012)]{xiao12}
H.~Xiao and W.~B. Wu.
\newblock Covariance matrix estimation for stationary time series.
\newblock \emph{The Annals of Statistics}, 40\penalty0 (1):\penalty0 466--493,
  2012.

\bibitem[Xie and Bentler(2003)]{xie03}
J.~Xie and P.~M. Bentler.
\newblock Covariance structure models for gene expression microarray data.
\newblock \emph{Structural Equation Modeling: A Multidisciplinary Journal},
  10\penalty0 (4):\penalty0 566--582, 2003.

\bibitem[Youssef(2013)]{youssef13}
P.~Youssef.
\newblock {Estimating the covariance of random matrices}.
\newblock \emph{Electronic Journal of Probability}, 18:\penalty0 1--26, 2013.

\bibitem[Zhang and Schneider(2010)]{zhang10}
Y.~Zhang and J.~Schneider.
\newblock Learning multiple tasks with a sparse matrix-normal penalty.
\newblock In \emph{Advances in Neural Information Processing Systems},
  volume~23. Curran Associates, Inc., 2010.

\bibitem[Zhivotovskiy(2024)]{zhivotovskiy21}
N.~Zhivotovskiy.
\newblock {Dimension-free bounds for sums of independent matrices and simple
  tensors via the variational principle}.
\newblock \emph{Electronic Journal of Probability}, 29:\penalty0 1--28, 2024.

\end{thebibliography}

\appendix

\section{Auxiliary results: Kronecker product and its properties}
\label{sec:kronecker}

In this section, we present some facts about the Kronecker product, which are extensively used in our proofs.
Let us recall that, for any matrices $A \in \R^{p \times q}$ and $B \in \R^{r \times s}$, their Kronecker product $A \otimes B$ is a matrix of size $pr \times qs$, defined as
\[
    A \otimes B =
    \begin{pmatrix}
        a_{11} B & \dots & a_{1q} B \\
        \vdots & \ddots & \vdots \\
        a_{p1} B & \dots & a_{pq} B 
    \end{pmatrix}.
\]
The basic property of the Kronecker product is the following identity (also referred to as the mixed product property):
\begin{equation}
    \label{eq:mixed_product}
    (A \otimes B) (C \otimes D) = (AC) \otimes (BD).
\end{equation}
This immediately implies that
\[
    (A \otimes B)^{-1} = A^{-1} \otimes B^{-1}
    \quad \text{and} \quad
    (A \otimes B)^k = (A^k) \otimes (B^k)
    \quad \text{for all $k \in \mathbb N$.}
\]
Another corollary of \eqref{eq:mixed_product} is that if $A \in \R^{p \times p}$ and $B \in \R^{r \times r}$ are symmetric matrices with eigenvalues $\lambda_1 \geq \lambda_2 \geq \dots \geq \lambda_p$ and $\mu_1 \geq \mu_2 \geq \dots \geq \mu_r$, respectively, then
\[
    \left\{ \lambda_i \mu_j : 1 \leq i \leq p, \, 1 \leq j \leq r \right\}
\]
are the eigenvalues of $(A \otimes B)$. Indeed, if $A = U^\top \Lambda U$ and $B = V^\top M V$ are the eigenvalue decompositions of $A$ and $B$, respectively, then
\[
    A \otimes B
    = (U^\top \Lambda U) \otimes (V^\top M V)
    = (U^\top \otimes V^\top) \big[ (\Lambda U) \otimes (M V) \big]
    = (U \otimes V)^\top (\Lambda \otimes M) (U \otimes V)
\]
is the eigenvalue decomposition of $(A \otimes B)$.
This fact yields that
\begin{equation}
    \label{eq:kronecker_eigenvalues}
    \Tr(A \otimes B) = \Tr(A) \Tr(B),
    \quad
    \|A \otimes B\| = \|A\| \|B\|,
    \quad
    \text{and}
    \quad
    \|A \otimes B\|_{\Fr} = \|A\|_{\Fr} \|B\|_{\Fr}.
\end{equation}

We also extensively exploit the properties of the Kronecker product related to interaction with the vectorization operator $\rmvec$ (see our notation in Section \ref{sec:introduction}). Let $A \in \R^{p \times q}$ and $B \in \R^{r \times s}$ be arbitrary matrices (not necessarily symmetric) and let $U \in \R^{s \times q}$ be a matrix with columns $\bu_1, \dots, \bu_q$, so we have $\rmvec(U) = (\bu_1^\top, \dots, \bu_q^\top)^\top \in \R^{sq}$. Then it is straightforward to check that
\begin{equation}
    \label{eq:kronecker_vec}
    (A \otimes B) \rmvec(U) = \rmvec(B U A^\top).
\end{equation}
As a consequence, for any matrices $A \in \R^{p \times q}$, $B \in \R^{r \times s}$, $U \in \R^{s \times q}$, and $V \in \R^{r \times p}$, it holds that
\begin{equation}
    \label{eq:kronecker_trace}
    \Tr(V^\top B U A^\top) = \rmvec(V)^\top \rmvec(B U A^\top) = \rmvec(V)^\top (A \otimes B) \rmvec(U).
\end{equation}
In our proofs, we use the representation \eqref{eq:kronecker_trace} to compute expectation the bilinear (in $U$ and $V$) form $\Tr(V^\top B U A^\top)$ when $U$ and $V$ are random matrices.

\section{Proof of Lemma \ref{lem:linearization}}
\label{sec:lem_linearization_proof}

Since all the elements of $\Gamma$ are i.i.d. standard Gaussian random variables, the random vector $\rmvec(\Sigma^{1/2} \Gamma \Sigma^{1/2})$ has a Gaussian distribution $\cN(\bzero, \Sigma \otimes \Sigma)$ (see, for instance, the paper of \cite{leng18} for Gaussian matrix models). Then it holds that
\begin{align*}
    &
    \E_\Gamma \, e^{\sqrt{\lambda} \Tr(\Gamma^\top \Sigma^{1/2} A \Sigma^{1/2})} \1\left( \|\Sigma^{1/2} \Gamma \Sigma^{1/2}\|_{\Fr} \leq \z \sqrt{\lambda} +  \sqrt{\beta} \Tr(\Sigma) \right)
    \\&
    = C_\Sigma \int\limits_{\|y\| \leq \z \sqrt{\lambda} +  \sqrt{\beta} \Tr(\Sigma)} e^{\sqrt{\lambda} y^\top \rmvec(A) - y^\top (\Sigma \otimes \Sigma)^{-1} y / 2} \rmd y,
\end{align*}
where $C_\Sigma = (2\pi)^{-d^2/2} \det(\Sigma)^{-1}$ is the normalizing constant. Using the equality
\begin{align*}
    &
    \sqrt{\lambda} y^\top \rmvec(A) - \frac{y^\top (\Sigma \otimes \Sigma)^{-1} y}2
    \\&
    = \frac{\lambda \rmvec(A)^\top (\Sigma \otimes \Sigma) \rmvec(A)}2 - \frac12 \|(\Sigma^{1/2} \otimes \Sigma^{1/2})^{-1} y - \sqrt{\lambda} (\Sigma \otimes \Sigma)^{1/2} \rmvec(A)\|^2
\end{align*}
and making the substitution
\[
    x = y - \sqrt{\lambda} (\Sigma \otimes \Sigma) \rmvec(A)
\]
under the integral, we obtain that
\begin{align*}
    &
    \E_\Gamma \, e^{\sqrt{\lambda} \Tr(\Gamma^\top \Sigma^{1/2} A \Sigma^{1/2})} \1\left( \|\Sigma^{1/2} \Gamma \Sigma^{1/2}\|_{\Fr} \leq \z \sqrt{\lambda} +  \sqrt{\beta} \Tr(\Sigma) \right)
    \\&
    = C_\Sigma \int\limits_{\|x + \sqrt{\lambda} (\Sigma \otimes \Sigma) \rmvec(A) \| \leq \z \sqrt{\lambda} +  \sqrt{\beta} \Tr(\Sigma)} e^{\lambda \Tr(A^\top \Sigma A \Sigma) / 2 - x^\top (\Sigma \otimes \Sigma)^{-1} x / 2} \rmd x
    \\&
    = e^{\lambda \|\Sigma^{1/2} A \Sigma^{1/2}\|_{\Fr}^2 / 2} \p\left( \|\Sigma^{1/2} \Gamma \Sigma^{1/2} + \sqrt{\lambda} \Sigma A \Sigma\|_{\Fr} \leq \z \sqrt{\lambda} +  \sqrt{\beta} \Tr(\Sigma) \right).
\end{align*}
Due to the property \eqref{eq:kronecker_vec} of the Kronecker product, it holds that
\[
    \|(\Sigma \otimes \Sigma) \rmvec(A) \|
    = \|\rmvec(\Sigma A \Sigma)\|
    = \|\Sigma A \Sigma\|_{\Fr} \leq \z,
\]
where the last inequality is guaranteed by the conditions of the lemma.
This yields that the ball $\cB\big((\Sigma \otimes \Sigma) \rmvec(A), \z \sqrt{\lambda} +  \sqrt{\beta} \Tr(\Sigma) \big)$ contains $\cB\big(\bzero, \sqrt{\beta} \Tr(\Sigma) \big)$ and, hence,
\begin{align*}
    &
    e^{\lambda \|\Sigma^{1/2} A \Sigma^{1/2}\|_{\Fr}^2 / 2} \; \p\left( \|\Sigma^{1/2} \Gamma \Sigma^{1/2} + \sqrt{\lambda} \Sigma A \Sigma\|_{\Fr} \leq \z \sqrt{\lambda} +  \sqrt{\beta} \Tr(\Sigma) \right)
    \\&
    \geq e^{\lambda \|\Sigma^{1/2} A \Sigma^{1/2}\|_{\Fr}^2 / 2} \; \p\left( \|\Sigma^{1/2} \Gamma \Sigma^{1/2}\|_{\Fr} \leq \z \sqrt{\lambda} +  \sqrt{\beta} \Tr(\Sigma) - \sqrt{\lambda} \|\Sigma A \Sigma\|_{\Fr} \right)
    \\&
    \geq e^{\lambda \|\Sigma^{1/2} A \Sigma^{1/2}\|_{\Fr}^2 / 2} \; \p\left( \|\Sigma^{1/2} \Gamma \Sigma^{1/2}\|_{\Fr} \leq \sqrt{\beta} \Tr(\Sigma) \right).
\end{align*}
It only remains to note that, due to \eqref{eq:kronecker_eigenvalues} and the Markov inequality, we have
\[
    \p\left( \|\Sigma^{1/2} \Gamma \Sigma^{1/2}\|_{\Fr} \geq \sqrt{\beta} \Tr(\Sigma) \right)
    \leq \frac{\E \|\Sigma^{1/2} \Gamma \Sigma^{1/2}\|_{\Fr}^2}{\beta \Tr(\Sigma)^2}
    = \frac{\Tr(\Sigma \otimes \Sigma)}{\beta \Tr(\Sigma)^2} = \frac1\beta,
\]
and then
\[
    e^{\lambda \|\Sigma^{1/2} A \Sigma^{1/2}\|_{\Fr}^2 / 2} \p\left( \|\Sigma^{1/2} \Gamma \Sigma^{1/2}\|_{\Fr} \leq \sqrt{\beta} \Tr(\Sigma) \right)
    \geq \frac{\beta - 1}\beta e^{\lambda \|\Sigma^{1/2} A \Sigma^{1/2}\|_{\Fr}^2 / 2}.
\]
\myendproof

\section{Proof of Lemma \ref{lem:derivatives}}
\label{sec:derivatives}

Let us recall that, for any Borel function $f : \R^d \times \mathbb \R^{d \times d} \rightarrow \mathbb \R$, $\sfP_{U} f(\bxi, U)$ denotes the expectation
\[
	\sfP_{U} f(\bxi, U) = \frac{\E f(\bxi, U) e^{\bxi^\top U \bxi}}{\E e^{\bxi^\top U \bxi}}.
\]
One can compute the derivative of $\sfP_U f(\bxi, U)$ according to the formula
\begin{align}
    \label{eq:f_derivative}
    \left\langle \nabla \sfP_U f(\bxi, U), V \right\rangle
    &\notag
    \equiv \Tr\left( \nabla \sfP_U f(\bxi, U)^\top V \right)
    \\&
    = \sfP_U \left\langle \frac{\partial f(\bxi, U)}{\partial U}, V \right\rangle + \sfP_U f(\bxi, U) \, \bxi^\top V \bxi - \sfP_U f(\bxi, U) \; \sfP_U \bxi^\top V \bxi
    \\&\notag
    = \sfP_U \left\langle \frac{\partial f(\bxi, U)}{\partial U}, V \right\rangle + \sfP_U \Big( f(\bxi, U) - \sfP_U f(\bxi, U) \Big) \left(\bxi^\top V \bxi - \sfP_U \bxi^\top V \bxi \right),
\end{align}
which holds for all $U, V \in \R^{d \times d}$ whenever the expectations in the left- and the right-hand sides are well-defined. It is straightforward to check that
\[
    \left\langle \nabla \varphi(U), V \right\rangle = \Tr \left(\nabla \varphi(U)^\top V \right) = \sfP_U \left( \bxi^\top V \bxi \right).
\]
Then, using \eqref{eq:f_derivative}, we can compute
\[
    \left\langle \nabla^2 \varphi(U), V^{\otimes 2} \right\rangle
    \equiv \Big\langle \nabla \left\langle \nabla \varphi(U), V \right\rangle, V \Big\rangle
    = \sfP_U \left(\bxi^\top V \bxi - \sfP_U \bxi^\top V \bxi \right)^2.
\]
For the higher order derivatives, we similarly have
\begin{align*}
    \left\langle \nabla^3 \varphi(U), V^{\otimes 3} \right\rangle
    &
    = -2 \sfP_U \left(\bxi^\top V \bxi - \sfP_U \bxi^\top V \bxi \right) \left\langle \frac{\rmd}{\rmd U} \sfP_U \bxi^\top V \bxi, V \right \rangle + \sfP_U \left(\bxi^\top V \bxi - \sfP_U \bxi^\top V \bxi \right)^3
    \\&
    = 0 + \sfP_U \left(\bxi^\top V \bxi - \sfP_U \bxi^\top V \bxi \right)^3
    = \sfP_U \left(\bxi^\top V \bxi - \sfP_U \bxi^\top V \bxi \right)^3
\end{align*}
and
\begin{align*}
    \left\langle \nabla^4 \varphi(U), V^{\otimes 4} \right\rangle
    &
    = -3 \sfP_U \left(\bxi^\top V \bxi - \sfP_U \bxi^\top V \bxi \right)^2 \left\langle \frac{\rmd}{\rmd U} \sfP_U \bxi^\top V \bxi, V \right \rangle + \sfP_U \left(\bxi^\top V \bxi - \sfP_U \bxi^\top V \bxi \right)^4
    \\&
    = -3 \sfP_U \left(\bxi^\top V \bxi - \sfP_U \bxi^\top V \bxi \right)^2 \left\langle \nabla^2 \varphi(U), V^{\otimes 2} \right \rangle + \sfP_U \left(\bxi^\top V \bxi - \sfP_U \bxi^\top V \bxi \right)^4
    \\&
    = \sfP_U \left(\bxi^\top V \bxi - \sfP_U \bxi^\top V \bxi \right)^4 - 3 \left( \sfP_U \left(\bxi^\top V \bxi - \sfP_U \bxi^\top V \bxi \right)^2 \right)^2.
\end{align*}
According to Assumption \ref{as:derivatives}, we have
\[
    \sfP_U \left(\bxi^\top V \bxi - \sfP_U \bxi^\top V \bxi \right)^4 \leq \tau^2 \; \|V\|_{\Fr}^4
\]
for all $V \in \R^{d \times d}$ and all $U \in \R^{d \times d}$ such that $\|U\|_{\Fr} \leq \rho_{\max}$. Then, due to the H\"{o}lder inequality, we have
\[
    0
    \leq \left\langle \nabla^2 \varphi(U), V^{\otimes 2} \right\rangle
    = \sfP_U \left(\bxi^\top V \bxi - \sfP_U \bxi^\top V \bxi \right)^2
    \leq \left( \sfP_U \left(\bxi^\top V \bxi - \sfP_U \bxi^\top V \bxi \right)^4 \right)^{1/2}
    \leq \tau \|V\|_{\Fr}^2
\]
and
\[
    \left| \left\langle \nabla^3 \varphi(U), V^{\otimes 3} \right\rangle
    \right|
    = \left| \sfP_U \left(\bxi^\top V \bxi - \sfP_U \bxi^\top V \bxi \right)^3 \right|
    \leq \left( \sfP_U \left(\bxi^\top V \bxi - \sfP_U \bxi^\top V \bxi \right)^4 \right)^{3/4}
    \leq \tau^{3/2} \; \|V\|_{\Fr}^3.
\]
Finally, for the fourth derivative, using the Cauchy-Schwarz inequality, we obtain that
\begin{align*}
    \left\langle \nabla^4 \varphi(U), V^{\otimes 4} \right\rangle
    &
    = \sfP_U \left(\bxi^\top V \bxi - \sfP_U \bxi^\top V \bxi \right)^4 - 3 \left( \sfP_U \left(\bxi^\top V \bxi - \sfP_U \bxi^\top V \bxi \right)^2 \right)^2
    \\&
    \geq -2 \sfP_U \left(\bxi^\top V \bxi - \sfP_U \bxi^\top V \bxi \right)^4
    \geq -2 \tau^2 \|V\|_{\Fr}^4
\end{align*}
and
\begin{align*}
    \left\langle \nabla^4 \varphi(U), V^{\otimes 4} \right\rangle
    &
    = \sfP_U \left(\bxi^\top V \bxi - \sfP_U \bxi^\top V \bxi \right)^4 - 3 \left( \sfP_U \left(\bxi^\top V \bxi - \sfP_U \bxi^\top V \bxi \right)^2 \right)^2
    \\&
    \leq \sfP_U \left(\bxi^\top V \bxi - \sfP_U \bxi^\top V \bxi \right)^4
    \leq \tau^2 \|V\|_{\Fr}^4
\end{align*}
for any $V \in \R^{d \times d}$ and any $U \in \R^{d \times d}$, such that $\|U\|_{\Fr} \leq \rho_{\max}$.

\myendproof

\section{Proof of Lemma \ref{lem:restricted_moment}}
\label{sec:lem_restricted_moment_proof}

Let $\eps$ be a Rademacher random variable, which is independent of $\bX_1, \dots, \bX_n$ and of $\Gamma$. The symmetry of the Gaussian distribution yields that
\begin{align*}
    &
    \E e^{\Phi(\sqrt{\lambda} \Sigma^{1/2} \Gamma \Sigma^{1/2})} \1\left( \|\Sigma^{1/2} \Gamma \Sigma^{1/2}\|_{\Fr} \leq \z \sqrt{\lambda} +  \sqrt{2} \, \Tr(\Sigma) \right)
    \\&
    = \E e^{\Phi(\eps\sqrt{\lambda} \Sigma^{1/2} \Gamma \Sigma^{1/2})} \1\left( \|\Sigma^{1/2} \Gamma \Sigma^{1/2}\|_{\Fr} \leq \z \sqrt{\lambda} +  \sqrt{2} \, \Tr(\Sigma) \right).
\end{align*}
Due to the conditions of the theorem, we can use Taylor's expansion of $\Phi(\eps \sqrt{\lambda} \Sigma^{1/2} \Gamma \Sigma^{1/2})$ on the event $\{\|\Sigma^{1/2} \Gamma \Sigma^{1/2}\|_{\Fr} \leq \z \sqrt{\lambda} + \sqrt{2} \, \Tr(\Sigma) \}$:
\begin{align}
    \label{eq:taylor}
    &\notag
    \left| \Phi\left(\eps \sqrt{\lambda} \Sigma^{1/2} \Gamma \Sigma^{1/2} \right) - \frac{\lambda}2 \E_{\bxi} \Tr\left((\bxi \bxi^\top - I) \Sigma^{1/2} \Gamma \Sigma^{1/2} \right)^2 - \frac{\lambda^{3/2} \eps}{6 \sqrt{n}} \E_{\bxi} \Tr\left((\bxi \bxi^\top - I) \Sigma^{1/2} \Gamma \Sigma^{1/2} \right)^3 \right|
    \\&
    \leq \frac{\tau^2 \lambda^2}{24 n} \left\|\Sigma^{1/2} \Gamma \Sigma^{1/2} \right\|_{\Fr}^4.
\end{align}
The notation $\E_{\bxi}$ stands for the expectation with respect to $\bxi$ only (that is, conditioned on $\Gamma$).
Moreover, Hoeffding's lemma implies that
\begin{align}
    \label{eq:hoeffding}
    \E_{\eps} \exp\left\{ \frac{\lambda^{3/2} \eps}{6 \sqrt{n}} \E_{\bxi} \Tr\left((\bxi \bxi^\top - I) \Sigma^{1/2} \Gamma \Sigma^{1/2} \right)^3 \right\}
    &\notag
    \leq \exp\left\{ \frac{\lambda^{3} }{72 n} \left[ \E_{\bxi} \Tr\left((\bxi \bxi^\top - I) \Sigma^{1/2} \Gamma \Sigma^{1/2} \right)^3 \right]^2 \right\}
    \\&
    \leq \exp\left\{ \frac{\lambda^{3} \tau^3}{72 n} \left\| \Sigma^{1/2} \Gamma \Sigma^{1/2} \right\|_{\Fr}^6 \right\}.
\end{align}
Combining the inequalities \eqref{eq:taylor} and \eqref{eq:hoeffding}, we obtain that
\begin{align*}
    &
    \E e^{\Phi(\sqrt{\lambda} \Sigma^{1/2} \Gamma \Sigma^{1/2})} \1\left( \|\Sigma^{1/2} \Gamma \Sigma^{1/2}\|_{\Fr} \leq \z \sqrt{\lambda} +  \sqrt{2} \, \Tr(\Sigma) \right)
    \\&
    \leq \E_\Gamma \exp\left\{ \frac{\lambda}2 \E_{\bxi} \left[ \Tr\left((\bxi \bxi^\top - I) \Sigma^{1/2} \Gamma \Sigma^{1/2} \right) \right]^2 + \frac{\lambda^{3} \tau^3}{72 n} \left\| \Sigma^{1/2} \Gamma \Sigma^{1/2} \right\|_{\Fr}^6 + \frac{\tau^2 \lambda^2}{24 n} \left\|\Sigma^{1/2} \Gamma \Sigma^{1/2} \right\|_{\Fr}^4 \right\}
    \\&\quad
    \cdot \1\left( \|\Sigma^{1/2} \Gamma \Sigma^{1/2}\|_{\Fr} \leq \z \sqrt{\lambda} +  \sqrt{2} \Tr(\Sigma) \right)
    \\&
    \leq \E_\Gamma \exp\left\{ \frac{\lambda}2 \E_{\bxi} \left[ \Tr\left( (\bxi \bxi^\top - I) \Sigma^{1/2} \Gamma \Sigma^{1/2} \right) \right]^2 + \frac{1}2 G(\lambda, \z) \; \|\Sigma^{1/2} \Gamma \Sigma^{1/2}\|_{\Fr}^2 \right\}.
\end{align*}
Next, we use the properties \eqref{eq:kronecker_trace} of the Kronecker product to represent the power of the exponent in the following form:
\begin{align*}
    &
    \frac{\lambda}2 \; \E_{\bxi} \left[ \Tr\left( (\bxi \bxi^\top - I) \Sigma^{1/2} \Gamma \Sigma^{1/2} \right) \right]^2 + \frac{1}2 G(\lambda, \z) \|\Sigma^{1/2} \Gamma \Sigma^{1/2}\|_{\Fr}^2
    \\&
    = \frac{\lambda}2 \; \rmvec(\Gamma)^\top \left( \E_{\bxi} \rmvec(\Sigma^{1/2} (\bxi \bxi^\top - I) \Sigma^{1/2}) \rmvec(\Sigma^{1/2} (\bxi \bxi^\top - I) \Sigma^{1/2})^\top \right) \rmvec(\Gamma)
    \\&\quad
    + \frac{1}2 G(\lambda, \z) \; \rmvec(\Gamma)^\top \left( \Sigma \otimes \Sigma \right) \rmvec(\Gamma)
    \\&
    = \frac{\lambda}2 \; \rmvec(\Gamma)^\top \left( \E_{\bX} \rmvec(\bX \bX^\top - \Sigma) \rmvec(\bX \bX^\top - \Sigma)^\top \right) \rmvec(\Gamma)
    \\&\quad
    + \frac{1}2 G(\lambda, \z) \; \rmvec(\Gamma)^\top \left( \Sigma \otimes \Sigma \right) \rmvec(\Gamma).
\end{align*}
Then it holds that
\begin{align}
    \label{eq:gauss_exp_moment}
    &\notag
    \E_\Gamma \exp\left\{ \frac{\lambda}2 \; \E_{\bxi} \left[ \Tr\left( (\bxi \bxi^\top - I) \Sigma^{1/2} \Gamma \Sigma^{1/2} \right) \right]^2 + \frac{1}2 G(\lambda, \z) \|\Sigma^{1/2} \Gamma \Sigma^{1/2}\|_{\Fr}^2 \right\}
    \\&
    = \det\left(I_{d^2} - \lambda \E_{\bX} \rmvec(\bX \bX^\top - \Sigma) \rmvec(\bX \bX^\top - \Sigma)^\top - G(\lambda, \z) (\Sigma \otimes \Sigma) \right)^{-1/2},
\end{align}
provided that the operator norm of
\begin{equation}
    \label{eq:matrix}
    \lambda \; \E_{\bX} \rmvec(\bX \bX^\top - \Sigma) \rmvec(\bX \bX^\top - \Sigma)^\top + G(\lambda, \z) (\Sigma \otimes \Sigma)
\end{equation}
is strictly less than $1$. Note that
\begin{align*}
    &
    \left\| \E_{\bX} \rmvec(\bX \bX^\top - \Sigma) \rmvec(\bX \bX^\top - \Sigma)^\top \right\|
    \\&
    = \sup\limits_{\|U\|_{\Fr} = 1} \rmvec(U)^\top \left( \E \rmvec(\bX \bX^\top - \Sigma) \rmvec(\bX \bX^\top - \Sigma)^\top \right) \rmvec(U)
    \\&
    = \sup\limits_{\|U\|_{\Fr} = 1} \E \left( \rmvec(U)^\top \rmvec(\bX \bX^\top - \Sigma) \right)^2
    \\&
    = \sup\limits_{\|U\|_{\Fr} = 1} \E \left[ \bX^\top U \bX - \Tr(U^\top \Sigma) \right]^2
    = \kappa.
\end{align*}
Then it holds that
\begin{align*}
    &
    \left\| \lambda \; \E_{\bX} \rmvec(\bX \bX^\top - \Sigma) \rmvec(\bX \bX^\top - \Sigma)^\top + G(\lambda, \z) (\Sigma \otimes \Sigma) \right\|
    \\&
    \leq \lambda \left\| \E_{\bX} \rmvec(\bX \bX^\top - \Sigma) \rmvec(\bX \bX^\top - \Sigma)^\top \right\| + G(\lambda, \z) \left\| \Sigma \otimes \Sigma \right\|
    \\&
    = \lambda \kappa + G(\lambda, \z) \|\Sigma\|^2.
\end{align*}
Since, according to the conditions of the lemma,
the right-hand side does not exceed $1/2$, we can apply the inequality
\begin{equation}
    \label{eq:log_det_bound}
    -\log \det(I_{d^2} - A) \leq \Tr(A) + 2 \|A\|_{\Fr}^2,
\end{equation}
which holds for any matrix $A \in \R^{d^2 \times d^2}$, whose maximal eigenvalue does not exceed $1/2$. Note that the squared Frobenius norm of the matrix \eqref{eq:matrix} is not greater than
\begin{align}
    \label{eq:frobenuis_bound}
    &\notag
    \left\| \lambda \; \E \rmvec(\bX \bX^\top - \Sigma) \rmvec(\bX \bX^\top - \Sigma)^\top + G(\lambda, \z) (\Sigma \otimes \Sigma) \right\|_{\Fr}^2
    \\&
    \leq 2 \lambda^2 \left\| \E \rmvec(\bX \bX^\top - \Sigma) \rmvec(\bX \bX^\top - \Sigma)^\top \right\|_{\Fr}^2
    + 2 G(\lambda, \z)^2 \left\| \Sigma \otimes \Sigma \right\|_{\Fr}^2
    \\&\notag
    = 2 \lambda^2 \ttv^2 + 2 G(\lambda, \z)^2 \left\| \Sigma \right\|_{\Fr}^4.
\end{align}
Taking into account \eqref{eq:gauss_exp_moment}, \eqref{eq:log_det_bound}, and \eqref{eq:frobenuis_bound}, we finally obtain that
\begin{align*}
    &
    \E_\Gamma \exp\left\{ \frac{\lambda}2 \E_H \Tr(H^\top \Sigma^{1/2} \Gamma \Sigma^{1/2})^2 + \frac{1}2 G(\lambda, \z) \|\Sigma^{1/2} \Gamma \Sigma^{1/2}\|_{\Fr}^2 \right\}
    \\&
    \leq \exp\left\{ \frac{\lambda \ttm}2 + 2 \lambda^2 \ttv^2 + \frac12 G(\lambda, \z) \Tr(\Sigma)^2 + 2 G(\lambda, \z)^2 \left\| \Sigma \right\|_{\Fr}^4 \right\}
\end{align*}
for any $\lambda > 0$, which fulfils \eqref{eq:admissible_lambda}.

\myendproof

\section{Proof of Lemma \ref{lem:mean}}
\label{sec:lem_mean_proof}

It is easy to observe that
\begin{align*}
    \ttm
    &
    = n \; \E \|\widehat \Sigma - \Sigma\|_{\Fr}^2
    = n \Tr\left[ (\widehat \Sigma - \Sigma)^2 \right]
    \\&
    = \frac1n \sum\limits_{i = 1}^n \sum\limits_{j = 1}^n \E \Tr\left[ (\bX_i \bX_i^\top - \Sigma) (\bX_j \bX_j^\top - \Sigma) \right]
    \\&
    = \frac1n \sum\limits_{i = 1}^n \E \Tr\left[ (\bX_i \bX_i^\top - \Sigma)^2 \right]
    = \E \Tr\left[ (\bX \bX^\top - \Sigma)^2 \right].
\end{align*}
Let us represent the expectation of $\Tr\left[ (\bX \bX^\top - \Sigma)^2 \right]$ in the following way:
\begin{align*}
    \E \Tr\left[ (\bX \bX^\top - \Sigma)^2 \right]
    &
    = \E \Tr\left[ \bX \bX^\top \bX \bX^\top - \Sigma \bX \bX^\top - \bX \bX^\top \Sigma + \Sigma^2 \right]
    = \E \|\bX\|^4 - \Tr(\Sigma^2)
    \\&
    = \E \left( \|\bX\|^2 - \Tr(\Sigma) \right)^2 + \left( \Tr(\Sigma) \right)^2 - \Tr(\Sigma^2).
\end{align*}
Obviously, the first term in the right-hand side is non-negative. This yields that
\[
    \ttm \geq \left( \Tr(\Sigma) \right)^2 - \Tr(\Sigma^2).
\]
To prove the upper bound, note that, due to Assumption \ref{as:fourth_moment} and the Cauchy-Schwarz inequality, we have
\[
    \E \left( \|\bX\|^2 - \Tr(\Sigma) \right)^2
    \leq \sqrt{\E \left( \|\bX\|^2 - \Tr(\Sigma) \right)^4}
    \leq \alpha \|\Sigma\|_{\Fr}^2 = \alpha \; \Tr(\Sigma^2).
\]
Thus, we obtain that
\[
    \ttm \leq \left( \Tr(\Sigma) \right)^2 + (\alpha - 1) \Tr(\Sigma^2).
\]
\myendproof

\section{Proof of Lemma \ref{lem:variance}}
\label{sec:lem_variance_proof}

\noindent
\textbf{The bound on $\ttv^2$.}\quad
We start with the upper bound on $\ttv^2$.
Let us introduce $\bfeta = \rmvec(\bX\bX^\top - \Sigma)$ for brevity and let $\widetilde \bX$, $\widetilde\bfeta$ be independent copies of $\bX$ and $\bfeta$, respectively. Then it holds that
\begin{align*}
    \ttv^2
    &
    = \left\| \E \bfeta \bfeta^\top \right\|_{\Fr}^2
    = \Tr\left[ \E \left( \bfeta \bfeta^\top \right) \E \left( \widetilde \bfeta \widetilde \bfeta^\top \right) \right]
    = \left[ \E \left( \bfeta^\top \widetilde \bfeta \right)^2 \right]
    \\&
    = \E \left( \Tr \left[ (\bX\bX^\top - \Sigma) (\widetilde \bX \widetilde \bX^\top - \Sigma) \right] \right)^2.
\end{align*}
Let $\widetilde \E$ stand for the expectation with respect to $\widetilde \bX$ (conditioned on $\bX$). Due to Assumption \ref{as:fourth_moment} and the Cauchy-Schwarz inequality, we have
\begin{align}
    \label{eq:expectation_trace_bound}
    \E \left( \Tr \left[ (\bX\bX^\top - \Sigma) (\widetilde \bX \widetilde \bX^\top - \Sigma) \right] \right)^2 
    &\notag
    = \E \; \widetilde \E \left( \widetilde\bX^\top (\bX \bX^\top - \Sigma) \widetilde \bX - \Tr\left[ \Sigma (\bX \bX^\top - \Sigma) \right] \right)^2
    \\&
    \leq \E \sqrt{ \widetilde \E \left( \widetilde\bX^\top (\bX \bX^\top - \Sigma) \widetilde \bX - \Tr\left[ \Sigma (\bX \bX^\top - \Sigma) \right] \right)^4 }
    \\&\notag
    \leq \alpha \; \E \left\|\Sigma^{1/2} (\bX \bX^\top - \Sigma) \Sigma^{1/2} \right\|_{\Fr}^2. 
\end{align}
Before we prove a bound on the expectation in the right-hand side, we use the linearity of the trace to rewrite it in the following form:
\begin{align*}
    \E \left\|\Sigma^{1/2} (\bX \bX^\top - \Sigma) \Sigma^{1/2} \right\|_{\Fr}^2
    &
    = \E \Tr\left[ (\bX \bX^\top - \Sigma) \Sigma (\bX \bX^\top - \Sigma) \Sigma \right]
    \\&
    = \E (\bX^\top \Sigma \bX)^2 - 2 \E \Tr\left( \bX \bX^\top \Sigma^3 \right) + \Tr(\Sigma^4)
    \\&
    = \E (\bX^\top \Sigma \bX)^2 - \Tr(\Sigma^4)
\end{align*}
Now, we can use Assumption \ref{as:fourth_moment} again:
\begin{align}
    \label{eq:expectation_fr_norm_bound}
    \E \left\|\Sigma^{1/2} (\bX \bX^\top - \Sigma) \Sigma^{1/2} \right\|_{\Fr}^2
    &\notag
    = \E \left(\bX^\top \Sigma \bX  - \E \bX^\top \Sigma \bX \right)^2 + \left( \E \bX^\top \Sigma \bX \right)^2 - \Tr(\Sigma^4)
    \\&
    \leq \alpha \left\|\Sigma^2 \right\|_{\Fr}^2 + \left( \E \bX^\top \Sigma \bX \right)^2 - \Tr(\Sigma^4)
    \\&\notag
    = \left( \alpha - 1 \right) \Tr(\Sigma^4) + \left( \Tr(\Sigma^2) \right)^2.
\end{align}
The inequalities \eqref{eq:expectation_trace_bound} and \eqref{eq:expectation_fr_norm_bound} yield that
\[
    \ttv^2 \leq \alpha \left( \Tr(\Sigma^2) \right)^2 + (\alpha^2 - \alpha) \Tr(\Sigma^4).
\]

\bigskip

\noindent
\textbf{The bound on $\kappa$.}\quad It only remains to derive the upper bound on $\kappa$. Let us fix any $U \in \R^{d\times d}$ satisfying the equality $\|U\|_{\Fr} = 1$ and show that
\[
    \E \left[ \bX^\top U \bX - \Tr(U^\top \Sigma) \right]^2 \leq \alpha \, \|\Sigma\|^2.
\]
Applying the Cauchy-Schwarz inequality and using Assumption \ref{as:fourth_moment}, we obtain that
\begin{align*}
    \E \left[ \bX^\top U^\top \bX - \Tr(U^\top \Sigma) \right]^2
    \leq \sqrt{\E \big( \bX^\top U^\top \bX - \Tr(U^\top \Sigma) \big)^4}
    \leq \alpha \left\| \Sigma^{1/2} U \Sigma^{1/2} \right\|_{\Fr}^2.
\end{align*}
The bound on $\kappa$ follows from the property \eqref{eq:kronecker_trace} of the Kronecker product. Indeed, we have
\[
    \left\| \Sigma^{1/2} U \Sigma^{1/2} \right\|_{\Fr}^2
    = \Tr\left( U^\top \Sigma U \Sigma \right)
    = \rmvec(U)^\top (\Sigma \otimes \Sigma) \rmvec(U)
    \leq \left\| \Sigma \otimes \Sigma \right\| \|\rmvec(U)\|^2.
\]
Since $\rmvec(U)$ is just a reshaping of $U$, we have
$\|\rmvec(U)\| = \|U\|_{\Fr} \leq 1$.
Besides, the operator norm of the Kronecker product of symmetric matrices is equal to the product of operator norms (see eq. \eqref{eq:kronecker_eigenvalues}):
\[
    \left\| \Sigma \otimes \Sigma \right\| = \|\Sigma\|^2.
\]
Hence, it holds that $\kappa \leq \alpha \|\Sigma\|^2$.

\myendproof

\section{Proof of Lemma \ref{lem:kth_moment}}
\label{sec:lem_kth_moment_proof}

First, note that
\[
    \left| \sfP_{\i U} \left( \bxi^\top V \bxi - \sfP_{\i U} \bxi^\top V \bxi \right)^k \right|
    = \frac{\left| \E \left( \bxi^\top V \bxi - \sfP_{\i U} \bxi^\top V \bxi \right)^k e^{\i \bxi^\top U \bxi} \right|}{\left| \E e^{\i \bxi^\top U \bxi} \right|}
    \leq \beta \; \E \left| \bxi^\top V \bxi - \sfP_{\i U} \bxi^\top V \bxi \right|^k.
\]
Using the inequality $(a + b)^k \leq 2^{k - 1} a^k + 2^{k - 1} b^k$, which holds for any $a, b > 0$ and any $k \geq 1$, we obtain that
\begin{equation}
    \label{eq:jensen}
    2^{- k + 1} \; \E \left| \bxi^\top V \bxi - \sfP_{\i U} \bxi^\top V \bxi \right|^k
    \leq \E \left| \bxi^\top V \bxi - \E \bxi^\top V \bxi \right|^k + \left| \E \bxi^\top V \bxi - \sfP_{\i U} \bxi^\top V \bxi \right|^k.
\end{equation}
In view of Assumption \ref{as:fourth_moment}, the first term in the right hand side does not exceed
\[
    \E \left| \bxi^\top V \bxi - \E \bxi^\top V \bxi \right|^k
    \leq \left( \E \left| \bxi^\top V \bxi - \E \bxi^\top V \bxi \right|^4 \right)^{k / 4}
    \leq \alpha^{k / 2} \|V\|_{\Fr}^k.
\]
On the other hand, due to the Cauchy-Schwarz inequality, the second term in the right-hand side of \eqref{eq:jensen} is not greater than
\begin{align*}
    \left| \E \bxi^\top V \bxi - \sfP_{\i U} \bxi^\top V \bxi \right|^k
    &
    = \frac{\left| \E \bxi^\top V \bxi \; \E e^{\i \bxi^\top U \bxi} - \E \, \bxi^\top V \bxi \, e^{\i \bxi^\top U \bxi} \right|^k}{\left| \E e^{\i \bxi^\top U \bxi} \right|^k}
    \\&
    = \frac{\left| \E (\bxi^\top V \bxi - \E \, \bxi^\top V \bxi) \left( \E e^{\i \bxi^\top U \bxi} - e^{\i \bxi^\top U \bxi} \right) \right|^k}{\left| \E e^{\i \bxi^\top U \bxi} \right|^k}
    \\&
    \leq \beta^k \left( \E \left( \bxi^\top V \bxi - \E \, \bxi^\top V \bxi \right)^2 \cdot \E \left| e^{\i \bxi^\top U \bxi} - \E e^{\i \bxi^\top U \bxi} \right|^2 \right)^{k / 2}.
\end{align*}
Taking into account that
\[
    \E \left| e^{\i \bxi^\top U \bxi} - \E e^{\i \bxi^\top U \bxi} \right|^2
    = \E \left(e^{\i \bxi^\top U \bxi} - \E e^{\i \bxi^\top U \bxi} \right) \left( e^{-\i \bxi^\top U \bxi} - \E e^{-\i \bxi^\top U \bxi} \right)
    = 1 - \left| \E e^{-\i \bxi^\top U \bxi} \right|^2
    \leq 1,
\]
we conclude that
\[
    \left| \E \bxi^\top V \bxi - \sfP_{\i U} \bxi^\top V \bxi \right|^k
    \leq \beta^k \alpha^{k/2} \|V\|_{\Fr}^k.
\]
Hence, it holds that
\[
    \left| \E \bxi^\top V \bxi - \sfP_{\i U} \bxi^\top V \bxi \right|^k
    \leq 2^{k - 1} \beta \left( 1 + \beta^k \right) \alpha^{k/2} \|V\|_{\Fr}^k. 
\]
\myendproof

\section{Proof of Lemma \ref{lem:exp_moment_bound}}
\label{sec:lem_exp_moment_bound_proof}

Let us rearrange the expression
\[
    -\frac{\lambda}2 \; \E_{\bX} \left( \bX^\top \Gamma \bX - \Tr(\Gamma \Sigma) \right)^2 + \frac{5 \alpha^2 \lambda^2 \rho^2}{n} \left\|\Sigma^{1/2} \Gamma \Sigma^{1/2} \right\|_{\Fr}^2
\]
in the following form:
\[
    -\frac{\lambda}2 \; \rmvec(\Gamma)^\top \left( \E_{\bX} (\bX \bX^\top  - \Sigma) \otimes (\bX \bX^\top  - \Sigma) \right) \rmvec(\Gamma) + \frac{5 \alpha^2 \lambda^2 \rho^2}{n} \; \rmvec(\Gamma)^\top \left( \Sigma \otimes \Sigma \right) \rmvec(\Gamma).
\]
Then it holds that
\begin{align}
    \label{eq:gamma_exp_moment}
    &\notag
    \E_\Gamma \exp\left\{ -\frac{\lambda}2 \; \E_{\bX} \left( \bX^\top \Gamma \bX - \Tr(\Gamma \Sigma) \right)^2 + \frac{5 \alpha^2 \lambda^2 \rho^2}{n} \left\|\Sigma^{1/2} \Gamma \Sigma^{1/2} \right\|_{\Fr}^2 \right\}
    \\&
    = \det\left( I_{d^2} + \lambda \; \E_{\bX} \left[ (\bX \bX^\top  - \Sigma) \otimes (\bX \bX^\top  - \Sigma) \right] - \frac{10 \alpha^2 \lambda^2 \rho^2}{n} \; \left( \Sigma \otimes \Sigma \right) \right)^{-1/2}.
\end{align}
Similarly to Lemma \ref{lem:restricted_moment}, we use the inequality
\[
        -\log \det(I_{d^2} - A) \leq \Tr(A) + 2 \|A\|_{\Fr}^2,
\]
which holds for any matrix $A \in \R^{d^2 \times d^2}$, whose maximal eigenvalue does not exceed $1/2$.
Due to the conditions of the lemma, we have
\[
    \frac{10 \alpha^2 \lambda^2 \rho^2 \|\Sigma\|^2}{n} \leq \frac12.
\]
Then it holds that
\begin{align*}
    &
    -\frac12 \log \det\left( I_{d^2} + \lambda \; \E_{\bX} \left[ (\bX \bX^\top  - \Sigma) \otimes (\bX \bX^\top  - \Sigma) \right] - \frac{10 \alpha^2 \lambda^2 \rho^2}{n} \; \left( \Sigma \otimes \Sigma \right) \right)
    \\&
    \leq -\frac{\lambda}2 \E_{\bX} \Tr \left[ (\bX \bX^\top  - \Sigma) \otimes (\bX \bX^\top  - \Sigma) \right] + \frac{5 \alpha^2 \lambda^2 \rho^2}{n} \;  \Tr\left( \Sigma \otimes \Sigma \right)
    \\&\quad
    + \left\| \lambda \; \E_{\bX} \left[ (\bX \bX^\top  - \Sigma) \otimes (\bX \bX^\top  - \Sigma) \right] - \frac{10 \alpha^2 \lambda^2 \rho^2}{n} \; \left( \Sigma \otimes \Sigma \right) \right\|_{\Fr}^2.
\end{align*}
The last term in the right-hand side does not exceed
\[
    2 \left\| \lambda \; \E_{\bX} \left[ (\bX \bX^\top  - \Sigma) \otimes (\bX \bX^\top  - \Sigma) \right] \right\|_{\Fr}^2
    + 2 \left\| \frac{10 \alpha^2 \lambda^2 \rho^2}{n} \; \left( \Sigma \otimes \Sigma \right) \right\|_{\Fr}^2
    = 2 \lambda^2 \ttv^2 + \frac{200 \alpha^4 \lambda^4 \rho^4 \|\Sigma\|_{\Fr}^4}{n^2}.
\]
Taking into account that
\[
    \E_{\bX} \Tr \left[ (\bX \bX^\top  - \Sigma) \otimes (\bX \bX^\top  - \Sigma) \right]
    = \E_{\bX} \left[ \Tr(\bX \bX^\top  - \Sigma) \right]^2 = \ttm,
\]
we obtain the inequality
\begin{align*}
    &
    -\frac12 \log \det\left( I_{d^2} + \lambda \; \E_{\bX} \left[ (\bX \bX^\top  - \Sigma) \otimes (\bX \bX^\top  - \Sigma) \right] - \frac{10 \alpha^2 \lambda^2 \rho^2}{n} \; \left( \Sigma \otimes \Sigma \right) \right)
    \\&
    \leq -\frac{\lambda \ttm}2 + 2 \lambda^2 \ttv^2 + \frac{5 \alpha^2 \lambda^2 \rho^2}{n} \left( \Tr(\Sigma) \right)^2 + \frac{200 \alpha^4 \lambda^4 \rho^4 \|\Sigma\|_{\Fr}^4}{n^2},
\end{align*}
which yields the assertion of Lemma \ref{lem:exp_moment_bound} in view of \eqref{eq:gamma_exp_moment}.

\myendproof

\section{Proof of Lemma \ref{lem:rho_inequalities}}
\label{sec:lem_rho_inequalities_proof}

It is easy to observe that, according to the definition of $\rho$,
\[
    \frac{2 \rho^2}{\|\Sigma\|_{\Fr}^2}
    \geq 4 \left( \frac{\ttr(\Sigma)^2}{\alpha} + \ttr(\Sigma^2) \right)
    > \frac{\ttr(\Sigma)^2}{\alpha} + \ttr(\Sigma^2),
\]
so it only remains to check the inequality
\[
    \frac{(\rho^2 - \left[ \Tr(\Sigma) \right]^2 )^2}{\|\Sigma\|_{\Fr}^4 + \|\Sigma\|^2 (\rho^2 - \left[ \Tr(\Sigma) \right]^2 )} \geq \frac{\ttr(\Sigma)^2}{\alpha} + \ttr(\Sigma^2).
\]
Let us denote $\tta = \ttr(\Sigma)^2 / \alpha + \ttr(\Sigma^2)$ for brevity and consider the inequality
\[
    \frac{\ttx^2}{\|\Sigma\|_{\Fr}^4 + \|\Sigma\|^2 \ttx} \geq \tta,
    \quad \ttx \geq 0,
\]
which is equivalent to
\[
    \ttx^2 - \tta \|\Sigma\|^2 \ttx - \tta \|\Sigma\|_{\Fr}^4 \geq 0,
    \quad \ttx \geq 0.
\]
Its solution is given by
\[
    \ttx \geq \frac{\tta \|\Sigma\|^2}2 + \frac12 \sqrt{\frac{\tta^2 \|\Sigma\|^4}4 + 4 \tta \|\Sigma\|_{\Fr}^4},
\]
so we have to check that
\[
    \rho^2 - \left[\Tr(\Sigma)\right]^2 \geq \frac{\tta \|\Sigma\|^2}2 + \frac12 \sqrt{\frac{\tta^2 \|\Sigma\|^4}4 + 4 \tta \|\Sigma\|_{\Fr}^4}.
\]
We will show even a tighter bound, namely,
\[
    \rho^2 - \left[\Tr(\Sigma)\right]^2
    \geq \tta \|\Sigma\|^2 + \|\Sigma\|_{\Fr}^2 \sqrt{\tta}.
\]
By the definition of $\rho$ (see eq. \eqref{eq:rho}), it holds that
\[
    \frac{\rho^2 - \left[\Tr(\Sigma)\right]^2}{\|\Sigma\|^2}
    = 2 \ttr(\Sigma^2) \left( \frac{\ttr(\Sigma)^2}{\alpha} + \ttr(\Sigma^2) \right).
\]
Since the effective rank of any non-zero symmetric matrix is at least $1$, we have
\[
    \ttr(\Sigma^2) \geq 1
    \quad \text{and} \quad
    \frac{\ttr(\Sigma)^2}{\alpha} + \ttr(\Sigma^2) \geq \left( \frac{\ttr(\Sigma)^2}{\alpha} + \ttr(\Sigma^2) \right)^{1/2}.
\]
This implies that
\begin{align*}
    \frac{\rho^2 - \left[\Tr(\Sigma)\right]^2}{\|\Sigma\|^2}
    &
    = 2 \ttr(\Sigma^2) \left( \frac{\ttr(\Sigma)^2}{\alpha} + \ttr(\Sigma^2) \right)
    \\&
    \geq \left( \frac{\ttr(\Sigma)^2}{\alpha} + \ttr(\Sigma^2) \right) + \ttr(\Sigma^2) \left( \frac{\ttr(\Sigma)^2}{\alpha} + \ttr(\Sigma^2) \right)^{1/2}
    \\&
    = \tta + \ttr(\Sigma^2) \sqrt{\tta}
    \\&
    = \frac{\tta \|\Sigma\|^2 + \|\Sigma\|_{\Fr}^2 \sqrt{\tta}}{\|\Sigma\|^2}.
\end{align*}
\myendproof

\section{Proof of Lemma \ref{lem:squared_fr_norm_variance_bound}}
\label{sec:lem_squared_fr_norm_variance_bound_proof}

We split the proof into several steps for the ease of exposure.

\medskip

\noindent
\textbf{Step 1: the Efron-Stein inequality.}
\quad
The key idea of the proof is to bound the variance, applying the Efron-Stein inequality to $\|\widehat\Sigma - \Sigma\|_{\Fr}^2$. Let $\widetilde\bX_1, \dots, \widetilde\bX_n$ be i.i.d. copies of the random vector $\bX$, which are independent of $\bX_1, \dots, \bX_n$. For any $i \in \{1, \dots, n\}$, let us introduce
\[
    \widehat \Sigma_{-i} = \frac1n \sum\limits_{j \neq i} \bX_j \bX_j^\top
    \quad \text{and} \quad
    \widehat \Sigma_i = \frac1n \widetilde\bX_i \widetilde\bX_i^\top + \widehat \Sigma_{-i}.
\]
According to the Efron-Stein inequality (see, for instance, \cite[Theorem 3.1]{boucheron13}), it holds that
\[
    \Var\left( \|\widehat\Sigma - \Sigma\|_{\Fr}^2 \right)
    \leq \frac12 \sum\limits_{i = 1}^n \E \left( \|\widehat\Sigma - \Sigma\|_{\Fr}^2 - \|\widehat\Sigma_i - \Sigma\|_{\Fr}^2 \right)^2.
\]
For any $i \in \{1, \dots, n\}$, we have
\begin{align*}
    \left\|\widehat\Sigma - \Sigma \right\|_{\Fr}^2 - \left\|\widehat\Sigma_i - \Sigma \right\|_{\Fr}^2
    &
    = \left\|\frac1n \bX_i \bX_i^\top + \widehat\Sigma_{-i} - \Sigma \right\|_{\Fr}^2 - \left\|\frac1n \widetilde\bX_i \widetilde\bX_i^\top + \widehat\Sigma_{-i} - \Sigma \right\|_{\Fr}^2
    \\&
    = \left\|\frac1n \big(\bX_i \bX_i^\top - \Sigma\big) + \widehat\Sigma_{-i} - \frac{n-1}n \Sigma \right\|_{\Fr}^2 - \left\|\frac1n \big( \widetilde\bX_i \widetilde\bX_i^\top - \Sigma \big) + \widehat\Sigma_{-i} - \frac{n-1}n \Sigma \right\|_{\Fr}^2
    \\&
    = \frac1{n^2} \left\|\bX_i \bX_i^\top - \Sigma\right\|_{\Fr}^2 - \frac1{n^2} \left\|\widetilde\bX_i \widetilde\bX_i^\top - \Sigma\right\|_{\Fr}^2 
    \\&\quad
    + \frac2n \Tr\left[ (\bX_i \bX_i^\top - \widetilde\bX_i \widetilde\bX_i^\top) \left(\widehat\Sigma_{-i} - \frac{n-1}n \Sigma \right) \right].
\end{align*}
Then, applying the Cauchy-Schwarz inequality, we obtain that
\begin{align}
    \label{eq:squared_fr_norm_variance_cauchy-schwarz}
    \Var\left( \|\widehat\Sigma - \Sigma\|_{\Fr}^2 \right)
    &\notag
    \leq \frac1{n^4} \sum\limits_{i = 1}^n \E \left( \left\|\bX_i \bX_i^\top - \Sigma\right\|_{\Fr}^2 - \left\|\widetilde\bX_i \widetilde\bX_i^\top - \Sigma\right\|_{\Fr}^2 \right)^2
    \\&\quad\notag
    + \frac4{n^2} \sum\limits_{i = 1}^n \E \left( \Tr\left[ (\bX_i \bX_i^\top - \widetilde\bX_i \widetilde\bX_i^\top) \left(\widehat\Sigma_{-i} - \frac{n-1}n \Sigma \right) \right] \right)^2
    \\&
    = \frac1{n^3} \; \E \left( \left\|\bX_1 \bX_1^\top - \Sigma\right\|_{\Fr}^2 - \left\|\widetilde\bX_1 \widetilde\bX_1^\top - \Sigma\right\|_{\Fr}^2 \right)^2
    \\&\quad\notag
    + \frac4{n} \; \E \left( \Tr\left[ (\bX_1 \bX_1^\top - \widetilde\bX_1 \widetilde\bX_1^\top) \left(\widehat\Sigma_{-1} - \frac{n-1}n \Sigma \right) \right] \right)^2.
\end{align}
It remains to bound the two terms in the right-hand side.

\medskip

\noindent
\textbf{Step 2: a bound on $\E \left( \|\bX_1 \bX_1^\top - \Sigma\|_{\Fr}^2 - \|\widetilde\bX_1 \widetilde\bX_1^\top - \Sigma\|_{\Fr}^2 \right)^2$.}
\quad
Since the random variables $\|\bX_1 \bX_1^\top - \Sigma\|_{\Fr}^2$ and $\|\widetilde\bX_1 \widetilde\bX_1^\top - \Sigma\|_{\Fr}^2$ are independent and identically distributed, it holds that
\begin{equation}
    \label{eq:squared_diff_squared_fr_norms}
    \E \left( \left\|\bX_1 \bX_1^\top - \Sigma\right\|_{\Fr}^2 - \left\|\widetilde\bX_1 \widetilde\bX_1^\top - \Sigma\right\|_{\Fr}^2 \right)^2
    = 2 \; \Var\left( \left\|\bX \bX^\top - \Sigma\right\|_{\Fr}^2 \right).
\end{equation}
We use the following lemma to bound the variance of $\|\bX \bX^\top - \Sigma\|_{\Fr}^2$.
\begin{Lem}
    \label{lem:zeta_squared_fr_norm_variance}
    Let $\bzeta \in \R^d$ be a centered random vector with a nondegenerate covariance matrix $\Omega$. Suppose that $\bzeta$ satisfies Assumption \ref{as:fourth_moment}, that is,
    \[
        \E \left( \bzeta^\top \Omega^{-1/2} V \Omega^{-1/2} \bzeta - \Tr(V) \right)^4 \leq \alpha^2 \|V\|_{\Fr}^4
        \quad \text{for all $V \in \R^{d \times d}$.}
    \]
    Then it holds that
    \[
        \Var\left( \left\| \bzeta \bzeta^\top - \Omega \right\|_{\Fr}^2 \right)
        \leq 4 \alpha^2 \Tr(\Omega^2)^2 + 24 \alpha \Tr(\Omega)^2 \Tr(\Omega^2)
    \]
    and
    \[
        \E \left\| \bzeta \bzeta^\top - \Omega \right\|_{\Fr}^2
        \leq \Tr(\Omega)^2 + (\alpha - 1) \Tr(\Omega^2).
    \]
\end{Lem}
We postpone the proof of Lemma \ref{lem:zeta_squared_fr_norm_variance} to Appendix \ref{sec:lem_zeta_squared_fr_norm_variance_proof}  and proceed with the following inequality:
\begin{equation}
    \label{eq:x_squared_fr_norm_variance}
    \Var\left( \left\|\bX \bX^\top - \Sigma\right\|_{\Fr}^2 \right)
    \leq 4 \alpha^2 \Tr(\Sigma^2)^2 + 24 \alpha \Tr(\Sigma)^2 \Tr(\Sigma^2).
\end{equation}
It remains to bound the second term in the right-hand side of \eqref{eq:squared_fr_norm_variance_cauchy-schwarz}.

\medskip

\noindent
\textbf{Step 3: a bound on $\E \left( \Tr\left[ (\bX_1 \bX_1^\top - \widetilde\bX_1 \widetilde\bX_1^\top) \left(\widehat\Sigma_{-1} - \frac{n-1}n \Sigma \right) \right] \right)^2$.}
\quad
Note that
\begin{align*}
    &
    \E \left( \Tr\left[ (\bX_1 \bX_1^\top - \widetilde\bX_1 \widetilde\bX_1^\top) \left(\widehat\Sigma_{-1} - \frac{n-1}n \Sigma \right) \right] \right)^2
    \\&
    = \E \, \E \left[ \left( \Tr\left[ (\bX_1 \bX_1^\top - \widetilde\bX_1 \widetilde\bX_1^\top) \left(\widehat\Sigma_{-1} - \frac{n-1}n \Sigma \right) \right] \right)^2 \;\Big\vert\; \bX_2, \dots, \bX_n \right]
    \\&
    = 2 \, \E \, \Var\left( \Tr\left[ \bX_1 \bX_1^\top  \left(\widehat\Sigma_{-1} - \frac{n-1}n \Sigma \right) \right] \;\Big\vert\; \bX_2, \dots, \bX_n \right)
    \\&
    = 2 \, \E \, \Var\left( \bX_1^\top  \left(\widehat\Sigma_{-1} - \frac{n-1}n \Sigma \right) \bX_1 \;\Big\vert\; \bX_2, \dots, \bX_n \right).
\end{align*}
According to Assumption \ref{as:fourth_moment}, the expression in the right-hand side does not exceed
\begin{align*}
    &
    2 \, \E \, \Var\left( \bX_1^\top  \left(\widehat\Sigma_{-1} - \frac{n-1}n \Sigma \right) \bX_1 \;\Big\vert\; \bX_2, \dots, \bX_n \right)
    \\&
    \leq 2 \alpha \, \E \left\| \Sigma^{1/2} \left( \widehat\Sigma_{-1} - \frac{n-1}n \Sigma \right) \Sigma^{1/2} \right\|_{\Fr}^2.
\end{align*}
Representing $\Sigma^{1/2} \big( \widehat\Sigma_{-1} - \frac{n-1}n \Sigma \big) \Sigma^{1/2}$ as a sum of independent centered random matrices, we obtain that
\begin{align*}
    \E \left\| \Sigma^{1/2} \big( \widehat\Sigma_{-1} - \frac{n-1}n \Sigma \big) \Sigma^{1/2} \right\|_{\Fr}^2
    &
    = \frac1{n^2} \sum\limits_{i = 2}^n \sum\limits_{j = 2}^n \E \; \Tr\left( \Sigma (\bX_i \bX_i^\top - \Sigma) \Sigma (\bX_j \bX_j^\top - \Sigma) \right) 
    \\&
    = \frac1{n^2} \sum\limits_{i = 2}^n \E \; \Tr\left( \Sigma (\bX_i \bX_i^\top - \Sigma) \Sigma (\bX_i \bX_i^\top - \Sigma) \right)
    \\&
    = \frac{n - 1}{n^2} \; \E \left\| \Sigma^{1/2} (\bX \bX^\top - \Sigma) \Sigma^{1/2} \right\|_{\Fr}^2.
\end{align*}
Thus, it holds that
\[
    \E \left( \Tr\left[ (\bX_1 \bX_1^\top - \widetilde\bX_1 \widetilde\bX_1^\top) \left(\widehat\Sigma_{-1} - \frac{n-1}n \Sigma \right) \right] \right)^2
    \leq \frac{2}{n} \; \E \left\| \Sigma^{1/2} (\bX \bX^\top - \Sigma) \Sigma^{1/2} \right\|_{\Fr}^2.
\]
Applying Lemma \ref{lem:zeta_squared_fr_norm_variance} to the random vector $\Sigma^{1/2} \bX$, we obtain that
\begin{equation}
    \label{eq:expectation_squared_trace bound}
    \E \left( \Tr\left[ (\bX_1 \bX_1^\top - \widetilde\bX_1 \widetilde\bX_1^\top) \left(\widehat\Sigma_{-1} - \frac{n-1}n \Sigma \right) \right] \right)^2
    \leq \frac{2\Tr(\Sigma^2)^2 + 2(\alpha - 1)\Tr(\Sigma^4)}n.
\end{equation}

\medskip

\noindent
\textbf{Step 4: final bound.}
\quad
Due to the inequalities \eqref{eq:squared_fr_norm_variance_cauchy-schwarz}, \eqref{eq:squared_diff_squared_fr_norms}, \eqref{eq:x_squared_fr_norm_variance}, and \eqref{eq:expectation_squared_trace bound} the variance of $\|\widehat\Sigma - \Sigma\|_{\Fr}^2$ is not greater than
\[
    \Var\left( \|\widehat\Sigma - \Sigma\|_{\Fr}^2 \right)
    \leq \frac{8 \Tr(\Sigma^2)^2 + 8\alpha \Tr(\Sigma^4)}{n^2}
    + \frac{8 \alpha^2 \Tr(\Sigma^2)^2 + 48 \alpha \Tr(\Sigma)^2 \Tr(\Sigma^2)}{n^3}.
\]
The proof is finished.

\myendproof

\section{Proof of Lemma \ref{lem:zeta_squared_fr_norm_variance}}
\label{sec:lem_zeta_squared_fr_norm_variance_proof}

Rewriting $\|\bzeta \bzeta^\top - \Omega\|_{\Fr}^2$ in the form
\[
    \left\|\bzeta \bzeta^\top - \Omega\right\|_{\Fr}^2
    = \Tr\left[(\bzeta \bzeta^\top - \Omega)(\bzeta \bzeta^\top - \Omega)\right]
    = \|\bzeta\|^4 - 2 \bzeta^\top \Omega \bzeta + \Tr(\Omega^2),
\]
we note that
\begin{equation}
    \label{eq:squared_fr_norm_variance_bound}
    \Var\left( \left\|\bzeta \bzeta^\top - \Omega\right\|_{\Fr}^2 \right)
    = \Var\left( \|\bzeta\|^4 - 2 \bzeta^\top \Omega \bzeta \right)
    \leq 2 \Var\left( \|\bzeta\|^4 \right) + 8 \Var\left( \bzeta^\top \Omega \bzeta \right).
\end{equation}
Here we used the Cauchy-Schwarz inequality to bound the variance of $\|\bzeta\|^4 - 2 \bzeta^\top \Omega \bzeta$. Due to Assumption \ref{as:fourth_moment}, it holds that
\begin{align}
    \label{eq:quadratic_form_variance_bound}
    \Var\left( \bzeta^\top \Omega \bzeta \right)
    &\notag
    = \E \left( \bzeta^\top \Omega \bzeta - \E \bzeta^\top \Omega \bzeta \right)^2
    = \E \left( \bzeta^\top \Omega \bzeta - \Tr(\Omega^2) \right)^2
    \\&
    \leq \sqrt{\E \left( \bzeta^\top \Omega \bzeta - \Tr(\Omega^2) \right)^4}
    \leq \alpha \left\|\Omega^2 \right\|_{\Fr}^2
    = \alpha \Tr(\Omega^4).
\end{align}
On the other hand,
\begin{align*}
    \Var\left( \|\bzeta\|^4 \right)
    &
    = \Var \left( \big(\|\bzeta\|^2 - \Tr(\Omega) + \Tr(\Omega) \big)^2 \right)
    \\&
    = \Var \left( \big(\|\bzeta\|^2 - \Tr(\Omega)\big)^2 + 2 \Tr(\Omega) \big(\|\bzeta\|^2 - \Tr(\Omega)\big) + \Tr(\Omega)^2 \right)
    \\&
    = \Var \left( \big(\|\bzeta\|^2 - \Tr(\Omega)\big)^2 + 2 \Tr(\Omega) \big(\|\bzeta\|^2 - \Tr(\Omega)\big) \right).
\end{align*}
The Cauchy-Schwarz inequality implies that
\begin{align*}
    \Var\left( \|\bzeta\|^4 \right)
    &
    = \Var\left( \big(\|\bzeta\|^2 - \Tr(\Omega)\big)^2 + 2 \Tr(\Omega) \big(\|\bzeta\|^2 - \Tr(\Omega)\big) \right)
    \\&
    \leq 2 \Var\left( \big(\|\bzeta\|^2 - \Tr(\Omega)\big)^2 \right) + 8 \Tr(\Omega)^2 \, \Var\left( \|\bzeta\|^2 - \Tr(\Omega) \right)
    \\&
    \leq 2 \E \left( \|\bzeta\|^2 - \Tr(\Omega) \right)^4 + 8 \Tr(\Omega)^2 \, \E \left( \|\bzeta\|^2 - \Tr(\Omega) \right)^2.
\end{align*}
Thus, according to the conditions of the lemma,
\begin{align}
    \label{eq:fourth_power_variance_bound}
    \Var\left( \|\bzeta\|^4 \right)
    &\notag
    \leq 2 \E \left( \|\bzeta\|^2 - \Tr(\Omega) \right)^4 + 8 \Tr(\Omega)^2 \E \left( \|\bzeta\|^2 - \Tr(\Omega) \right)^2
    \\&
    \leq 2 \alpha^2 \left\|\Omega\right\|_\Fr^4 + 8 \alpha \Tr(\Omega)^2 \left\|\Omega\right\|_\Fr^2
    \\&\notag
    = 2 \alpha^2 \Tr(\Omega^2)^2 + 8 \alpha \Tr(\Omega)^2 \Tr(\Omega^2).
\end{align}
Hence, summing up the inequalities \eqref{eq:squared_fr_norm_variance_bound}, \eqref{eq:quadratic_form_variance_bound}, and \eqref{eq:fourth_power_variance_bound}, we obtain that
\begin{align*}
    \Var\left( \left\|\bzeta \bzeta^\top - \Omega\right\|_{\Fr}^2 \right)
    &
    \leq 4 \alpha^2 \Tr(\Omega^2)^2 + 16 \alpha \Tr(\Omega)^2 \Tr(\Omega^2) + 8 \alpha \Tr(\Omega^4)
    \\&
    \leq 4 \alpha^2 \Tr(\Omega^2)^2 + 24 \alpha \Tr(\Omega)^2 \Tr(\Omega^2),
\end{align*}
and the first part of the proof is finished. In order to get an upper bound on the expectation of $\left\|\bzeta \bzeta^\top - \Omega\right\|_{\Fr}^2$, we observe that
\begin{align*}
    \E \left\|\bzeta \bzeta^\top - \Omega\right\|_{\Fr}^2
    &
    = \E \left( \|\bzeta\|^4 - 2 \bzeta^\top \Omega \bzeta + \Tr(\Omega^2) \right)
    = \E \left( \|\bzeta\|^2 - \Tr(\Omega) + \Tr(\Omega) \right)^2 - \Tr(\Omega^2)
    \\&
    = \E \left( \|\bzeta\|^2 - \Tr(\Omega)\right)^2 + \Tr(\Omega)^2 - \Tr(\Omega^2)
    \leq \Tr(\Omega)^2 + (\alpha - 1) \Tr(\Omega^2),
\end{align*}
where the last line follows from Assumption \ref{as:fourth_moment} and the Cauchy-Schwarz inequality.

\myendproof

\end{document}